\documentclass[11pt,a4paper,reqno,oneside]{amsproc}

\usepackage{amssymb,euscript,lmodern}

\usepackage{xcolor}

\usepackage{hyperref}
\usepackage{bbm}

\allowdisplaybreaks[3]

\renewcommand{\mod}{\ \text{\rm mod}\ }

\renewcommand{\tilde}{\widetilde}
\renewcommand{\Tilde}{\widetilde}

\newcommand{\loc}{\mathrm{loc}}

\newcommand{\cA}{\mathcal{A}}
\newcommand{\cH}{\mathcal{H}}

\newcommand{\RR}{\mathbb{R}}
\newcommand{\CC}{\mathbb{C}}
\newcommand{\NN}{\mathbb{N}}
\newcommand{\cB}{\mathcal{B}}
\newcommand{\cD}{\mathcal{D}}
\newcommand{\cG}{\mathcal{G}}
\newcommand{\cE}{\mathcal{E}}
\newcommand{\cM}{\mathcal{M}}
\newcommand{\cL}{\mathcal{L}}
\newcommand{\cS}{\mathcal{S}}

\newcommand{\rmi}{\mathrm{i}}

\newcommand{\dd}{\mathrm{d}}
\newcommand{\half}{{\frac{1}{2}}}
\newcommand{\nhalf}{{\frac{N}{2}}}

\newcommand{\dsum}{\mathop{\dot{+}}}

\renewcommand{\Tilde}{\widetilde}

\DeclareMathOperator{\supp}{\mathrm{supp}}
\DeclareMathOperator{\sgn}{\mathrm{sgn}}
\DeclareMathOperator{\ran}{\mathrm{ran}}
\DeclareMathOperator{\dom}{\mathrm{dom}}
\DeclareMathOperator{\Op}{Op}

\DeclareMathOperator{\spec}{spec}
\DeclareMathOperator{\res}{res}
\DeclareMathOperator{\pv}{p.v.}

\newcommand{\spece}{\spec_\mathrm{ess}}
\newcommand{\specp}{\spec_\mathrm{p}}

\newtheorem{theorem}{\bf Theorem}
\newtheorem{proposition}[theorem]{\bf Proposition}
\newtheorem{lemma}[theorem]{\bf Lemma}
\newtheorem{corollary}[theorem]{\bf Corollary}

\theoremstyle{definition}
\newtheorem{definition}[theorem]{\bf Definition}
\newtheorem{remark}[theorem]{\bf Remark}

\numberwithin{equation}{section}
\numberwithin{theorem}{section}

\title[]{Essential spectrum of singularly perturbed Dirac operators on bounded domains}

\author[J. Behrndt]{\large \sc Jussi Behrndt}
\address{ (J.~Behrndt) Institut f\"ur Angewandte Mathematik, Technische Universit\"at Graz, Steyrergasse 30, 8010 Graz, Austria}
\email{behrndt@tugraz.at}

\author[B. Benhellal]{Badreddine Benhellal}
\address{(B. Benhellal) National Higher School of Mathematics, Scientific and Technology Hub of Sidi Abdellah, P.O. Box 75,  16093 Algiers, Algeria}
\email{badreddine.benhellal@nhsm.edu.dz}

\author[M. Holzmann]{Markus Holzmann}
\address{(M. Holzmann) Institut f\"ur Angewandte Mathematik, Technische Universit\"at Graz, Steyrergasse 30, 8010 Graz, Austria}
\email{holzmann@math.tugraz.at}

\author[K. Pankrashkin]{Konstantin Pankrashkin}
\address{(K. Pankrashkin)  Carl von Ossietzky Universit\"at Oldenburg,
	Fakult\"at V -- Mathematik und Naturwissenschaften,  Institut f\"ur Mathematik, Ammerl\"ander Heerstr. 114--118, 26129 Oldenburg, Germany}
\email{konstantin.pankrashkin@uol.de}

\begin{document}
	
	\begin{abstract} 
	In this paper self-adjoint Dirac operators on bounded domains in $\mathbb{R}^2$ and $\mathbb{R}^3$ with elliptic boundary conditions  perturbed by singular $\delta$-shell potentials supported on compact curves or surfaces are considered. While the unperturbed operator has purely discrete spectrum, it is shown that certain critical values of the interaction strengths lead to the appearence of  essential spectrum. More precisely, each connected component of the interaction support, where the coupling constants have such critical values, produces in dimension two a point  and in dimension three even a closed interval in the essential spectrum, where the size of the latter depends on the geometry and the coupling constants.
	\end{abstract}

\maketitle	

\setcounter{tocdepth}{2}

\tableofcontents

	\section{Introduction}

	Dirac operators with singular interactions supported on curves in $\RR^2$  or surfaces in $\RR^3$ 
	have received a lot of attention in the mathematical community in recent years. Such operators in $n$ dimensions,
	with $n\in\{2,3\}$, are formally given by 
	the first order matrix differential expression acting on vector functions $f$ as
	\begin{equation}\label{dintro}   
	f\mapsto -\rmi \sum_{k=1}^{n}\alpha_k\,\partial_k f + m \beta f + (\eta I_N + \tau \beta)\delta_\Sigma f,
	\end{equation}
	where $m\in\RR$ is a mass parameter and $\alpha_k,\beta$ denote the Pauli spin matrices if $n=2$ and the Dirac matrices if $n=3$, see Section~\ref{sec21}.	Furthermore, the singular interaction term $(\eta I_N + \tau \beta)\delta_\Sigma$ consists of an electrostatic interaction with strength $\eta\in\RR$, a Lorentz scalar interaction with strength $\tau\in\RR$, and 
	$\delta_\Sigma$ stands for the $\delta$-distribution supported on a compact curve or surface $\Sigma\subset\RR^n$, $n\in\{2,3\}$.  In order to rigorously interpret the above expression as a self-adjoint operator $D^{\eta,\tau}$ 
	in $L^2(\RR^n,\CC^N)$ (with $N=2$ if $n=2$ and $N=4$ if $n=3$), the formal singular interaction term is usually treated as a  transmission condition along $\Sigma$. The first most obvious mathematical problem  is to specify an appropriate operator domain ensuring the self-adjointness, which is followed by a qualitative and quantitative spectral analysis of the operator. We refer the interested reader to 
	the review paper \cite{BHSLS24} for an overview on general results and methods in the treatment of Dirac operators with singular interactions, 
	as well as for further references, and the references listed below.

	In the following we will focus on a particularly interesting spectral peculiarity of Dirac operators of the above
	type, namely, the emergence of a new portion of essential spectrum; cf. \cite{BEHT25,BH20,BHSLS24,BCF26,OV} and, in particular, \cite{BHOP20,BP24}. In the absence of interactions ($\eta=\tau=0$) and for $m\ne 0$, the unperturbed operator $D$ is just the free Dirac operator in $\RR^n$ with the 
	spectral gap $(-|m|,|m|)$, and it was seen early that adding a singular perturbation of the above type with $\eta^2-\tau^2\ne 4$ (the so-called non-critical coupling) only produces finitely many eigenvalues in the gap, while the essential spectrum remains unchanged (which corresponds to the usual intuition that compactly supported perturbations do not affect the essential spectrum) and one has the expected Sobolev regularity of functions in the operator domain; we refer to \cite{AMV14, AMV15, BEHL18, BEHL19, BB222, BPZ25, CLMT21,DES89,OP21,R21-1,R21-2,R23} for a small list of related contributions. However, the sitation changes for $\eta^2-\tau^2=4$: this parameter combination case is usually referred to as the critical couplings.
	It was first observed for $n=2$ in \cite{BHOP20} that if $\Sigma$ is a smooth closed loop, then
	   \begin{equation}\label{eff2}
		\spece D^{\eta,\tau}=\big(-\infty,-|m|\big]\cup\left\{-\frac{\tau}{\eta} m \right\} \cup\big[|m|,\infty\big),
	\end{equation}
	i.e. a new point of essential spectrum appears.
Later, it turned out that the effect of critical couplings in three dimensions is even more striking.
More precisely, if $\Sigma\subset\RR^3$ is a compact smooth surface and $\eta^2-\tau^2=4$,
then it was observed in \cite{BP24} that
\begin{equation}\label{eff3}
     \spece D^{\eta,\tau}
     =\big(-\infty,-|m|\big]\cup\left[-\frac{\tau}{\eta}m -\frac{\cA_\Sigma}{2|\eta|},-\frac{\tau}{\eta}m 
     +\dfrac{\cA_\Sigma}{2|\eta|}\right] \cup\big[|m|,\infty\big),
\end{equation}
where $\cA_\Sigma$ is a geometric parameter given by
\[
\cA_\Sigma:=\max_{x\in\Sigma}\big|\kappa_1(x)-\kappa_2(x)\big|,
\]    
with $\kappa_j(x)$ being the principal curvatures of $\Sigma$ at $x$. In particular, if $\Sigma$ is not a union of spheres, then $\cA_\Sigma>0$,
and one obtains an additional non-trivial interval in the essential spectrum. Moreover, in all dimensions a lack of Sobolev regularity in the operator domain was observed in the critical case.
\medskip

	The present paper is inspired by an enthusiastic discussion of one of the authors 
	with the world renowned physicists Yakir Aharonov in Cetraro, Italy, in June 2024,
	who became very interested in the spectral effect for Dirac operators described above in \eqref{eff2} and \eqref{eff3}. However, for a physical interpretation it is necessary, according to Aharonov, to better understand the origin of the new essential spectrum and, in particular, it should be clarified if such an effect also appears in ``finite boxes''. More precisely, assume that $\Omega\subset\RR^n$, $n\in\{2,3\}$, is a bounded domain and let
	$\Sigma\subset\Omega$ be a smooth compact hypersurface (curve for $n=2$ and surface for $n=3$).
	Impose a self-adjoint boundary condition on $\partial\Omega$ and model the singular interaction term $(\eta I_N + \tau \beta)\delta_\Sigma$ by the same transmission condition as in $\RR^n$. The question is then whether the new essential spectrum also appears for the perturbed Dirac operator at critical couplings. The present paper addresses this question and provides an affirmative answer.

	Let us describe our results and methods in more detail. Let $\Omega\subset\RR^n$ be a bounded Lipschitz domain, $n\in\{2,3\}$,
	and $\Sigma\subset\Omega$ be a compact smooth hypersurface with finitely many connected components $\Sigma_1,\dots,\Sigma_P$; cf. Section~\ref{sec31}.
	Choose $\eta,\tau\in\RR^P$, $\eta=(\eta_1,\dots,\eta_P)$, $\tau=(\tau_1,\dots,\tau_P)$, and consider
	the Dirac operator $A^{\eta,\tau}$ in $\Omega$ with some local self-adjoint boundary conditions on $\partial\Omega$ (see Section \ref{sec21}
	and, e.g., \cite{baer,baer2,BHM,BFSV,grosse,OV}) as well as the transmission conditions
\begin{gather*}
\rmi (\alpha\cdot\nu^\Sigma)(f_{p+}-f_{p-})+\dfrac{1}{2}(\eta_p I_N+\tau_p\beta)(f_{p+} + f_{p-})=0 \text{ on each }\Sigma_p, \quad p\in\{1,\dots,P\},
\end{gather*}
modeling the distributional perturbation term in \eqref{dintro};
	here $f_{p\pm}$ are one-sided traces on the two sides of $\Sigma_p$ (see Section \ref{sec31}) and $\eta_p,\tau_p\in\RR$, $p\in\{1,\dots,P\}$, are viewed as coupling parameters. Due to some technical reasons we assume $|\eta_p|\ne |\tau_p|$ for all $p$ (we comment on this case in Remark~\ref{remlast}). As an important preparation, we show that the operator $A^{\eta,\tau}$ is self-adjoint in $L^2(\Omega,\CC^N)$ on a suitable domain
	(Theorem~\ref{thm17}). This requires a preliminary study of the integral kernel of the resolvent of the unperturbed self-adjoint 
		Dirac operator in $L^2(\Omega,\CC^N)$, which is presented in Section~\ref{sec22}.
	In Section~\ref{sec4} we turn to our main point of interest, the spectral properties of $A^{\eta,\tau}$. For the non-critical case $\eta_p^2-\tau_p^2\ne 4$ we show that the functions in the operator domain
	are $H^1$-regular near $\Sigma_p$ (Corollary~\ref{cor21}). If $\eta_p^2-\tau_p^2\ne 4$ for all $p\in\{1,\dots,P\}$, then the functions
	in the operator domain are indeed $H^1$-regular near each $\Sigma_p$, however, the inclusion $\dom A^{\eta,\tau}\subset H^1(\Omega\setminus\Sigma,\CC^N)$ cannot be ensured as the boundary condition at $\partial\Omega$
	may still require lower regularity (see Corollary~\ref{cor22} for an accurate description of the operator domain in the non-critical case). By adapting the earlier analysis for the operators in the whole space using the tools of microlocal analysis we then arrive at our main result about the 
	essential spectrum of $A^{\eta,\tau}$ (Theorem~\ref{thm37}), which reads as follows:
\begin{equation}\label{mainintro}
			\spece A^{\eta,\tau}=\bigcup_{p=1}^P J_{p},
\end{equation}
			where
	\begin{align*}
		J_{p}&:=\begin{cases}
			\emptyset, & \eta_p^2-\tau_p^2\ne 4,\\[\medskipamount]
			\left\{-\dfrac{\tau_p}{\eta_p}m\right\}, & n=2,\ \eta_p^2-\tau_p^2= 4,\\[\bigskipamount]
			-\dfrac{\tau_p}{\eta_p}m + \dfrac{1}{\eta_p}\left[ -\dfrac{a_p}{2},\dfrac{a_p}{2}\right],& n=3,\ \eta_p^2-\tau_p^2= 4,
		\end{cases}\\
		\intertext{and}		
		a_p&:=\max_{x\in\Sigma_p}\big|\kappa_1(x)-\kappa_2(x)\big|\quad\text{for}\,\,\, n=3, 
	\end{align*}
	with $\kappa_j(x)$ being the principal curvatures of $\Sigma_p$ at $x$.
	In other words, each connected component $\Sigma_p$, $p\in\{1,\dots,P\}$, carrying a critical coupling produces its own contribution $J_{p}$ to the essential spectrum.  Note that the sets $J_{p}$ depend neither on $\Omega$ nor on the boundary condition on $\partial\Omega$, and that, in contrast to the whole-space setting \eqref{eff2}-\eqref{eff3}, they are completely visible in a bounded domain, in particular also for $m=0$.   
The strategy for proving the main result \eqref{mainintro} relies on characterizing the essential spectrum of $A^{\eta,\tau}$ by reducing the spectral problem to the boundary space $L^2(\Sigma, \CC^N)$ via the abstract framework of boundary triples; cf. Section~\ref{sec31}, Lemma~\ref{lem16}, and Appendix~\ref{sec-bt}. Within this setting, the problem translates into investigating the invertibility of an operator family parameterized by the corresponding Weyl function, which explicitly incorporates certain boundary integral operators of block matrix form, where the integral kernel of the unperturbed Dirac operator in $\Omega$ mentioned before plays an important role. By analyzing the corresponding blocks using pseudodifferential methods (see Appendix~\ref{appendix_pseudos}), computing their principal symbols in Section~\ref{section_separation}, and employing decoupling and perturbation techniques in Section~\ref{Section 4.2}, we are able to explicitly determine the exact intervals $J_p$ of the essential spectrum.

\section*{Acknowledgments}

The authors are grateful to Yakir Aharonov for stimulating discussions leading to the problem of the present paper.
The research of Jussi Behrndt is funded in part by the Austrian Science Fund (FWF) 10.55776/P 33568-N.
Konstantin Pankrashkin was partially supported by the Deutsche For\-schungsgemein\-schaft (German Research Foundation, DFG), project 491606144.
This article is based upon work from COST Action 24122 mSPACE, supported by COST (European Cooperation in Science and Technology), www.cost.eu.

\section{Dirac operators on bounded Lipschitz domains}

In this preparatory section we introduce a class of self-adjoint Dirac operators on bounded Lipschitz domains and 
we construct an integral representation 
for their resolvents in Theorem~\ref{theo1}, which is of independent interest.

\subsection{Boundary conditions and self-adjointness}\label{sec21}
For $n \in \{ 2, 3 \}$ we set
\[
N:=\begin{cases}
	2, & n=2,\\
	4, & n=3,
\end{cases}
\]
and let $\alpha_1,\dots,\alpha_{n},\beta$ be pairwise anticommuting Hermitian $N\times N$ matrices such that 
$\alpha_k^2=\beta^2=I_N$, where $I_N$ is the $N\times N$ identity matrix. In order to simplify some computations
below we make an explicit choice based on the Pauli matrices
\[
\sigma_1:=\begin{pmatrix}
	0 & 1\\
	1 & 0
\end{pmatrix},
\quad
\sigma_2:=\begin{pmatrix}
	0 & -\rmi\\
	\rmi & 0
\end{pmatrix},
\quad
\sigma_3:=\begin{pmatrix}
	1 & 0\\
	0 & -1
\end{pmatrix}.
\]
For $n=2$ we choose $\alpha_k:=\sigma_k$ for $k\in\{1,2\}$ and $\beta:=\sigma_3$, while
for $n=3$ we set
\[
\alpha_k:=\begin{pmatrix}
	0 & \sigma_k\\
	\sigma_k & 0
\end{pmatrix},
\quad k\in\{1,2,3\},
\quad
\beta:=\begin{pmatrix}
	I_2 & 0\\
	0 & -I_2
\end{pmatrix}.
\]
We will use the usual abbreviations
\[
\alpha\cdot x:=\sum_{k=1}^n \alpha_k x_k,
\quad
\sigma\cdot x:=\sum_{k=1}^n \sigma_k x_k
\quad
 \text{for }x=(x_1,\dots,x_n)\in\RR^n.
\]

The Dirac operator $D_m$ with a mass $m\in\RR$ is a matrix differential operator acting on differentiable vector functions $f:\RR^n\to \CC^N$ by
\begin{equation}
	\label{eqedm}
	D_m f:=-\rmi \alpha\cdot\nabla f + m\beta f= -\rmi\sum_{k=1}^n \alpha_k \partial_k f + m \beta f.
\end{equation}
In the sequel the mass $m$ will often play no essential role in our considerations and hence we will typically denote 
the operator in \eqref{eqedm} simply by $D$; whenever the mass term is relevant for an argument we shall indicate it by using the subindex $m$, that is,
$D=D_m$; in most cases this will be done for $m=0$.

From now on we will assume throughout this paper that 
\[
\Omega \subset\RR^n \text{ is a bounded domain (connected open set) with Lipschitz boundary}
\]
and in short we will say that $\Omega$ is a \textit{bounded Lipschitz domain}.
For $s\in[0,1]$ we consider the Sobolev-type spaces
\begin{equation} \label{def_H_alpha}
H^s_\alpha(\Omega,\CC^N):=\big\{f\in H^s(\Omega,\CC^N):\ D_0 f\in L^2(\Omega,\CC^N)\big\},
\end{equation}
which are complete if endowed with norms
\begin{equation}\label{normi}
\| f\|^2_{H_{\alpha}^{s}(\Omega,\CC^N)}:= \| f\|^2_{H^{s}(\Omega,\CC^N)}+ \| D_0 f\|^2_{L^{2}(\Omega,\CC^N)}.
\end{equation}
Observe that for $s\in[0,1]$ and $f\in H^s(\Omega,\CC^N)$ one has $D_0 f\in L^2(\Omega,\CC^N)$ if and only if $D_m f\in L^2(\Omega,\CC^N)$ for some 
$m\in\RR$ and that the norm in \eqref{normi} is equivalent to $\| f\|^2_{H^{s}(\Omega,\CC^N)}+ \| D_m f\|^2_{L^{2}(\Omega,\CC^N)}$.

We start with a standard result on the Dirichlet trace map that can be found in, e.g.
\cite[Lemma~4.1]{BHSLS24}.

\begin{lemma}\label{lem1}
The Dirichlet trace map
		$
		C^\infty(\overline{\Omega},\CC^N)\ni f\mapsto f|_{\partial\Omega}\in L^2(\partial\Omega,\CC^N)
		$
		admits a unique extension to a bounded linear operator 
		\[
		\gamma^{\partial \Omega}: H_{\alpha}^{\half}(\Omega, \CC^N)\mapsto L^2(\partial \Omega,\CC^N).
		\]
\end{lemma}

In the following, let $\cB\in L^\infty(\partial\Omega,\CC^{N\times N})$
be an a.e. defined matrix function 
such that the operator
\begin{align} \label{def_A_mu}
		A f=-\rmi \alpha\cdot\nabla f + m\beta f, \quad
		\dom A &= \bigl\{ f \in  H_{\alpha}^{\half}(\Omega, \CC^N): \cB \gamma^{\partial \Omega} f = 0  \text{ on }\partial\Omega\bigr\},
\end{align}
is self-adjoint in $L^2(\Omega,\CC^N)$. Then $\cB$ will be referred to as an \emph{elliptic boundary condition on $\partial\Omega$}. 
We note that 
the resolvent of the 
operator $A$ is compact, and hence the spectrum of $A$ is purely discrete. The matrix function $\cB$ will remain fixed throughout the rest of the text.

\begin{remark}
A large choice of elliptic boundary conditions on $\partial\Omega$ is available for smooth $\Omega$, which is discussed
in, e.g. \cite{baer,baer2,BHM,BFSV,grosse,OV}, and these results were extended for non-smooth $\partial \Omega$ in \cite{BHSLS24, BKMP26}.		
An important boundary condition suitable for any bounded Lipschitz $\Omega$ is given by the choice
\begin{equation} \label{def_B_mu}
	\cB = \cB_\pm := I_N\pm \rmi \beta \alpha \cdot \nu^{\partial \Omega},
\end{equation}
where $\nu^{\partial\Omega}=(\nu^{\partial\Omega}_1,\dots, \nu^{\partial\Omega}_n)$
is the a.e. defined outer unit normal on $\partial\Omega$. The boundary condition $\cB_\pm \gamma^{\partial \Omega}f=0$
is often referred to as MIT bag (for $+$) and anti-MIT bag (for $-$).
 We note that for this choice one has $\dom A\subset H^1(\Omega,\CC^N)$ if $\Omega$ is smooth \cite{BHM,BFSV,OV} or convex \cite{KPconvex}.
\end{remark}

For subsequent constructions we explicitly state the following assertion:
\begin{lemma}\label{lem2}
	One has 
	\begin{equation}
		\label{inkl1}
		H^1_0(\Omega,\CC^N)\subset \dom A\subset H^1_\loc(\Omega,\CC^N),
	\end{equation}
	and 
	\begin{equation*}
	  C^\infty_c(\Omega,\CC^N)\subset \ran(A - z)\bigl( \dom A\cap C^\infty(\Omega,\CC^N) \bigr),\quad z\in\res A.
	\end{equation*}
	In particular, $A$ is essentially self-adjoint	on 
	\begin{equation}
		\label{dominf}
		\dom_{\infty} A:=\dom A\cap C^\infty(\Omega,\CC^N).
	\end{equation}
\end{lemma}

\begin{proof}
	The left-hand side inclusion in \eqref{inkl1} is clear, and the one on the right-hand side follows by interior elliptic regularity.
	For any $g\in C^\infty_c(\Omega,\CC^N)$ and $z\in\res A$ there is a unique $f\in \dom A$ with $(A-z)f=g$. In particular, $(D-z)f=g$	in $\Omega$,
	and due to the interior elliptic regularity we obtain $f\in C^\infty(\Omega,\CC^N)$. Hence, $C^\infty_c(\Omega,\CC^N)\subset \ran(A-z)\dom_{\infty} A$,
	and the density of $C^\infty_c(\Omega,\CC^N)$ in $L^2(\Omega,\CC^N)$ implies that $\dom_\infty A$ is an operator core of $A$ for any $z\in\res A$; 
	thus $A$ restricted to $\dom_\infty A$ is essentially self-adjoint.
\end{proof}

\subsection{Integral kernels and fundamental solutions}\label{sec22}
Recall that the Euclidean Dirac operator $D$ satisfies the identity
\begin{equation}
	\label{dmdelta}
	(D - z) ( D + z) = \bigl( - \Delta + m^2 - z^2 \bigr) \otimes I_N,\quad z\in\CC,
\end{equation}
where $\Delta$ is the Laplacian in $\RR^n$. For any $\lambda\in\CC$, a fundamental solution $\cE_{\lambda}$
of $-\Delta-\lambda$ is given by
\begin{equation} \label{fundamental_solution_Laplace}
	\cE_\lambda:\ \mathbb{R}^n \setminus \{ 0 \} \rightarrow \mathbb{C}, \quad \cE_\lambda(x) = \begin{cases}
		\dfrac{1}{2\pi}K_0\big(-\rmi\sqrt{\lambda}|x|\big), & \text{if } n=2 \text{ and } \lambda\ne 0,\\[\bigskipamount]
		-\dfrac{1}{2\pi}\log|x|, & \text{if } n=2 \text{ and } \lambda=0, \\[\bigskipamount]
		\dfrac{e^{\rmi \sqrt{\lambda}|x|}}{4\pi|x|},  & \text{if } n=3,
	\end{cases}
\end{equation}
where the square root is fixed by $\Im \sqrt{w}> 0$ for $w \in \CC \setminus[0, \infty)$ and $\sqrt{w} \geq 0$ for $w\in[0,\infty)$,  and $K_j$  are the $j$-th modified Bessel functions of second kind.
It follows from \eqref{dmdelta} that the function
\[
E_{z}:=(D+z)(\cE_{z^2-m^2}\otimes I_N)
\]
is a fundamental solution of $D-z$ for any $z\in\CC$.
We set
\[
k(z):=\sqrt{z^2 - m^2};
\]
then a simple direct computation shows that for any $x \in \mathbb{R}^n  \setminus \{ 0 \}$
\begin{itemize}
	\item if $n=2$:
	\begin{equation} \label{fundamentla_solution_Dirac_2D}
		E_{z}(x)=\begin{cases}
			\dfrac{1}{2 \pi} k(z) K_1 \big( - \rmi k(z) | x |\big)  \dfrac{\alpha \cdot x}{ | x | }\\
			\qquad + \dfrac{1}{2 \pi} K_0 \big(- \rmi k(z) | x |\big) \big( z I_2 + m \beta\big), &z^2 \ne m^2, \\[\bigskipamount]
			-\dfrac{1}{2 \pi}\big( z I_2 + m \beta\big)\log|x| +\dfrac{\rmi}{2 \pi}  \dfrac{\alpha \cdot x}{ | x|^2 }, & z^2=m^2,	
		\end{cases}
	\end{equation}
	\item if $n=3$:
	\begin{equation} \label{fundamentla_solution_Dirac_3D}
		E_{z}(x)= \bigg( z I_4 + m \beta + \big( 1 - \rmi k(z) |x|\big) \dfrac{\rmi  \alpha \cdot x}{  | x |^2}  \bigg) \dfrac{1}{4 \pi | x |} e^{\rmi k(z) |x|}.
	\end{equation}
\end{itemize}
In particular, since $E_{z}$ is a fundamental solution for $D-z$ one has for all $y \in \mathbb{R}^n$ and any $\xi \in \mathbb{C}^N$
\begin{equation} \label{fundamental_solution_translation}
	(D-z)E_{ z}(\cdot -y)\xi=\delta_{y}\xi.
\end{equation}
Furthermore, as $E_{z}$ has a singularity of order $n-1$ at $x=0$ and it is smooth for $x \neq 0$, one has 
\begin{equation} \label{E_z_L_1_loc}
	E_{z} \in L^1_\loc(\mathbb{R}^n, \CC^{N \times N}).
\end{equation}
In our further considerations we will use the function
\begin{align}\label{SolFunSymmetric}
	\Upsilon_z(x):= E_z(x)- E_{\bar{z}}(-x)^\ast, \quad x\in \RR^n\setminus\{0\}, \quad z\in \CC.
\end{align}
Some basic facts about $\Upsilon_z$ are collected in the following lemma.

\begin{lemma} \label{lemma_Upsilon_z}
  Set $I_m:=(-\infty, -|m|)\cup(|m|,\infty)$. Then,    $\Upsilon_z \equiv 0$ for all $z \notin I_m$, while for
  $z\in I_m$ one has
  \begin{align}\label{SolFunSymmetric1}
	\Upsilon_z= 2 \rmi(D+z)(\Im \cE_{k^2(z)}\otimes I_N), \quad \Im \cE_{k(z)^2}(x)= \begin{cases}
		\dfrac{1}{4} J_0(k(z)|x|), \quad  & n=2, \\[\bigskipamount]
		\dfrac{\sin(k(z)|x|)}{4\pi|x|}, & n=3,	
	\end{cases}
\end{align}
where $J_0$ is the Bessel function of the first kind of order $0$. In particular, in all cases $\Upsilon_z$ extends to a real-analytic function on $\RR^n$, which is denoted by the same symbol, and $(D-z)\Upsilon_{z}=0$ in $\RR^n$. 
\end{lemma}
\begin{proof}
  First, using the convention for the complex square root one sees directly that $\Upsilon_z \equiv 0$ for $z \notin I_m$. 
  
  Next, for $z \in I_m$ one has $k(z) > 0$ and thus, $\Upsilon_z= 2 \rmi(D+z)(\Im \cE_{k^2(z)}\otimes I_N)$. To further study this function in dimension $n=2$, we use the identity $\pi \rmi J_0(t) = K_0(-\rmi t) - K_0(\rmi t)$ for $t > 0$, see, e.g., \cite[equation 10.27.9]{DLMF}, yielding for $x \neq 0$
  \begin{align*}
     \Im \cE_{k^2(z)}(x) &= \frac{1}{4 \pi \rmi} \Bigl( K_0(-\rmi k(z) |x|) - K_0(\rmi k(z) |x|) \Bigr)\\
     &= \frac{1}{4} J_0(k(z)|x|) =  \frac{1}{4} \sum_{j=0}^\infty (-1)^j \frac{(k(z)|x|)^{2j}}{4^j (j!)^2},
  \end{align*}
  where in the last step $J_0$ was represented as power series (see, e.g., \cite[equation 10.2.2]{DLMF}). From this it follows that $\Im \cE_{k^2(z)}$ and thus, also $\Upsilon_z$ can be extended to real-analytic functions on $\mathbb{R}^2$.
  
  Likewise, for $n=3$ one gets with Euler's identity and the power series representation of $\sin t$ that
  \begin{equation*}
    \Im \cE_{k^2(z)}(x)= \frac{\sin(k(z)|x|)}{4\pi|x|} = \frac{k(z)}{4 \pi} \sum_{j=0}^{\infty} (-1)^j \frac{(k(z) |x|)^{2j}}{(2j+1)!}
  \end{equation*}
  and therefore, $\Im \cE_{k^2(z)}$ and also $\Upsilon_z$ can be extended to real-analytic functions on $\mathbb{R}^3$.
  
  Finally, while one clearly has $(D-z)\Upsilon_{z}=0$ in $\RR^n$ for $z \notin I_m$, one also has for $z \in I_m$
  \begin{align*}
    (D-z)\Upsilon_{z} &= 2 \rmi(D-z)(D+z)(\Im \cE_{k^2(z)}\otimes I_N)\\
    & = 2 \rmi(-\Delta - k(z)^2)(\Im \cE_{k^2(z)}\otimes I_N) = 0 \quad \text{in } \mathcal{D}'(\mathbb{R}^n)
  \end{align*}
  and thus, due to the smoothness of $\Upsilon_{z}$ also $(D-z)\Upsilon_{z} = 0$ in $\mathbb{R}^n$. This finishes the proof.
\end{proof}

In the next theorem we provide an integral representation of the resolvent of the self-adjoint Dirac operator $A$ in \eqref{def_A_mu} with the help of the fundamental solution $E_z$. The proof is inspired by \cite{ammann, raulot}.

\begin{theorem}[Integral kernel of the resolvent of $A$]\label{theo1}
	For any $z\in\res A$	
	there exists a matrix function $T_{z}\in C^\infty(\Omega\times\Omega,\CC^{N\times N})$
	such that the matrix function $G_{z}$ defined by
	\[
	G_{z}: \bigl\{ (x,y) \in \Omega\times\Omega: x \neq y \bigr\} \rightarrow \CC^{N\times N}, \quad G_{z}(x,y) = E_z(x-y)+T_{z}(x,y),
	\]	
	satisfies $G_{z}\in L^1_\loc(\Omega\times\Omega,\CC^{N\times N})$ and

	\begin{align}\label{res}
		(A -z)^{-1} f(x)=\int_{\Omega} G_{z}(x,y) f(y)\,\dd y, \quad
		x\in\Omega,\ f\in C^\infty_c(\Omega,\CC^N).
	\end{align}
	Moreover, with $\Upsilon_z$ as in \eqref{SolFunSymmetric}, for all $z\in\res A$ one has
	\begin{gather}
		\label{mk}
		T_{z}(x,y)^\ast=T_{\bar z}(y,x) -\Upsilon_z(x-y)^\ast, \quad (x,y)\in\Omega\times\Omega,
	\end{gather}
	and
	\begin{gather}
		\label{mk2}
		G_{z}(x,y)^\ast=G_{\Bar z}(y,x),\quad (x,y)\in\Omega\times\Omega,\, x\not= y.
	\end{gather}
	In particular, $T_{z}(x,y)^\ast=T_{\bar z}(y,x)$ for all $(x,y)\in\Omega\times\Omega$ if $z\in(\res A)\setminus  I_m$. 
\end{theorem}

For the proof we need some preparations. 
Recall that a function $F: M 
\rightarrow L^2(\Omega, \CC^N)$ is differentiable in an open and bounded set $M\subset \mathbb{R}^n$ with respect to 
$x_j$, if there exists a function $\partial_{x_j} F: M \rightarrow 
L^2(\Omega, \CC^N)$ such that for any $x \in M$
\begin{equation*}
	\lim_{t \rightarrow 0} \left\| \frac{F(x+t e_j) - F(x)}{t} - 
	\partial_{x_j} F(x) \right\|_{L^2(\Omega, \CC^N)} = 0,
\end{equation*}
where $e_j$ denotes the $j$-th canonical basis vector in $\mathbb{R}^n$.
In the following we write $C^\infty(M, L^2(\Omega, \CC^N))$ for the set 
of all smooth $L^2(\Omega, \CC^N)$-valued functions. The following 
simple lemma will be useful to show the differentiability of 
$L^2(\Omega, \CC^N)$-valued functions.

\begin{lemma} \label{lemma_differentiable}
	Let $M \subset \mathbb{R}^n$ be open and bounded and let $f \in 
	C^\infty(M \times \overline{\Omega}, \CC^N)$. Then, the function $F: M 
	\rightarrow L^2(\Omega, \CC^N)$ defined by $F(x) := f(x, \cdot)$ is 
	differentiable and $\partial_{x_j} F(x) = (\partial_{x_j} f)(x, \cdot)$.
\end{lemma}
\begin{proof}
	First, as  $f \in C^\infty(M \times \overline{\Omega}, \CC^N)$, the 
	map $F: M \rightarrow L^2(\Omega, \CC^N)$ is well-defined. Moreover, fix $x\in M$ and choose $\delta>0$ with $\overline{B_\delta(x)}\subset M$. For $0<|t|<\delta$ and $y\in\Omega$ one has
	\begin{equation*}
		\dfrac{f(x+t e_j, y) - f(x,y)}{t} - (\partial_{x_j} f)(x, y)=\int_0^1\left(  (\partial_{x_j}f)(x+ st e_j,y)-(\partial_{x_j}f)(x,y)\right)\dd s,
	\end{equation*}
	and as $\partial_{x_j}f$ is uniformly continuous on the compact set $\overline{B_\delta(x)}\times \overline{\Omega}$ the right-hand side converges to $0$ uniformly in $y\in\Omega$, as $t \rightarrow 0$. 
	Therefore,
	\begin{equation*}
		\begin{split}
			\lim_{t \rightarrow 0} &\left\| \frac{F(x+t e_j) - F(x)}{t} - 
			(\partial_{x_j} f)(x, \cdot) \right\|_{L^2(\Omega, \CC^N)}^2 \\
			&= \lim_{t \rightarrow 0}\int_{\Omega} \left| \frac{f(x+t e_j,y) - f(x,y)}{t} - 
			(\partial_{x_j} f)(x, y) \right|^2 \dd y = 0,
		\end{split}
	\end{equation*}
	which shows that $F$ is differentiable and $\partial_{x_j} F(x) = 
	(\partial_{x_j} f)(x, \cdot)$. 
\end{proof}

\begin{lemma}\label{prop1}
	For any $z\in\res A$ and any $y\in\Omega$ there is a unique matrix function
	$\cM_y^z:\Omega\rightarrow \CC^{N\times N}$ such that for any $\xi\in \CC^N$ one has $\cM^z_y(\cdot) \xi\in H^\frac{1}{2}(\Omega,\CC^N)$ 
	and the function
	\[
	\cG^z_y:\Omega \setminus \{ y \}\rightarrow \CC^{N\times N},\quad x\mapsto \cG^z_y(x)=E_{\bar z}(x - y) +\cM_y^z(x)
	\]
	satisfies the identity $(D-\bar z)\cG^z_y\xi=\delta_y \xi$ in $\mathcal D'(\Omega,\CC^N)$ and the boundary condition for $A$, i.e. for every $\chi\in C^\infty(\overline\Omega)$ with $\chi=1$ near $\partial\Omega$ and $\chi=0$ near $y$ one has
	\[
	\chi \cG^z_y\xi\in H_{\alpha}^\frac{1}{2}(\Omega, \CC^N)\quad\text{and}\quad\cB \gamma^{\partial\Omega}(\chi \cG^z_y\xi)=0 \text{ on }\partial\Omega.
	\]
	Moreover, 
	\begin{equation}
		\label{gzy}
		(y\mapsto \cM^z_y\xi) \in C^\infty(\Omega, L^2(\Omega,\CC^N)).
	\end{equation}
\end{lemma}

\begin{proof} We first prove the uniqueness. Let $z\in\res A$, $y\in\Omega$, $\xi\in\CC^N$, and suppose that there are $\cM_y^z$ and $\tilde{\cM}^z_y$ such that the matrix functions
	\[
	x\mapsto \cG^z_y(x)=E_{\bar z}(x-y)+ \cM_y^z(x), \quad x\mapsto\Tilde{\cG}^z_y(x)=E_{\bar z}(x-y)+ \Tilde{\cM}^z_y(x),
	\quad x\not= y,
	\]
	satisfy the assumptions. Then $(\cM_y^z-\Tilde{\cM}^z_y)\xi\in H^\frac{1}{2}(\Omega,\CC^N)$ and from
	\begin{equation}\label{jussi1}
		(D-\Bar z)(\cM_y^z-\Tilde{\cM}^z_y)\xi=(D-\Bar z)(\cG_y^z-\Tilde{\cG}^z_y)\xi=\delta_y\xi-\delta_y\xi=0
	\end{equation}
	we conclude $D(\cM_y^z-\Tilde{\cM}^z_y)\xi\in L^2(\Omega,\CC^N)$ and also $D_0(\cM_y^z-\Tilde{\cM}^z_y)\xi\in L^2(\Omega,\CC^N)$. Hence we have
	$(\cM_y^z-\Tilde{\cM}^z_y)\xi\in H_\alpha^\frac{1}{2}(\Omega,\CC^N)$; cf. \eqref{def_H_alpha}. Furthermore, for  $\chi\in C^\infty(\overline\Omega)$ with $\chi=1$ near $\partial\Omega$ and $\chi=0$ near $y$
	we have by assumption
	$$\chi (\cM_y^z-\Tilde{\cM}^z_y)\xi=\chi(\cG_y^z-\Tilde{\cG}^z_y)\xi\in H_{\alpha}^\frac{1}{2}(\Omega, \CC^N)$$
	and
	\[
	\cB \gamma^{\partial\Omega}\big(\chi (\cM_y^z-\Tilde{\cM}^z_y)\xi\big)
	=\cB \gamma^{\partial\Omega}\big(\chi(\cG_y^z-\Tilde{\cG}^z_y)\xi\big)=0,
	\]
	which shows $\cB \gamma^{\partial\Omega}\big((\cM_y^z-\Tilde{\cM}^z_y)\xi\big)=0$. Therefore,
	$(\cM_y^z-\Tilde{\cM}^z_y)\xi\in\dom A$ and \eqref{jussi1} implies $(A-\Bar z)(\cM_y^z-\Tilde{\cM}^z_y)\xi=0$.
	Since the operator $A -\Bar z$ is injective and $\xi\in\CC^N$ is arbitrary, it follows that $\cM_y^z=\Tilde{\cM}^z_y$.

	Now we prove the existence. Let $y_0\in\Omega$ and $r>0$ with $B_{2r}(y_0)\subset\Omega$. Let $\rho\in C^\infty_c(\Omega)$  such that
	$    \rho= 1$ in $B_{r}(y_0)$ and $\supp \rho\subset B_{2r}(y_0)$.
	Let $y\in B_{r}(y_0)$, consider the matrix function $x\mapsto R(x):=-\rmi \alpha\cdot\nabla\rho(x)$, and note that~\eqref{fundamental_solution_translation}
	implies
	\begin{equation}\label{dff}
		\begin{split}
			(D-\bar z)\big( \rho E_{\bar z}(\cdot -y)\xi\big)&= \rho (D-\bar z)E_{\bar z}(\cdot -y)\xi + R E_{\bar z}(\cdot -y)\xi\\
			&=\delta_y \xi + R E_{\bar z}(\cdot-y)\xi
		\end{split}
	\end{equation}
	for any $\xi\in\CC^N$.
	Since the matrix function $R$ is smooth with support in $B_{2r}(y_0)$ and vanishes in
	$B_r(y_0)$, we conclude that the second term on the right-hand side in \eqref{dff}
	belongs to $L^2(\Omega;\CC^N)$. Therefore, there exists a matrix function $F^z_y$
	with
	\[
	F^z_y\xi\in\dom A, \quad (A -\bar z)F^z_y\xi=R E_{\bar z}(\cdot -y)\xi,
	\]
	and \eqref{dff} takes the form
	\[
	(D-\bar z)\big[\big(\rho E_{\bar z}(\cdot -y)-F^z_y\big)\xi\big]=\delta_y \xi.
	\]
	Hence, the matrix function $\cG_y^z:= \rho E_{\bar z}(\cdot - y) -F_y^z$
	satisfies by construction the identity $(D-\bar z)\cG^z_y\xi=\delta_y \xi$.
	In particular, for $\chi\in C^\infty(\overline\Omega)$ with $\chi=1$ near $\partial\Omega$ and $\chi=0$ near $y$ 
	one has $D (\chi \cG^z_y \xi) \in L^2(\Omega, \mathbb{C}^N)$ and $D_0 (\chi \cG^z_y \xi) \in L^2(\Omega, \mathbb{C}^N)$. Furthermore, as $\chi \rho E_{\bar z}(\cdot - y)$ is smooth and $F_y^z \xi \in \dom A \subset H^{\half}_\alpha(\Omega, \mathbb{C}^N)$ we get $\chi \cG^z_y \xi \in H^{\half}_\alpha(\Omega, \mathbb{C}^N)$. Next, since
	$\rho E_{\bar z}(\cdot - y)$ vanishes near $\partial\Omega$,  one has
	\[
	\cB\gamma^{\partial\Omega} (\chi \cG_y^z\xi)=-\cB \gamma^{\partial\Omega}(\chi F^z_y\xi)=0 \text{ on } \partial\Omega.
	\]
	By denoting $\cM^z_y:=(\rho-1)E_{\bar z}(\cdot -y)-F^z_y$ we have $\cG^z_y = E_{\bar z}(\cdot-y) +\cM^z_y$.
	As $(\rho-1)E_{\bar z}(\cdot -y)$
	is smooth on $\overline{\Omega}$ and $F^z_y\xi\in\dom A$, it follows that $\cM^z_y\xi\in H^\frac{1}{2}_\alpha(\Omega,\CC^N)$ and, in particular, 
	$\cM^z_y\xi\in H^\frac{1}{2}(\Omega,\CC^N)$.

	It remains to show~\eqref{gzy}, which is essentially a consequence of Lemma~\ref{lemma_differentiable}. In fact, with 
	$R=-\rmi \alpha\cdot\nabla\rho$ and  $\xi\in\mathbb{C}^N$
	consider the function
	\begin{equation*}
		F_1: B_{\frac{r}{2}}(y_0) \rightarrow L^2(\Omega,\CC^N),\qquad F_1(y) := R E_{\bar z}(\cdot -y)\xi,
	\end{equation*}
	and observe that
	\begin{equation*}
		(x,y)\mapsto f_1(x,y):= R(x) E_{\bar z}(x-y) \xi \in C^\infty(\overline{\Omega} \times \overline{B_{\frac{r}{2}}(y_0)}, \mathbb{C}^N).
	\end{equation*}
	Then it follows from Lemma~\ref{lemma_differentiable} that $\partial_{y_i}F_1(y)=(\partial_{y_i} f_1)(\cdot,y)=\partial_{y_i}R E_{\bar z}(\cdot-y)\xi$.
	By iterating this procedure one concludes
	$(\partial^\gamma F_1)(y) = \partial^\gamma_y RE_{\bar z}(\cdot-y)\xi \in L^2(\Omega, \mathbb{C}^N)$ for any $\gamma\in\NN^n_0$, which shows
	\begin{equation} \label{smooth1}
		F_1 \in C^\infty\bigl(B_{\frac{r}{2}}(y_0), L^2(\Omega,\CC^N)\bigr).
	\end{equation}
	
	Since the resolvent $(A-\bar z)^{-1}$ is bounded in $L^2(\Omega,\mathbb{C}^N)$, one gets that 
	\begin{equation*}
		y\mapsto F^z_y \xi=(A -\bar z)^{-1} F_1(y) \in C^\infty\bigl(B_{\frac{r}{2}}(y_0), L^2(\Omega,\CC^N)\bigr).
	\end{equation*}
	Next, in a similar way as in~\eqref{smooth1}, one shows for the map $F_2: B_{\frac{r}{2}}(y_0) \rightarrow L^2(\Omega,\CC^N)$ defined by $F_2(y) := (\rho-1)E_{\bar z}(\cdot -y)\xi$ that
	\begin{align*}
		F_2 \in C^\infty\bigl(B_{\frac{r}{2}}(y_0), L^2(\Omega,\CC^N)\bigr).
	\end{align*}
	Therefore, taking the definition $\cM^z_y=(\rho-1)E_{\bar z}(\cdot -y)-F^z_y$ into account, one concludes that 
	\begin{equation*}
		(y\mapsto \cM^z_y\xi) \in C^\infty(B_{\frac{r}{2}}(y_0), L^2(\Omega,\CC^N)).
	\end{equation*}
	Since this is true for any fixed $y_0 \in \Omega$, the uniqueness implies that $(y\mapsto \cM^z_y\xi)$ is smooth in $\Omega$, i.e. that \eqref{gzy} is true.
\end{proof}

\begin{lemma}\label{lemma1}
	For any $f\in\dom_\infty A= \dom A\cap C^\infty(\Omega,\CC^N)$ there holds
	\begin{align*}
		\int_{\Omega} (\cG_y^z)^\ast(x) (A - z)f(x)\dd x=f(y) \quad \text{ for all }
		y\in\Omega.
	\end{align*}
\end{lemma}
\begin{proof}
	Let $f\in\dom_\infty A$ and $y\in\Omega$. Pick an arbitrary $\xi\in\CC^N$ and set
	\[
	J :=\int_{\Omega} \big\langle \xi, (\cG_{y}^z)^\ast(x) (A- z)f(x)\big\rangle_{\CC^N}\dd x.
	\]
	By linearity, it suffices to show that $J=\langle \xi,f(y)\rangle_{\CC^N}$.
	
	Let $r>0$ with $B_{2r}(y)\subset\Omega$ and $\rho\in C^\infty_c(\Omega)$  such that
	$\rho= 1$ in $B_{r}(y)$ and $\supp \rho\subset B_{2r}(y)$. Then $J=J_1+J_2$ with
	\begin{align*}
		J_1 &:=\int_{\Omega} \big\langle \xi, (\cG_{y}^z)^\ast(x) (A- z)(\rho f)(x)\big\rangle_{\CC^N}\dd x,\\
		J_2 &:=\int_{\Omega} \big\langle \xi, (\cG_{y}^z)^\ast(x) (A- z)\big[(1-\rho) f\big](x)\big\rangle_{\CC^N}\dd x.
	\end{align*}
	Due to $\rho f\in C^\infty_c(\Omega,\CC^N)$ we compute $J_1$ as follows:
	\begin{align*}
		J_1 &= \big\langle \cG_{y}^{z}\xi, \overline{(D -z)(\rho f)}\,\big\rangle_{\cD'(\Omega,\CC^N), \cD(\Omega,\CC^N)}\\
		&= \big\langle (D - \Bar z)\cG_{y}^{z} \xi, \overline{\rho f}\big\rangle_{\cD'(\Omega,\CC^N), \cD(\Omega,\CC^N)}\\
		&= \big\langle \delta_{y} \xi, \overline{\rho f}\big\rangle_{\cD'(\Omega,\CC^N), \cD(\Omega,\CC^N)}\\
		& = \big\langle \xi,  (\rho f)(y)\big\rangle_{\CC^N}\\
		& =\big\langle \xi, f(y)\big\rangle_{\CC^N},
	\end{align*}
	and it remains to show that $J_2=0$. For that, choose $\Tilde \rho\in C^\infty(\overline\Omega)$ such that  $\Tilde \rho =0$ on $B_{\frac{r}{2}}(y)$ and $\Tilde\rho=1$ on $\Omega\setminus B_{\frac{3r}{4}}(y)$; as $1-\rho=0$ on $B_r(y)$, one has $\Tilde \rho=1$ on $\supp(1-\rho)$ and $\nabla \Tilde\rho$ and $1-\rho$ have disjoint support. Then
	\[
	J_2=\int_{\Omega} \big\langle \Tilde\rho(x)\cG_{y}^z(x)\xi, (D- z)\big[(1-\rho) f\big](x)\big\rangle_{\CC^N}\dd x.
	\]
	It follows from Lemma~\ref{prop1} and $\Tilde \rho =0 $ near $y$  (and $\Tilde\rho=1$ near $\partial\Omega$)  that $\Tilde\rho\cG_y^z(\cdot)\xi\in\dom A$ for any $\xi\in\CC^N$;  moreover, $(1-\rho)f\in\dom A$, since $(1-\rho)f\in H_\alpha^{\frac{1}{2}}(\Omega,\CC^N)$ and $\gamma^{\partial\Omega}((1-\rho)f)=\gamma^{\partial\Omega}f$.  Thus
	\begin{align*}
		J_2&=\int_{\Omega} \big\langle \tilde\rho(x)\mathcal{G}_y^z(x)\xi, (D- z)\big[(1-\rho)f\big](x)\big\rangle_{\CC^N}\dd  x\\
		&=\big\langle \Tilde\rho\cG_y^z\xi, (A- z)\big[(1-\rho)f\big] \big\rangle_{L^2(\Omega,\CC^N)}\\
		&=\big\langle (A- \bar z)(\Tilde\rho\cG_y^z\xi), (1-\rho) f\big\rangle_{L^2(\Omega,\CC^N)}.
	\end{align*}
	From $\Tilde \rho =0 $ near $y$ we also have $\Tilde\rho \delta_y=0$ and hence
	\begin{align*}
		(A- \bar z)(\Tilde\rho\cG_y^z\xi)&=\Tilde\rho (D- \bar z)\cG_y^z\xi + (-\rmi\alpha\cdot\nabla\Tilde\rho) \cG_y^z\xi\\
		&= \Tilde\rho \delta_y \xi + (-\rmi\alpha\cdot\nabla\Tilde\rho)\cG_y^z\xi\\
		&=(-\rmi\alpha\cdot\nabla\Tilde\rho)\cG_y^z\xi.
	\end{align*}
	Therefore, using that the functions $-\rmi\alpha\cdot\nabla\Tilde\rho$ and $1-\rho$ have disjoint supports we obtain  $J_2=0$.
\end{proof}

\begin{proof}[\bf Proof of Theorem \ref{theo1}]	
	For $z\in\res A$ and $x,y\in\Omega$ define
	\[
	G_z(x,y):=\cG^z_x(y)^*,\,\,\,x\not= y,\quad\text{and}\quad  T_z(x,y):=\cM^z_x(y)^*-\Upsilon_z(x-y)
	\]
	with $\cG^z_x$ and $\cM^z_x$ from Lemma~\ref{prop1} (now viewed as functions in the $y$-variable)  and $\Upsilon_z$ from \eqref{SolFunSymmetric}. As $\cG_x^z(y)= E_{\bar{z}}(y-x)+ \cM_x^z(y)$ and $E_{\bar{z}}(y-x)^\ast=E_z(x-y)-\Upsilon_z(x-y)$ by \eqref{SolFunSymmetric}, we conclude for all $x,y\in\Omega$ with $x \neq y$ that
	\begin{equation} \label{decomposition_G}
	G_z(x,y)=E_z(x-y)+T_z(x,y).
	\end{equation}
	In view of \eqref{gzy},  for all $j,k\in\{1, \ldots,N\}$ and each compact $\omega\subset\Omega$ it holds that
	\begin{equation} \label{T_z_L_2_loc}
		\begin{split}
			\int_\omega \int_{\Omega}\big|(\cM_x^z)_{j,k}(y)\big|^2\,\dd y\,\dd x = \int_\omega \big\|(\cM^z_x)_{j,k}\big\|^2_{L^2(\Omega)}\dd x <\infty;
		\end{split}
	\end{equation}
	here we use that the continuous map $\omega \ni x\mapsto (\cM_x^z)_{j,k}\in L^2(\Omega)$ defines an element of $L^2(\omega; L^2(\Omega)) \simeq L^2(\omega\times\Omega)$, which also gives a meaning to $ (x,y)\mapsto \cM_x^z(y)$ as a measurable function on $\Omega\times\Omega$. As $\Upsilon_z$ is smooth  by Lemma~\ref{lemma_Upsilon_z},  $(T_z)_{j,k}\in L^2_\loc(\Omega\times\Omega) \subset L^1_\loc(\Omega\times\Omega)$. Since $E_z \in L^1_\loc(\mathbb{R}^n, \CC^{N \times N})$ by~\eqref{E_z_L_1_loc}, we conclude that $G_z \in L^1_\loc(\Omega \times \Omega, \CC^{N \times N})$.
	Since   $C^\infty_c(\Omega,\CC^N)\subset (A - z) \dom_\infty A$   by Lemma~\ref{lem2}, one gets with Lemma~\ref{lemma1} 
	\[
	(A-z)^{-1}g(x)=\int_\Omega G_z(x,y)g(y)\,\dd y,\qquad g\in C^\infty_c(\Omega,\CC^N),\ x\in\Omega.
	\]
	Due to the last identity and the local integrability for any $f,g\in C^\infty_c(\Omega,\CC^N)$ one uses Fubini's theorem
	to show that
	\begin{equation*}
		\begin{split}
			\int_\Omega\int_\Omega \big\langle f(x),
			G_z(x,y) g(y)\big\rangle_{\CC^N}\dd y\,\dd x&=\big\langle f,(A-z)^{-1}g\big\rangle_{L^2(\Omega,\CC^N)}\\
			&=\big\langle (A-\Bar z)^{-1}f,g\big\rangle_{L^2(\Omega,\CC^N)}\\
			&=
			\int_\Omega\int_\Omega \big\langle G_{\Bar z}(y,x) f(x),
			g(y)\big\rangle_{\CC^N}\dd x\,\dd y\\
			&=\int_\Omega\int_\Omega \big\langle f(x),
			G_{\Bar z}(y,x)^*g(y)\big\rangle_{\CC^N}\dd x\,\dd y\\
			&=\int_\Omega\int_\Omega \big\langle f(x),
			G_{\Bar z}(y,x)^*g(y)\big\rangle_{\CC^N}\dd y\,\dd x,
		\end{split}
	\end{equation*}
	which yields $G_z(x,y)=  G_{\Bar z}(y,x)^*$ for a.e. $(x,y)$,
	Moreover, in view of~\eqref{decomposition_G} and the equality $E_{\Bar z}(y-x)^*=E_z(x-y)-\Upsilon_z(x-y)$, see \eqref{SolFunSymmetric}, we also obtain
	\begin{equation}
		\label{ksym}
		T_z(x,y)=  T_{\Bar z}(y,x)^* -\Upsilon_z(x-y)  \text{ for a.e. } (x,y).
	\end{equation}
	It remains to show the smoothness of $T_z$ (which then also implies that \eqref{ksym} is valid for all $(x,y) \in \Omega \times \Omega$, i.e. 
	\eqref{mk} and hence also \eqref{mk2} hold): As this smoothness is of crucial importance for all subsequent constructions, we provide a detailed argument. To show the claim, we verify that all distributional derivatives of any component of $T_z$ with respect to $x$- and $y$-variables are regular distributions, and then use an elliptic regularity argument to deduce the claimed smoothness.

	Pick any $j,k\in\{1,\dots,N\}$ and set $h:=(T_z)_{j,k}\in L^2_\loc(\Omega\times\Omega)$ and 
	\[
	H: \Omega \rightarrow L^2(\Omega), \quad H(x) := h(x,\cdot).
	\]
	By Lemma~\ref{prop1}, Lemma \ref{lemma_differentiable}, and the smoothness of  $\Upsilon_z$ from Lemma~\ref{lemma_Upsilon_z}, the map $H$ is well-defined and $H \in C^\infty(\Omega, L^2(\Omega))$. Let $i\in\{1,\dots,n\}$ and $e_i$ be the $i$-th canonical basis vector in $\RR^n$. Then for any $\varphi\in C^\infty_c(\Omega \times \Omega)$ and $\Phi \in C_c^\infty(\Omega, L^2(\Omega))$  defined by $\Phi(x):=\varphi(x,\cdot)$ we have $\partial_{x_i} \Phi(x)=\partial_{x_i}\varphi (x,\cdot)$ (see Lemma~\ref{lemma_differentiable}) and we compute
	using integration by parts (where $\dd s$  denotes the usual $(n-1)$-dimensional Hausdorff measure on $\partial\Omega$ and $(\nu^{\partial\Omega})_i$
	is the $i$-th component of the unit normal field pointing outwards) 
	\begin{equation}
		\label{lll}
		\begin{split}
			\int_\Omega \int_\Omega (\partial_{x_i}H(x))&(y)\overline{\varphi(x,y)}\,\dd y\,\dd x =\int_\Omega\langle \partial_{x_i} H(x), \Phi(x) \rangle_{L^2(\Omega)}\dd x\\
			&=\int_\Omega \partial_{x_i}\langle H(x), \Phi(x) \rangle_{L^2(\Omega)}\dd x-\int_\Omega \langle H(x),\partial_{x_i} \Phi(x) \rangle_{L^2(\Omega)}\dd x\\
			&=\int_{\partial\Omega}(\nu^{\partial\Omega})_i \langle H(x),\Phi(x)\rangle_{L^2(\Omega)}\dd  s(x) -\int_\Omega \langle H(x),\partial_{x_i}\Phi(x)\rangle_{L^2(\Omega)}\dd x\\
			&=-\int_\Omega \langle H(x),\partial_{x_i}\Phi(x)\rangle_{L^2(\Omega)}\dd x\\
			&=-\int_\Omega \int_\Omega h(x,y)\overline{\partial_{x_i}\varphi (x,y)}\,\dd y\,\dd x,
		\end{split}
	\end{equation}
	where in the last but one step it was used that $\Phi(x)=0$ for all $x\in\partial\Omega$. Therefore, 
	$h$ admits a weak derivative with respect to $x_i$ given by
	$\partial_{x_i} h(x,y)=(\partial_{x_i} H(x))(y)$.
	Since $\partial_{x_i} H \in C^\infty(\Omega, L^2(\Omega))$ one shows as in~\eqref{T_z_L_2_loc} that $\partial_{x_i}h\in L^2_\loc(\Omega\times\Omega)$. By iterating this procedure we conclude
	$\partial^\gamma_x h\in L^2_\loc(\Omega\times\Omega)$ for any $\gamma\in\NN^n_0$.
	
	The preceding discussion shows that
	\[
	\partial^\gamma_x (T_z)_{j,k}\in L^2_\loc(\Omega\times\Omega)
	\text{ for all } z\in\CC\setminus\spec A,\ j,k\in\{1,\dots,N\},\ \gamma\in\NN^n_0.
	\]
	At the same time, by \eqref{ksym} we have  $(T_{z})_{j,k}(x,y)=\overline{(T_{\Bar z})_{k,j}(y,x)} -(\Upsilon_z)_{j,k}(x-y)$; as $\Upsilon_{z}$ is smooth, the preceding discussion applied to $T_{\Bar z}$ yields 
	\[
	\partial^\gamma_y (T_z)_{j,k}\in L^2_\loc(\Omega\times\Omega)
	\text{ for all } z\in\CC\setminus\spec A,\ j,k\in\{1,\dots,N\},\ \gamma\in\NN^n_0.
	\]
	In particular, for any $p\in\NN$ and any $j,k\in\{1,\dots,N\}$ we have
	\[
	\Big[(\Delta_{x})^p+ (\Delta_{y})^p\Big] (T_z)_{j,k}\in L^2_\loc(\Omega\times\Omega),
	\]
	and,  since the constant coefficient operator $(\Delta_{x})^p+ (\Delta_{y})^p$ is elliptic of order $2p$ in $\RR^{2n}$, the interior elliptic regularity implies $(T_z)_{j,k}\in H^{2p}_\loc(\Omega\times \Omega)$.
	As $p\in\NN$ is arbitrary, the Sobolev embedding theorem yields $(T_z)_{j,k}\in C^\infty(\Omega\times\Omega)$.
\end{proof}

Finally, we provide an asymptotic expansion of the Green's function $E_{z}$ from~\eqref{fundamentla_solution_Dirac_2D} \& \eqref{fundamentla_solution_Dirac_3D} which will be useful for our further analysis. The technical proof is presented in Appendix~\ref{appendix_pseudos}.

\begin{lemma} \label{lemma_expansion1}
	Define the functions $w^\pm: \mathbb{R}^n \setminus \{ 0 \} \rightarrow \mathbb{C}^{N/2 \times N/2}$ by 
	\begin{equation} \label{def_w_pm}
		w^\pm(x)=\begin{cases}
			\dfrac{\rmi}{2\pi}\dfrac{1}{x_1\pm \rmi x_2}, & n=2,\\[\bigskipamount]
			\dfrac{\rmi (\sigma\cdot x)}{4\pi |x|^3}, & n=3,
		\end{cases}
	\end{equation}
	let $\cE_0$ be as in \eqref{fundamental_solution_Laplace}, set for $z \in \mathbb{C}$
	\begin{equation*}
		\cH_{z}(x) = \begin{pmatrix} (z+m) \cE_0(x)\otimes I_\nhalf & w^+(x) \\ w^-(x) & (z-m) \cE_0(x)\otimes I_\nhalf \end{pmatrix},\quad x \in \mathbb{R}^n \setminus \{ 0 \},
	\end{equation*}
	and define with $\chi \in C_c^\infty(\mathbb{R}^n)$
	\begin{equation*}
		\mathcal{K}_{z} := \chi \big(E_{z} - \cH_{z} \bigr).
	\end{equation*}
	Then, the following is true:
	\begin{itemize}
      \item[(a)] $\mathcal{K}_{z} \in L^1(\mathbb{R}^n, \CC^N) \cap L^\infty(\mathbb{R}^n, \CC^N)$;
      \item[(b)] $\mathcal{K}_z$ is smooth outside the origin and rapidly decreasing (with all of its derivatives) at infinity;
      \item[(c)] $\mathcal{K}_z$ has an expansion of the form
	\begin{equation} \label{expansion_K_z}
		\mathcal{K}_{z}(x) \sim \sum_{\ell \geq 0} \bigl( q_\ell(x) + p_\ell(x) \log |x| \bigr),
	\end{equation}
	where each $q_\ell$ is a homogeneous distribution of degree $3+\ell-n$ that is smooth outside the origin and $p_\ell$ is a polynomial homogeneous of degree $3+\ell-n$;
	\item[(d)]  for $n=2$ the function $\mathcal{K}_{z}$ is continuous and for every $R>0$ there is $c>0$ such that $|\nabla \mathcal K_z(x)|\leq c(1+|\log|x||)$ for all $0<|x|\leq R$;
	\item[(e)] for $n=3$ and every $R>0$ there is $c>0$ such that $|\nabla \mathcal K_z(x)|\leq c|x|^{-1}$ for all $0<|x|\leq R$.
	\end{itemize}
\end{lemma}

The asymptotic expansion in~\eqref{expansion_K_z} is understood in the following sense: For any $L \in \mathbb{N}_0$ there exists $\ell_0 \in \mathbb{N}_0$ such that
\begin{equation*}
	\mathcal{K}_{z}(x) - \sum_{\ell = 0}^{\ell_0} \bigl( q_\ell(x) + p_\ell(x) \log |x| \bigr) \in C^L(\mathbb{R}^n).
\end{equation*}

\section{Dirac operator with a shell interaction on a bounded Lipschitz domain}

In this section, we rigorously define the Dirac operator $A^{\eta,\tau}$ with an electrostatic $\delta$-shell interaction of strength $\eta$ and 
a
Lorentz scalar $\delta$-shell interaction of strength $\tau$ on a compact surface $\Sigma\subset\Omega$, as well as a self-adjoint boundary condition on $\partial\Omega$.
First, we establish the geometric setup, function spaces, and trace operators along $\Sigma$. We then introduce the  operator $A^{\eta,\tau}$ by interpreting the singular interaction as a transmission condition. Finally, we formulate this setup within the abstract framework of boundary triples and Weyl functions, providing the essential foundation for the subsequent spectral analysis.

\subsection{Geometric assumptions, function spaces, and trace operators}\label{sec31}

Let $\Sigma\subset\Omega$ be a compact $C^\infty$-smooth hypersurface (i.e. a not necessarily connected curve for $n=2$ or surface for $n=3$)
with finitely many connected components $\Sigma_p$, $p\in\{1,\dots,P\}$.
Due to the Jordan-Brouwer separation theorem, for each $p$ one has $\RR^n\setminus\Sigma_p=U_p^+\cup U_p^-$,
where $U_p^+$ and $U_p^-$ are disjoint $C^\infty$-domains such that $U^+_p$ is bounded and $U_p^-$ is unbounded. In the following we will only make use of the parts of $U_p^\pm$ that 
belong to $\Omega$ and we use the notation $\Omega_p^\pm:=U_p^\pm \cap\Omega$. 
We denote by $\nu^\Sigma$ the unit normal field on $\Sigma$ which points towards $\Omega_p^-$ on each $\Sigma_p$,
i.e. $\nu^{\Sigma}$ is the exterior unit normal for each $\Omega_p^+$.
For a function $f\in L^2(\Omega)$ we use the notation $f_{p\pm}$ for the restriction of $f$ onto $\Omega_p^\pm$.

For $s\in[0,1]$ consider the Sobolev-type spaces $H^s_\alpha(\Omega_p^\pm,\CC^N)$ defined as in~\eqref{def_H_alpha} equipped with the corresponding Hilbert space norms, and let 
\begin{equation*}
\begin{split}
 H^s_\alpha(\Omega\setminus\Sigma_p,\CC^N)&=H^s_\alpha(\Omega_p^+,\CC^N)\times H^s_\alpha(\Omega_p^-,\CC^N)\\
 &=\bigl\{f\in H^s(\Omega\setminus\Sigma_p,\CC^N):\ D_0 f\in \mathcal D'(\Omega\setminus\Sigma_p)\cap L^2(\Omega,\CC^N)\bigr\},
 \end{split}
\end{equation*}
where $D_0 f\in \mathcal D'(\Omega\setminus\Sigma_p)$ is the distribution 
\begin{equation*}
 C_c^\infty(\Omega\setminus\Sigma_p, \CC^N)\ni \varphi\mapsto \big\langle D_0f,  
 \varphi\big\rangle_{\cD'(\Omega \setminus \Sigma_p,\CC^N), \cD(\Omega \setminus \Sigma_p,\CC^N)},
\end{equation*}
and the requirement $D_0 f\in \mathcal D'(\Omega\setminus\Sigma_p)\cap L^2(\Omega,\CC^N)$ means that this distribution is 
represented by a function $g \in L^2(\Omega, \CC^N)$, i.e. 
\begin{equation*}
  \big\langle D_0f,  
 \varphi\big\rangle_{\cD'(\Omega \setminus \Sigma_p,\CC^N), \cD(\Omega \setminus \Sigma_p,\CC^N)} = \int_\Omega g \cdot \varphi \, \dd x, \quad \varphi \in C_c^\infty(\Omega\setminus\Sigma_p, \CC^N);
\end{equation*}
we will again denote this representing function by $D_0 f$, i.e. $D_0f = g$.
In particular, elliptic regularity  implies that $f \in H^s_\alpha(\Omega\setminus\Sigma_p,\CC^N)$ is locally $H^1$-regular in $\Omega\setminus\Sigma_p$ (but not in all of $\Omega$).
As $\Sigma_p$ is $C^\infty$-smooth it follows from \cite[Lemma~4.1 and Corollary~4.6]{BHSLS24} that there
exist unique continuous trace operators 
\begin{equation}\label{gammap}
\gamma^{\Sigma_p}_\pm: H^s_\alpha(\Omega\setminus\Sigma_p,\CC^N)\mapsto H^{s-\half}(\Sigma_p,\CC^N)
\end{equation}
such that for $f=(f_{p+},f_{p-})\in H^s_\alpha(\Omega\setminus\Sigma_p,\CC^N)$ one has
$\gamma^{\Sigma_p}_+ f =f_{p+}|_{\Sigma_p}$ if $f_{p+}\in C^\infty(\overline{\Omega_p^+})$,
and $\gamma^{\Sigma_p}_- f =f_{p-}|_{\Sigma_p}$ if $f_{p-}\in C^\infty(\overline{\Omega_p^-})$.

In the following we work with traces of functions  
having different Sobolev regularity properties near
$\Sigma=\bigcup_{p=1}^P\Sigma_p$ and near $\partial\Omega$. To take these aspects into account we first define (in the same way as above) for 
$s\in[0,1]$ the Sobolev-type spaces
\begin{equation}
 H^s_\alpha(\Omega\setminus\Sigma,\CC^N)=\bigl\{f\in H^s(\Omega\setminus\Sigma,\CC^N):\ D_0 f\in \mathcal D'(\Omega\setminus\Sigma)\cap L^2(\Omega,\CC^N)\bigr\},
\end{equation}
and the function spaces
\begin{equation}
 H^{s,\Sigma}_\alpha(\Omega\setminus\Sigma,\CC^N):=\Big\{f\in L^2(\Omega,\CC^N): \chi f \in H^s_\alpha(\Omega\setminus\Sigma,\CC^N)\,\text{for any }\chi\in C^\infty_c(\Omega)\Big\}
 \end{equation}
and  
\begin{equation}\label{ueblerraum} 
\begin{split}
 H^{s,\partial\Omega}_\alpha(\Omega\setminus \Sigma,\CC^N):=\Big\{f\in L^2(\Omega,\CC^N)&: \chi f \in H^s_\alpha(\Omega,\CC^N)\\ &\text{ for any }\chi\in C^\infty(\overline \Omega)
	 \text{ with }\chi=0 \text{ near }\Sigma\Big\}.
\end{split}
	 \end{equation}
Roughly speaking, the functions in $H^{s,\Sigma}_\alpha(\Omega\setminus\Sigma,\CC^N)$	 
 are locally $H^1$-regular in $\Omega\setminus\Sigma$ (since $D_0 (\chi f)\in \mathcal D'(\Omega\setminus\Sigma)\cap L^2(\Omega,\CC^N)$)
and $H^s$-regular near $\Sigma$. Similarly, the functions in $H^{s,\partial\Omega}_\alpha(\Omega\setminus\Sigma,\CC^N)$	 
 are locally $H^1$-regular in $\Omega\setminus\Sigma$
and $H^s$-regular near $\partial\Omega$, but since for $g \in H^s_\alpha(\Omega \setminus \Sigma, \CC^N)$ and $\chi\in C^\infty(\overline \Omega)$ with $\chi=0$ near $\Sigma$ one always has $\chi g \in H^s_\alpha(\Omega)$, no statement about the regularity of $f \in H^{s,\partial\Omega}_\alpha(\Omega\setminus \Sigma,\CC^N)$ near $\Sigma$ can be made. Making use of the traces $\gamma^{\Sigma_p}_\pm$ at $\Sigma_p$ in 
\eqref{gammap} it is clear that there
exist unique continuous trace operators 
\begin{equation}\label{gamma111}
\gamma^\Sigma_\pm: H^{s,\Sigma}_\alpha(\Omega\setminus\Sigma,\CC^N)\mapsto H^{s-\half}(\Sigma,\CC^N),
\quad \gamma^\Sigma_\pm f\vert_{\Sigma_p}=\gamma^{\Sigma_p}_\pm (\chi_p f_{p\pm}),\quad p\in\{1,\dots, P\},
\end{equation}
where $\chi_p\in C^\infty_c(\Omega)$ with $\chi_p=1$ near $\Sigma_p$ and $\chi_p=0$ near 
$\bigcup_{j\ne p}\Sigma_j$.
 On the Lipschitz boundary $\partial\Omega$ we choose $s=\half$ and observe that there is a unique 
 continuous trace operator
\begin{equation}\label{gamma222}
\gamma^{\partial\Omega}: H^{\half,\partial\Omega}_\alpha(\Omega\setminus\Sigma,\CC^N)\mapsto L^2(\partial\Omega,\CC^N)
\end{equation}
that acts on $f\in H^{\half,\partial\Omega}_\alpha(\Omega\setminus\Sigma,\CC^N)$ in the same way as 
the trace operator in Lemma~\ref{lem1} when applied to $\chi_{\partial\Omega} f$,
where $\chi_{\partial\Omega}\in C^\infty(\overline\Omega)$ with $\chi_{\partial\Omega}=1$ near $\partial\Omega$ and $\chi_{\partial\Omega}=0$ near $\Sigma$ (so that $\chi_{\partial\Omega} f\in H^{\half}_\alpha(\Omega,\CC^N)$); for the trace operator in \eqref{gamma222} we use the same symbol as in Lemma~\ref{lem1} 
and 
it will be clear from the context on which space the trace operator $\gamma^{\partial\Omega}$ is considered.

Sometimes we will also work with locally $H^1$-regular functions in $\Omega\setminus\Sigma$
which are $H^s$-regular near $\Sigma_p$ for some 
$p\in\{1,\dots,P\}$ and $s \in [0, 1]$. Thus, we consider the spaces 
\begin{equation} \label{def_H_s_Sigma_p}
 \begin{split}
  H^{s,\Sigma_p}_\alpha(\Omega\setminus\Sigma,\CC^N):=\Big\{f\in L^2(\Omega,\CC^N)&:\chi f \in H^s_\alpha(\Omega\setminus\Sigma_p,\CC^N)\\ &\text{ for any }\chi\in C^\infty_c\big(\Omega\setminus\bigcup_{j\ne p}\Sigma_j\big)\Big\}.
 \end{split}
\end{equation}
Note that the functions in $H^{s,\Sigma_p}_\alpha(\Omega\setminus\Sigma,\CC^N)$	 
 are locally $H^1$-regular in $\Omega\setminus\Sigma$
and $H^s$-regular near $\Sigma_p$, but since for $g \in H^s_\alpha(\Omega \setminus \Sigma, \CC^N)$ and $\chi\in C^\infty_c(\Omega\setminus\bigcup_{j\ne p}\Sigma_j)$ one always has $\chi g \in H^s_\alpha(\Omega \setminus \Sigma_p,\CC^N)$, no statement about the regularity of $f \in H^{s,\Sigma_p}_\alpha(\Omega\setminus \Sigma,\CC^N)$ near $\bigcup_{j\ne p}\Sigma_j$  (and $\partial \Omega$) can be made.
Observe that 
\begin{equation*}
 H^{s,\Sigma}_\alpha(\Omega\setminus\Sigma,\CC^N)=\bigcap_{p=1}^P H^{s,\Sigma_p}_\alpha(\Omega\setminus\Sigma,\CC^N).
\end{equation*}
For each  $p\in\{1,\dots,P\}$ there
exist unique continuous trace operators 
\begin{equation}\label{gammap2}
\gamma^{\Sigma_p}_\pm: H^{s,\Sigma_p}_\alpha(\Omega\setminus\Sigma,\CC^N)\mapsto H^{s-\half}(\Sigma_p,\CC^N),
\end{equation}
where we use the same symbol as in \eqref{gammap} to avoid further technicalities. We note for completeness that the action of \eqref{gammap2} on 
$f\in H^{s,\Sigma_p}_\alpha(\Omega\setminus\Sigma,\CC^N)$
is defined as the action of 
the trace operator in \eqref{gammap} applied to $\chi_p f$, where $\chi_p\in C^\infty_c(\Omega)$ with $\chi_p=1$ near $\Sigma_p$ and $\chi_p=0$ 
near $\partial \Omega$ and $\bigcup_{j\ne p}\Sigma_j$ (so that $\chi_p f\in H^s_\alpha(\Omega\setminus\Sigma_p,\CC^N)$). It is also clear that the trace operator in \eqref{gamma111} applied to functions 
in $H^{s,\Sigma}_\alpha(\Omega\setminus\Sigma,\CC^N)$ and restricted to $\Sigma_p$ coincides with \eqref{gammap2}.

The next lemma is essentially a consequence of \cite[Lemma~2.4]{BFSV} for $n=2$ and 
\cite[Proposition~2.16]{OV} for $n=3$.

\begin{lemma}\label{lemreg}
Let $p\in\{1,\dots,P\}$ and $f\in H^{0,\Sigma_p}_\alpha(\Omega\setminus\Sigma,\CC^N)$. 
If $\gamma^{\Sigma_p}_\pm f \in H^\half(\Sigma_p,\CC^N)$, then $f\in  H^{1,\Sigma_p}_\alpha(\Omega\setminus\Sigma,\CC^N)$.
\end{lemma}

Finally, we introduce one more trace operator that acts on the domain of the self-adjoint Dirac operator $A$ in \eqref{def_A_mu},
\begin{equation}\label{gamma-sigma}
\gamma^\Sigma: \dom A= \big\{ f \in  H_{\alpha}^{\half}(\Omega, \CC^N): \cB \gamma^{\partial \Omega} f = 0  \big\} \rightarrow H^\half(\Sigma,\CC^N),
\end{equation}
which is well defined since $\dom A\subset H^1_\loc(\Omega, \CC^N)$; cf. Lemma~\ref{lem2}.

\begin{lemma}\label{lem9}
	The trace operator \eqref{gamma-sigma} is surjective and bounded with respect to the graph norm of $A$.
\end{lemma}

\begin{proof}
Let $U$ be a smooth open neighborhood of $\Sigma$ with $\overline{U}\subset \Omega$.
Then $U$ is an extension domain, and each function in $H^1(U, \CC^N)$ can be extended to a function in $H^1_0(\Omega,\CC^N)\subset\dom A$. It follows that the map
\[
\dom A\ni f\mapsto f|_U\in H^1(U,\CC^N)
\]
is surjective. Furthermore, the map
\[
H^1(U,\CC^N)\ni f\mapsto f|_\Sigma\in H^\half(\Sigma,\CC^N)
\]
is surjective due to the standard Sobolev trace theorem.
Hence, the trace operator $\gamma^\Sigma$ in \eqref{gamma-sigma} is surjective being the composition of these two maps.

In order to show the boundedness let us pick a function $\chi\in C^\infty_c(\Omega)$ with $\chi=1$ on $U$, then for any
$f\in \dom A$ one has $\chi f\in H^1_0(\Omega,\CC^N)$ and hence the $0$-extension onto $\RR^n$ satisfies $\widetilde{\chi f}\in H^1 (\RR^n,\CC^N)$. Due to the Sobolev trace theorem there is a $c_1>0$
such that
\[
\|\gamma^\Sigma g\|_{H^\half(\Sigma,\CC^N)}\le c_1 \|g\|_{H^1(U,\CC^N)} \text{ for all }g\in H^1(U,\CC^N).
\]
The uniform ellipticity of the Dirac operator in $\RR^n$ implies that there is a $c_2>0$ such that
for any $h\in H^1(\RR^n,\CC^N)$ one has
\[
\|h\|_{H^1(\RR^n,\CC^N)}\le c_2 \Big( \|h\|_{L^2(\RR^n,\CC^N)}+\|D h\|_{L^2(\RR^n,\CC^N)}\Big).
\]	
Then for any $f\in \dom A$ one has
\begin{align*}
	\|\gamma^\Sigma f\|_{H^\half(\Sigma,\CC^N)}&\le c_1\big\|f|_U\|_{H^1(U,\CC^N)}
	\le c_1 \|\widetilde{\chi f}\|_{H^1(\RR^n,\CC^N)}\\
	&\le c_1c_2 \Big( \|\widetilde{\chi f}\|_{L^2(\RR^n,\CC^N)}+\|D (\widetilde{\chi f})\|_{L^2(\RR^n,\CC^N)}\Big)\\
	&=c_1 c_2 \Big( \|\chi f\|_{L^2(\Omega,\CC^N)}+\|D (\chi f)\|_{L^2(\Omega,\CC^N)}\Big)\\
	&=c_1 c_2 \Big( \|\chi f\|_{L^2(\Omega,\CC^N)}+\big\|\chi D f -\rmi (\alpha\cdot\nabla \chi) f\big\|_{L^2(\Omega,\CC^N)}\Big)\\
	&\le c_1 c_2 c_3\Big( \|f\|_{L^2(\Omega,\CC^N)}+ \|D f\|_{L^2(\Omega,\CC^N)}\Big)
\end{align*}
for $c_3:=\|\chi\|_\infty +\|\alpha\cdot\nabla\chi\|_\infty$,
which gives the sought claim.
\end{proof}

\subsection{The operator $A^{\eta,\tau}$}

Let $\Sigma$ and $\Sigma_p$, $p\in\{1,\dots,P\}$, be as in the beginning of the previous subsection, and let
$\eta,\tau:\Sigma\to\RR$ be functions that are constant on each $\Sigma_p$. We denote
\[
\eta_p:=\eta|_{\Sigma_p}\in\RR\quad\text{and}\quad \tau_p:=\tau|_{\Sigma_p}\in\RR,\quad p\in\{1,\dots,P\}.
\]
Our next goal is to rigorously define a Dirac operator in $\Omega$ with the same elliptic boundary condition $\cB\gamma^{\partial\Omega}f=0$ 
on $\partial\Omega$ as $A$ in \eqref{def_A_mu},
an electrostatic $\delta$-shell interaction of strength $\eta$ and a Lorentz scalar $\delta$-shell interaction of strength $\tau$, i.e. 
a Dirac operator acting as
\[
f\mapsto D f \text{ in } \cD'(\Omega\setminus\Sigma,\CC^N)\cap L^2(\Omega,\CC^N)
\]
on the functions $f$ satisfying the same boundary condition on $\partial\Omega$ as $A$ 
and the transmission condition
\begin{gather*}
\rmi (\alpha\cdot\nu^\Sigma)(f_{p+}-f_{p-})+\dfrac{1}{2}(\eta_p I_N+\tau_p\beta)(f_{p+} + f_{p-})=0 \text{ on each }\Sigma_p.
\end{gather*}

\begin{definition} \label{adefjussi}
The operator $A^{\eta,\tau}$ in $L^2(\Omega,\CC^N)$ is defined as 
\begin{equation}\label{dom-pert}
 \begin{split}
  A^{\eta,\tau} f &= D f ,\\
  \dom A^{\eta,\tau} &=\Big\{
	f\in H^{0,\Sigma}_\alpha(\Omega\setminus\Sigma,\CC^N)\cap H^{\half,\partial\Omega}_\alpha(\Omega\setminus\Sigma,\CC^N): \cB\gamma^{\partial\Omega}f=0\text{ on }\partial\Omega,\\
	 & \qquad\qquad\qquad \rmi (\alpha\cdot\nu^\Sigma)\big(\gamma^\Sigma_+f-\gamma^\Sigma_-f\big)+\dfrac{1}{2}(\eta I_N+\tau\beta)(\gamma^\Sigma_+f + \gamma^\Sigma_-f)=0 \text{ on }\Sigma\Big\},
 \end{split}
\end{equation}
where $\gamma^\Sigma_\pm$ and $\gamma^{\partial\Omega}$ are the trace operators in \eqref{gamma111} (for $s=0$) and \eqref{gamma222}, respectively.
\end{definition}

It will turn out in Theorem~\ref{thm17} with the help of extension theory methods for symmetric operators (see Appendix~\ref{sec-bt}) that the Dirac operator 
$A^{\eta,\tau}$ is self-adjoint in $L^2(\Omega,\CC^N)$.
For these considerations it is useful to introduce the minimal and maximal Dirac operators defined by
\begin{equation*}
 A_{min} f=D f,\quad \dom A_{min}=H^1_0(\Omega\setminus\Sigma,\CC^N),
\end{equation*}
and
\begin{equation*}
 A_{max} f=D f,\quad \dom A_{max}=H^0_\alpha(\Omega\setminus\Sigma,\CC^N).
\end{equation*}
It can be shown in the same way as in \cite[Proposition~3.1]{BH20} that $A_{min}$ is a densely defined closed symmetric operator  in $L^2(\Omega,\CC^N)$ such that
$A_{min}^*=A_{max}$ holds. Besides the minimal and maximal operator we shall make use of the restriction $S$ of the self-adjoint Dirac operator $A$ onto
the kernel of the trace map $\gamma^\Sigma$ in \eqref{gamma-sigma}, which is given by
\begin{equation}\label{smjaja}
S f=D f,\quad \dom S=\big\{f\in \dom A:\ \gamma^\Sigma f=0\big\}=\dom A\vert_{\ker\gamma^\Sigma}.
\end{equation}
It follows from Lemma~\ref{lem9} that $S$ is a densely defined closed symmetric operator in $L^2(\Omega,\CC^N)$, and it is also clear that the following 
inclusions hold
\begin{equation}\label{subsets}
 A_{min}\subset S\subset A=A^* \subset S^*\subset A_{min}^*=A_{max}.
\end{equation}
Furthermore, it is not difficult to see that also the inclusion $S\subset A^{\eta,\tau}$ holds. In the next lemma we specify the domain of the operator
$S^*$.

\begin{lemma}\label{lem10}
	The adjoint of the operator $S$ in \eqref{smjaja} is given by
	\begin{equation}\label{star1}
 \begin{split}
  S^* f &= D f ,\\
	\dom S^*&=\Big\{
	f\in H^{0,\Sigma}_\alpha(\Omega\setminus\Sigma,\CC^N)\cap H^{\half,\partial\Omega}_\alpha(\Omega\setminus\Sigma,\CC^N): \cB\gamma^{\partial\Omega}f=0\text{ on }\partial\Omega\Big\},
\end{split}
\end{equation}
and, in particular, the operator $A^{\eta,\tau}$ in Definition~\ref{adefjussi} is a restriction of $S^*$.
\end{lemma}

\begin{proof}
It is clear from \eqref{subsets} that $S^*$ acts as $D$ and that $\dom S^*$ is contained in $\dom A_{max}=H^0_\alpha(\Omega\setminus\Sigma,\CC^N)$.

Let us first show the inclusion $\subset$ for the operator domain in \eqref{star1}. Consider an arbitrary $f\in \dom S^*$. 
We shall often use that for $\chi\in C^\infty( \overline{\Omega})$ one has
\begin{equation}\label{malsomalso}
 D (\chi f)=\chi D f-\rmi (\alpha\cdot\nabla\chi) f\in \cD'(\Omega\setminus\Sigma,\CC^N)\cap L^2(\Omega,\CC^N).
\end{equation}
In particular, choosing $\chi\in C_c^\infty( \Omega)$ in~\eqref{malsomalso}
shows $\chi f\in  H^0_\alpha(\Omega\setminus\Sigma,\CC^N)$ and hence 
 $f\in H^{0,\Sigma}_\alpha(\Omega\setminus\Sigma,\CC^N)$.
Now let $\chi\in C^\infty(\overline \Omega)$ with $\chi=0$ near $\Sigma$. Using \eqref{malsomalso} we compute for $g\in \dom A$
\begin{equation}
	\label{tmp001}
\begin{aligned}
	\big\langle D(\chi f),g\big\rangle_{L^2(\Omega,\CC^N)}&=\big\langle \chi D f-\rmi (\alpha\cdot\nabla\chi)f,g\big\rangle_{L^2(\Omega,\CC^N)}\\
	&=\langle  D f,\overline\chi g \rangle_{L^2(\Omega,\CC^N)}+\big\langle f, \rmi (\alpha\cdot\nabla\overline \chi)g\big\rangle_{L^2(\Omega,\CC^N)}.
\end{aligned}
\end{equation}
We further have  $\overline \chi g\in \dom S$, hence,
\begin{align*}
	\langle  D f,\overline\chi g \rangle_{L^2(\Omega,\CC^N)}&=\langle  S^* f,\overline\chi g \rangle_{L^2(\Omega,\CC^N)}
	=\big\langle  f,S (\overline\chi g) \big \rangle_{L^2(\Omega,\CC^N)}\\
	&=\big\langle  f,D (\overline\chi g) \big \rangle_{L^2(\Omega,\CC^N)}
=\big\langle  f,\overline\chi D g-\rmi (\alpha\cdot\nabla\overline\chi) g \big \rangle_{L^2(\Omega,\CC^N)}\\
&=\big\langle \chi f, A g\big \rangle_{L^2(\Omega,\CC^N)}
-\big\langle f, \rmi (\alpha\cdot\nabla\overline \chi)g\big\rangle_{L^2(\Omega,\CC^N)},
\end{align*}
and the substitution into \eqref{tmp001} yields
\[
	\big\langle D(\chi f),g\big\rangle_{L^2(\Omega,\CC^N)}=\big\langle \chi f, A g\big \rangle_{L^2(\Omega,\CC^N)}.
\]
As $g\in \dom A$ was arbitrary, this means $\chi f\in \dom A^*=\dom A$, i.e. 
\[
\chi f\in H^\half_\alpha(\Omega,\CC^N) \quad\text{and}\quad \cB \gamma^{\partial\Omega} (\chi f)=0 \text{ on }\partial\Omega.
\]
Since this holds for any $\chi\in C^\infty( \overline{\Omega})$ with $\chi = 0$ near $\Sigma$, 
one obtains $f\in H^{\half,\partial\Omega}_\alpha(\Omega\setminus\Sigma,\CC^N)$; cf. \eqref{ueblerraum}. 
The identity $\cB \gamma^{\partial\Omega} (\chi f)=0$ with $\chi\in C^\infty(\overline \Omega)$ such that $\chi=1$ near $\partial\Omega$ gives the required boundary condition for $f$ on $\partial\Omega$. The inclusion $\subset$ in \eqref{star1}
is now proved.

It remains to show the inclusion $\supset$ in \eqref{star1}. Let $f$ be an arbitrary function from the set on the right-hand side of \eqref{star1}. Let $\chi\in C^\infty_c(\Omega)$ with $\chi=1$ near $\Sigma$. 
We represent $$f=\chi f +(1-\chi) f$$ and show that both summands are in $\dom S^*$. 

Note first that $(1-\chi)f\in H^\half_\alpha(\Omega,\CC^N)$ 
and $(1-\chi)f=f$ near $\partial\Omega$, and thus
we have $(1-\chi)f\in \dom A\subset\dom S^*$. Hence, it suffices  to show $\chi f\in\dom S^*$.
For this fix some $\rho \in C^\infty_c(\Omega)$ with $\rho=1$ on an open neighborhood of $\supp \chi$.
Then, due to the locality of $D$ for any $g\in\dom S$ one obtains
\[
\langle \chi f, S g\rangle_{L^2(\Omega,\CC^N)}=\langle \chi f, D g\rangle_{L^2(\Omega,\CC^N)}=\big\langle \chi f, D (\rho g)\big\rangle_{L^2(\Omega,\CC^N)}.
\]
We have $\rho g\in H^1_0(\Omega\setminus\Sigma,\CC^N)=\dom A_{min}$. 
Due to $\chi f\in H^0_\alpha(\Omega\setminus\Sigma,\CC^N)=\dom A_{max}$ and $A_{max} = A_{min}^*$ we have   
\begin{align*}
	\big\langle \chi f, D (\rho g)\big\rangle_{L^2(\Omega,\CC^N)}&= \langle \chi f, A_{min} (\rho g)\rangle_{L^2(\Omega,\CC^N)}\\
	&=\langle A_{max}(\chi f), \rho g\rangle_{L^2(\Omega,\CC^N)}\\
	&=\big \langle D(\chi f), \rho g\rangle_{L^2(\Omega,\CC^N)}\\
	&=\big \langle D(\chi f), g\rangle_{L^2(\Omega,\CC^N)},
\end{align*}
which yields $\langle \chi f, S g\rangle_{L^2(\Omega,\CC^N)}=\langle D (\chi f),  g\rangle_{L^2(\Omega,\CC^N)}$ for any $g\in\dom S$
and shows $\chi f\in \dom S^*$. This completes the proof.
\end{proof}

\subsection{Definition using boundary triples}

In this subsection we use the abstract concept of boundary triples from Appendix~\ref{sec-bt}
to study the operator $A^{\eta,\tau}$ in Definition~\ref{adefjussi}. More precisely,
we use the construction of a boundary triple from Theorem~\ref{thm13} by starting with
the operator $Q:=A$ and its symmetric restriction $S$.
Then the operator $A^{\eta,\tau}$ can be represented as in Theorem~\ref{thm12} .
The following considerations can be viewed
as a generalization of the treatment in \cite{BHOP20} and, in particular,
we refer the reader there for a more detailed exposition for the case $n=2$.

Recall that $\Sigma_p$, $p\in\{1,\dots,P\}$, denote the connected components of the compact $C^ \infty$-smooth hypersurface $\Sigma$ on which the transmission conditions are imposed. For each $p$ we choose a classical elliptic (scalar) pseudodifferential operator $\Lambda_p$ of order $\half$ on $\Sigma_p$, and we additionally require $\Lambda_p$ to be injective and symmetric in $L^2(\Sigma_p)$; a specific choice of $\Lambda_p$ will be made later. Then
\begin{equation}\label{nimmdas}
\Lambda:=\bigoplus_{p=1}^P\Lambda_p:\ H^s(\Sigma)\to H^{s-\half}(\Sigma)
\end{equation}
is an isomorphism for any $s\in\RR$, and we denote
its natural extension $\Lambda\otimes I_N$ on $\CC^N$-valued functions by the same symbol.
Due to Lemma~\ref{lem9}
the map
\[
\gamma:=\Lambda \gamma^\Sigma :\ \dom A\to L^2(\Sigma,\CC^N)
\]
is surjective and bounded with respect to the graph norm.
Since $\ker\gamma=\ker\gamma^\Sigma$ it is clear from
\eqref{smjaja} that
\begin{equation}\label{smjaja2}
S f=D f,\quad \dom S=\big\{f\in \dom A:\ \gamma f=0\big\}=\dom A\vert_{\ker\gamma}.
\end{equation}
Following \eqref{eq-gz1a} we now define
\[
\Gamma_z:=\big(\gamma(A-\Bar z)^{-1}\big)^*,\quad z\in\res A.
\]
In order to determine $\Gamma_z$ more explicitly consider first
$\varphi\in C^\infty(\Sigma,\CC^N)$ and $f\in C^\infty_c(\Omega,\CC^N)$, and note that
\begin{align*}
    \big[\gamma^\Sigma (A-\Bar z)^{-1}f\big](x)=\int_\Omega G_{\Bar z}(x,y)f(y)\dd y \,\,\text{ for any }x\in \Sigma
\end{align*}
by virtue of Theorem~\ref{theo1}. We have $\Lambda \varphi\in C^\infty (\Sigma,\CC^N)$, and due to the properties of the integral kernel $G_z$ stated in Theorem~\ref{theo1}, the function
\[
\Sigma\times \Omega\ni (x,y)\mapsto G_{\Bar z}(x,y)\in \CC^{N\times N}
\]
is integrable on $\Sigma\times \supp f$. Using \eqref{mk2} and Fubini's theorem we compute
\begin{align*}
    \big\langle \varphi, \gamma (A-\Bar z)^{-1}f\big\rangle_{L^2(\Sigma,\CC^N)}&
    =\big\langle \Lambda\varphi, \gamma^\Sigma (A-\Bar z)^{-1}f\big\rangle_{L^2(\Sigma,\CC^N)}\\
    &=\int_\Sigma \bigg\langle(\Lambda\varphi)(x),\int_\Omega G_{\Bar z}(x,y)f(y)\,\dd y\bigg\rangle_{\CC^N}\dd s(x)\\
    &=\int_\Sigma \int_\Omega \big\langle G_{\Bar z}(x,y)^*(\Lambda\varphi)(x),f(y)\Big\rangle_{\CC^N}\dd y\,\dd s(x)\\
    &=\int_\Omega\int_\Sigma \big\langle G_{z}(y,x)(\Lambda\varphi)(x),f(y)\big\rangle_{\CC^N}\dd s(x)\,\dd y\\
    &=\int_\Omega \bigg\langle\int_\Sigma G_{z}(y,x) (\Lambda\varphi)(x)\,\dd s(x),f(y)\bigg\rangle_{\CC^N}\dd y\\
    &=\langle \Phi_z \Lambda\varphi, f\rangle_{L^2(\Omega,\CC^N)},
\end{align*}
where the mapping $\Phi_z$ is defined by
\begin{equation}\label{phi1}
\Phi_z \psi(x):=\int_\Sigma G_z(x,y)\psi(y)\,\dd s(y),\quad \psi\in C^\infty(\Sigma,\CC^N),\quad x\in \Omega \setminus \Sigma.
\end{equation}
As $C^\infty_c(\Omega,\CC^N)$ is dense in $L^2(\Omega,\CC^N)$, the resulting identity
\[
\big\langle \varphi, \gamma (A-\Bar z)^{-1}f\big\rangle_{L^2(\Sigma,\CC^N)}=\langle \Phi_z \Lambda\varphi, f\rangle_{L^2(\Omega,\CC^N)}
\]
holds for all $\varphi\in C^\infty(\Sigma,\CC^N)$ and $f\in L^2(\Omega,\CC^N)$, and hence we conclude
\begin{equation}\label{gamma111z}
 \Gamma_z\varphi =\Phi_z \Lambda\varphi,\quad \varphi\in C^\infty(\Sigma,\CC^N).
\end{equation}
By the abstract theory (see Appendix~\ref{sec-bt}), the map $\Gamma_z$ 
extends uniquely to a bounded linear operator $L^2(\Sigma,\CC^N)\to \ker(S^*-z)\subset \dom S^*$
equipped with the graph norm. Using the fact that $\Lambda:L^2(\Sigma,\CC^N)\to H^{-\half}(\Sigma,\CC^N)$
is an isomorphism, we conclude that the mapping $\Phi_z$ in \eqref{phi1} defined on $C^\infty(\Sigma,\CC^N)$  extends uniquely to a bounded linear operator
\[
\Phi_z: H^{-\half}(\Sigma,\CC^N)\to L^2(\Omega,\CC^N),\qquad  \ran \Phi_z\subset\ker (S^*-z).
\]
As the $L^2$-norm and the graph norm of $S^*$ are equivalent on $\ker (S^*-z)$
we can also view $\Phi_z$ as a bounded linear operator in the context
$\Phi_z: H^{-\half}(\Sigma,\CC^N)\to H^0_\alpha(\Omega\setminus\Sigma,\CC^N)
$,
where $\dom S^*\subset\dom A_{max}=H^0_\alpha(\Omega\setminus\Sigma,\CC^N)$; cf. \eqref{subsets}.
We conclude that the identity \eqref{gamma111z} extends to
\begin{equation}\label{eq-gbz}
 \Gamma_z\varphi =\Phi_z \Lambda\varphi,\quad \varphi\in L^2(\Sigma,\CC^N).
\end{equation}
Similarly, considering $\Phi_z$ as a bounded operator from $L^2(\Sigma,\CC^N)$ to $L^2(\Omega,\CC^N)$
we conclude that its adjoint $\Phi_z^*$ given by
\[
\Phi_z^* f(x)=\int_\Omega G_{\Bar z}(x,y) f(y)\,\dd y,\quad f\in C^\infty_c(\Omega,\CC^N),\quad x\in \Sigma,
\]
extends uniquely to a bounded linear operator
\[
\Phi_z^*: L^2(\Omega,\CC^N)\to H^\half(\Sigma,\CC^N).
\]

Following the construction in Appendix~\ref{sec-bt} we note that for
$z\in\res A$ each function $f\in\dom S^*$ can be uniquely represented in the form
$f=f_z+\Phi_z\Lambda \xi$ with $f_z\in \dom A$ and $\xi \in L^2(\Sigma,\CC^N)$
(where $\xi$ does not depend on $z$).
In order to construct a boundary triple with the help of Theorem~\ref{thm13}
it is necessary to
determine $\xi$ directly from $f$. This will be based on the fact that the function $\Phi_{z}\Lambda\xi$
has a jump along $\Sigma$, while $f_z$ has no jump along $\Sigma$ (as
$f_z\in\dom A\subset H^1_\loc(\Omega,\CC^N)$).

\begin{lemma}\label{lem14} For $z\in\res A$ the following is true:  
 \begin{itemize}
  \item[(a)] 	The singular integral operator  $C_{z}: C^\infty(\Sigma,\CC^N)\to C^\infty(\Sigma,\CC^N)$,
	\begin{equation} \label{def_C_z}
	C_{z}\varphi (x)=\pv \int_\Sigma G_{z}(x,y)\varphi(y)\,\dd s(y),\quad x\in\Sigma, \quad \varphi\in C^\infty(\Sigma,\CC^N),
	\end{equation}
	where $\pv$ stands for the Cauchy principal value, extends uniquely to a bounded linear operator in $H^r(\Sigma,\CC^N)$ for any $r\in\RR$.
\item[(b)] For any
	$\varphi\in H^{-\half}(\Sigma,\CC^N)$ one has
	\begin{equation} \label{Plemelj_Sokhotskii}
	\gamma^\Sigma_\pm \Phi_{z} \varphi=\mp\dfrac{\rmi}{2}(\alpha\cdot\nu^\Sigma) \varphi + C_{z}\varphi.
	\end{equation}
\end{itemize}
\end{lemma}

In Section~\ref{section_separation} we will analyze the operator $C_z$ further with pseudodifferential techniques, see Corollary~\ref{corollary_C_z_pseudo}. For the construction of the boundary triple used to define $A^{\eta, \tau}$ the results of the preceding lemma are sufficient.

\begin{proof}[Proof of Lemma~\ref{lem14}]
 	(a) Choose $R > 0$ such that $\Omega \subset B_R(0)$ and $\chi \in C_c^\infty(\mathbb{R}^n)$ such that $\chi \equiv 1$ on $B_{2R}(0)$. Note that for the function $\cH_{z}$ in Lemma~\ref{lemma_expansion1} one has $\cH_{z}- \cH_{0}=z \cE_0\otimes I_N$ with $\cE_0$ defined by~\eqref{fundamental_solution_Laplace}. Hence,
  \begin{equation*}
    \begin{split}
      E_{z}(x) - E_0(x) 
      &= E_{z}(x) - \cH_{z}(x) - (E_{0}(x) - \cH_{0}(x))  + z \cE_0(x) \otimes I_N,\quad x\in\RR^n\setminus\{0\},
    \end{split}
  \end{equation*}
  where $E_0$ is defined by~\eqref{fundamentla_solution_Dirac_2D} \& \eqref{fundamentla_solution_Dirac_3D}.
Therefore, using Theorem~\ref{theo1} and the function $\mathcal{K}_z$ from Lemma~\ref{lemma_expansion1}  one can write for $x,y \in \Sigma$, $x\not= y$,
\begin{equation} \label{decomposition_kernel}
  \begin{split}
  G_z(x,y) &= E_0(x-y) + \chi(x-y)(E_{z}(x-y) - E_0(x-y)) + T_z(x,y)   \\
      &= E_0(x-y) + \mathcal{K}_z(x-y) - \mathcal{K}_0(x-y)  + z \cE_0(x-y) \otimes I_N + T_z(x,y).
  \end{split}
\end{equation}
Based on the above decomposition of $G_z$ we consider the 
operator
  \begin{equation} \label{def_Cauchy_op_full_space}
  \widetilde{C} \varphi(x) = \pv \int_\Sigma E_{0}(x-y) \varphi(y) \dd s(y), \quad x \in \Sigma, \quad \varphi \in C^\infty(\Sigma, \CC^N),
\end{equation}
and the integral operators acting on $\varphi \in C^\infty(\Sigma, \CC^N)$ evaluated at $x \in \Sigma$ as
\begin{equation} \label{def_int_ops}
  \begin{split}
  K_{\Sigma,z} \varphi(x) &= \int_\Sigma \mathcal{K}_z(x-y) \varphi(y) \dd s(y),\\
  K_{\Sigma,0} \varphi(x) &= \int_\Sigma \mathcal{K}_0(x-y) \varphi(y) \dd s(y),\\
  S_\Sigma \varphi(x) &= \int_\Sigma \cE_0(x-y) \varphi(y) \dd s(y), \\
  \widetilde{T}_{\Sigma,z}  \varphi(x) &= \int_\Sigma T_z(x,y) \varphi(y) \dd s(y).
  \end{split}
\end{equation}
It is known that $\widetilde{C}$ extends to a bounded linear operator in $H^r(\Sigma,\CC^N)$ for any $r\in\RR$; for $n=2$ this is stated in \cite[Proposition~3.3]{BHOP20}, while for $n=3$ one has by \cite[Theorem~4.1]{BBZ} that $\widetilde{C}$ is an $h$-pseudodifferential operator of order $0$ in the sense of \cite[Definition~2.1]{BBZ} and hence, by \cite[Remark after Theorem~14.2]{zworski} this implies the claimed continuity.  
Next, consider the map 
\begin{equation*}
  \widetilde{K}_z \varphi(x) := \int_{\mathbb{R}^n} \mathcal{K}_z(x-y) \varphi(y) \dd y, \qquad x \in \mathbb{R}^n, \quad \varphi \in C_c^\infty(\mathbb{R}^n).
\end{equation*}
In the following, we use the notations for classical pseudodifferential operators introduced in Appendix \ref{appendix_pseudos}.
Taking the expansion of $\mathcal{K}_z$ from Lemma~\ref{lemma_expansion1} into account, one can directly apply \cite[Chapter~7, Proposition~2.8]{T1} which shows that $\widetilde{K}_z \in \Op \cS^{-3}_N(\mathbb{R}^n)$. Hence, following the standard construction of restricting an integral operator in $\mathbb{R}^n$ to $\Sigma$ one gets
$K_{\Sigma,z}  \in \Op \cS^{-2}_N(\Sigma)$, see the considerations leading to \cite[Chapter~7, equation~(11.4)]{T1}; 
in particular, the maps $K_{\Sigma,z}$ can be extended to a bounded operator $K_{\Sigma,z}: H^r(\Sigma, \CC^N) \rightarrow H^{r+2}(\Sigma, \CC^N)$ for any $r\in\RR$ and thus, 
$K_{\Sigma,z}$ gives rise to a bounded operator in $H^r(\Sigma, \CC^N)$. In the same way one shows that $K_{\Sigma,0}  \in \Op \cS^{-2}_N(\Sigma)$ and that $K_{\Sigma,0}$ can be extended to a bounded operator
in $H^r(\Sigma, \CC^N)$ for any $r\in\RR$.
Moreover, it is known from \cite[Chapter~7, Proposition~11.2]{T1} that $S_\Sigma\in \Op \cS^{-1}_1(\Sigma)$ and thus, it can also be extended to a bounded operator in $H^r(\Sigma, \CC^N)$, $r\in\RR$. Finally, as $T_z$ is smooth (see Theorem~\ref{theo1}), one has $\widetilde{T}_{\Sigma,z} \in \Op \cS^{-\infty}_N(\Sigma)$, implying that also $\widetilde{T}_{\Sigma,z}$ can be extended to a bounded operator in $H^r(\Sigma, \CC^N)$, $r\in\RR$. Thus, it follows from~\eqref{decomposition_kernel} and~\eqref{def_Cauchy_op_full_space} that
\begin{equation*}
  C_{z} = \widetilde{C} + K_{\Sigma,z} - K_{\Sigma,0} +z S_\Sigma \otimes I_N + \widetilde{T}_{\Sigma,z}
\end{equation*}
can be extended to a bounded operator in $H^r(\Sigma, \CC^N)$ for any $r\in\RR$. In particular, for any $\varphi\in C^\infty(\Sigma, \CC^N)$ one has $C_z \varphi\in \bigcap_{r\in \RR} H^r(\Sigma, \CC^N)=C^\infty(\Sigma, \CC^N)$ by the Sobolev embedding theorem, and thus $C_z$ maps $C^\infty(\Sigma, \CC^N)$ into itself.

(b) We show \eqref{Plemelj_Sokhotskii} for $\varphi \in C^\infty(\Sigma, \CC^N)$, the claim for $\varphi \in H^{-\half}(\Sigma, \CC^N)$ follows then by density and continuity. Again, the decomposition of $G_z$ in~\eqref{decomposition_kernel} and the operators defined in~\eqref{def_Cauchy_op_full_space} and \eqref{def_int_ops} will play an important role. Define the operator
\begin{equation*}
  \widetilde{\Phi} \varphi(x) = \int_\Sigma E_{0}(x-y) \varphi(y) \dd s(y), \quad x \in \Omega \setminus \Sigma, \quad \varphi \in C^\infty(\Sigma, \CC^N).
\end{equation*}
Then, by \cite[Proposition~4.4]{BHSLS24} one has 
 \begin{equation} \label{Plemelj}
   \gamma^\Sigma_\pm \widetilde{\Phi} \varphi=\mp\dfrac{\rmi}{2}(\alpha\cdot\nu^\Sigma) \varphi + \widetilde{C}\varphi.
 \end{equation}
 Next, consider the functions $f_1, f_2$ defined for $x \in \Omega \setminus \Sigma$ by
 \begin{equation*}
   \begin{split}
     f_1(x) &:= \int_\Sigma \bigl( \mathcal{K}_z(x-y) - \mathcal{K}_0(x-y)  + T_z(x,y) \bigr) \varphi(y) \dd s(y), \\
     f_2(x) &:= z \int_\Sigma \cE_0(x-y) \varphi(y) \dd s(y).
	\end{split}
 \end{equation*}
 Since $\mathcal{K}_z- \mathcal{K}_0$ is bounded on $\RR^n$ and continuous on $\RR^n\setminus\{0\}$ by  Lemma~\ref{lemma_expansion1}, since single points have surface measure zero on $\Sigma$, and since $T_z$ is smooth by Theorem~\ref{theo1}, one gets with dominated convergence that $f_1$ can be extended to a continuous function $\widetilde{f}_1$ in $\Omega$ satisfying
 \begin{equation*}
   \begin{split}
   \widetilde{f}_1(x) &= \int_\Sigma \bigl( \mathcal{K}_z(x-y) - \mathcal{K}_0(x-y)  + T_z(x,y) \bigr) \varphi(y) \dd s(y) \\
   &= \bigl( K_{\Sigma,z} \varphi - K_{\Sigma,0} \varphi + \widetilde{T}_{\Sigma,z} \varphi\bigr) (x), \qquad x \in \Sigma.
   \end{split}
 \end{equation*}
 Furthermore, by the gradient estimates in Lemma~\ref{lemma_expansion1}~(d) \& (e) one has, for $0<|x|\leq 2R$,
 \begin{align*}
  | \nabla(\mathcal{K}_z- \mathcal{K}_0)(x)|&\leq c\big(1+\big|\log|x|\big|\big), &&n=2,\\
  | \nabla(\mathcal{K}_z- \mathcal{K}_0)(x)|&\leq c|x|^{-1}, &&n=3.
  \end{align*}
  Since for any $r>0$ we have
 \begin{align*}
 \sup_{|x|<r}\int_{\Sigma}\big(1+\big|\log|x-y|\big|\big)\dd s(y)&<\infty, &&n=2, \\ \sup_{|x|<r}\int_{\Sigma}|x-y|^{-1}\dd s(y)&<\infty,&&n=3,
 \end{align*} and $T_z$ is smooth, the differentiation under the integral sign shows that $\nabla f_1$ is bounded on $\Omega_0\setminus\Sigma$ for any open set $\Omega_0$ with $\Sigma\subset\Omega_0$ and $\overline{\Omega_0}\subset\Omega$. Hence, $f_1$ belongs to $H^1$ on both sides of $\Sigma$ in a neighborhood of $\Sigma$, and, therefore, its trace $\gamma_\pm^\Sigma f_1$ coincide with $\widetilde{f}_1\vert_{\Sigma}$. 
 Moreover, by \cite[equation~(7.8)]{M00} one has
 \begin{equation*}
   \gamma^\Sigma_\pm f_2  =z (S_\Sigma \otimes I_N) \varphi.
 \end{equation*}
 Hence, we conclude
 \begin{equation*}
   \begin{split}
     \gamma^\Sigma_\pm \Phi_{z} \varphi &= \gamma^\Sigma_\pm \bigl( \widetilde{\Phi} \varphi + f_1 +f_2 \bigr) \\
     &= \mp\dfrac{\rmi}{2}(\alpha\cdot\nu^\Sigma) \varphi + \widetilde{C}\varphi 
     + K_{\Sigma,z} \varphi - K_{\Sigma,0} \varphi + \widetilde{T}_{\Sigma,z} \varphi +z (S_\Sigma \otimes I_N) \varphi \\
     &= \mp\dfrac{\rmi}{2}(\alpha\cdot\nu^\Sigma) \varphi+ C_{z}\varphi,
   \end{split}
 \end{equation*}
 where in the last step again~\eqref{decomposition_kernel} was used. This is the claimed formula.
\end{proof}

 Now we use Theorem~\ref{thm13} to construct an explicit boundary triple in the present situation.
The following assertion is the first step for the 
spectral study of $A^{\eta,\tau}$. 
We provide a short proof of Lemma~\ref{lem15} and refer the reader to 
\cite[Proof of Proposition~3.5]{BHOP20} for more details (in the case $n=2$, 
$\Sigma=\Sigma_1$ and $\Omega=\mathbb R^2$).

\begin{lemma}\label{lem15}
Fix some $\zeta\in \RR\cap\res A$, let $S^*$ be the adjoint of the densely defined closed symmetric operator $S=A|_{\ker\gamma}$ from Lemma~\ref{lem10}, and consider the 
mappings
\begin{equation*}
 \Pi_0,\Pi_1:\dom S^*\to L^2(\Sigma,\CC^N)
\end{equation*}
defined by
\begin{equation*}
 \begin{split}
  \Pi_0 f&:=\rmi \Lambda^{-1}(\alpha\cdot\nu^\Sigma)(\gamma^\Sigma_+ f-\gamma^\Sigma_- f),\\
	\Pi_1 f&:=\Lambda\bigg[\dfrac{1}{2}(\gamma^\Sigma_+ f + \gamma^\Sigma_- f)- C_\zeta \Lambda \Pi_0 f\bigg].
 \end{split}
\end{equation*}
Then $\big( L^2(\Sigma,\CC^N), \Pi_0, \Pi_1\big)$ is a boundary triple for $S^*$ such that $A=S^*|_{\ker\Pi_0}$.
For $z\in\res A$ the associated $\gamma$-field $\Gamma_z$ is given by \eqref{eq-gbz} and the associated Weyl function $M_z$ is given by
\begin{equation} \label{Weyl_function}
M_z:= \Lambda (C_z-C_\zeta)\Lambda.
\end{equation}
\end{lemma}

\begin{proof}
Essentially the assertions follow from Theorem~\ref{thm13} for the special case
$\zeta\in\RR$. More precisely, write $f\in \dom S^*$ in the form
\[
f=f_\zeta+\Gamma_\zeta \xi=f_\zeta+\Phi_\zeta\Lambda \xi,\quad f_\zeta\in \dom A,\quad \xi \in L^2(\Sigma,\CC^N),
\]
see \eqref{jabitte2} and \eqref{eq-gbz}.
Due to $f_\zeta\in \dom A\subset H^1_\loc(\Omega,\CC^N)$ we have $\gamma^\Sigma_+ f_\zeta-\gamma^\Sigma_-f_\zeta=0$. Therefore,
Lemma~\ref{lem14}(b) yields
\[
\gamma^\Sigma_+ f-\gamma^\Sigma_- f=-\rmi (\alpha\cdot\nu^\Sigma) \Lambda \xi,
\]
 and hence
\[ 
\xi=\rmi \Lambda^{-1}(\alpha\cdot\nu^\Sigma)(\gamma^\Sigma_+ f-\gamma^\Sigma_- f),
\]
which yields in the above expression for $\Pi_{0}$ according to Theorem~\ref{thm13}. Furthermore, as $\zeta=\overline{\zeta}$ Theorem~\ref{thm13}
gives
\begin{align*}
	\Pi_1 f&=\gamma f_\zeta=\Lambda \gamma^\Sigma f_\zeta=\dfrac{1}{2}\Lambda (\gamma^\Sigma_+ + \gamma^\Sigma_-)f_\zeta\\
	&=\dfrac{1}{2}\Lambda (\gamma^\Sigma_+ + \gamma^\Sigma_-)(f-\Phi_\zeta\Lambda\xi)
	=\dfrac{1}{2}\Lambda (\gamma^\Sigma_+ + \gamma^\Sigma_-)(f-\Phi_\zeta\Lambda\Pi_0 f)
\end{align*}
and another application of Lemma~\ref{lem14}(b) gives the required expression for $\Pi_1$. Now Theorem~\ref{thm13} implies that $\big( L^2(\Sigma,\CC^N), \Pi_0, \Pi_1\big)$ is a boundary triple for $S^*$ such that 
$A=S^*|_{\ker\Pi_0}$ and the associated $\gamma$-field $\Gamma_z$ is given by \eqref{eq-gbz}.

For the corresponding Weyl function $M_z$ we have
\[
M_z=\gamma(\Gamma_z-\Gamma_\zeta)=\dfrac{1}{2}\Lambda(\gamma^\Sigma_+ +\gamma^\Sigma_-)(\Phi_z-\Phi_\zeta)\Lambda,
\]
and the result follows from Lemma~\ref{lem14}(b).
\end{proof}

Using the above construction, we can finally represent the operator $A^{\eta,\tau}$ in Definition~\ref{adefjussi} via the boundary triple $\big( L^2(\Sigma,\CC^N), \Pi_0, \Pi_1\big)$  from Lemma~\ref{lem15}.

\begin{lemma}\label{lem16}
Assume that $|\eta_p|\ne |\tau_p|$ for all $p\in\{1,\dots,P\}$ and define
\begin{equation} \label{def_Theta}
\Theta^{\eta,\tau}:=-\Lambda\bigg[ \dfrac{1}{\eta^2-\tau^2}(\eta I_N -\tau\beta)+C_\zeta\bigg]\Lambda,
\end{equation}	
which is viewed as an operator on the maximal domain $\dom\Theta^{\eta,\tau}$ in $L^2(\Sigma,\CC^N)$.
Then
\begin{equation*} 
  \dom A^{\eta,\tau} = \big\{ f \in \dom S^*: \Pi_0 f\in \dom \Theta^{\eta,\tau},
\Pi_1 f= \Theta^{\eta,\tau}\Pi_0 f \big\},
\end{equation*}
i.e. the operator $A^{\eta,\tau}$ is the restriction of $S^*$ on the functions $f$ satisfying 
\[
\Pi_0 f\in \dom \Theta^{\eta,\tau}\quad\text{and}\quad
\Pi_1 f= \Theta^{\eta,\tau}\Pi_0 f.
\]
\end{lemma}	
\begin{proof}
	From the expression for the boundary maps $\Pi_0,\Pi_1$ for any $f\in\dom S^*$
	one obtains
	\[
	\gamma^\Sigma_+f-\gamma^\Sigma_-f=-\rmi(\alpha\cdot\nu^\Sigma)\Lambda\Pi_0 f
	\quad\text{and}\quad
	\frac{1}{2}
	(\gamma^\Sigma_+f+\gamma^\Sigma_-f)=\Lambda^{-1}\Pi_1f+C_\zeta\Lambda \Pi_0 f.
	\]
The direct substitution into the boundary condition at $\Sigma$ for $A^{\eta,\tau}$ in Definition~\ref{adefjussi} shows that $f\in \dom A^{\eta, \tau}$ if and only if
\[
\Lambda \Pi_0 f +(\eta I_N+\tau\beta)\bigl[\Lambda^{-1}\Pi_1+C_\zeta\Lambda \Pi_0 \bigr] f=0,
\]
which can be equivalently rewritten as
\[
(\eta I_N+\tau\beta)\Lambda^{-1}\Pi_1f=-\Big( I_N +(\eta I_N+\tau\beta)C_\zeta\Big)\Lambda \Pi_0f.
\]
The assumption $|\eta_p|\ne |\tau_p|$ for all $p\in\{1,\dots,P\}$ ensures that the
piecewise constant matrix function $\eta I_N+\tau\beta:\Sigma\rightarrow\CC^{N \times N}$ is invertible with
$(\eta I_N+\tau\beta)^{-1}=\tfrac{1}{\eta^2-\tau^2}(\eta I_N-\tau\beta)$,
which gives the required representation. 
\end{proof}

\section{Spectral properties of shell interactions}\label{sec4}

The main results of this section are Theorem~\ref{thm17} and Theorem~\ref{thm37}, where we show the self-adjointness and determine the essential spectrum of the operator $A^{\eta, \tau}$ in Definition~\ref{adefjussi}. Moreover, if the interaction is \textit{non-critical} on one part $\Sigma_p$, $p \in \{ 1, \dots, P \}$, i.e. $\eta_p^2 - \tau_p^2 \neq 4$, then it will turn out in Corollary~\ref{cor21} that the functions in $\dom A^{\eta, \tau}$ are $H^1$-regular near $\Sigma_p$. For both purposes, we will rewrite the parameter $\Theta^{\eta,\tau}$ given by~\eqref{def_Theta} as the sum of operators acting on $\Sigma_p$ and a perturbation term, followed by a more detailed analysis. In the next subsection, which is of preparatory nature, we first 
investigate the involved maps as pseudodifferential operators; see Appendix \ref{appendix_pseudos} for the associated notations and basic results.

\subsection{Separating connected components} \label{section_separation}

Let us make a specific choice of the operator $\Lambda_p \in \Op\cS^{\half}_1(\Sigma_p)$ in \eqref{nimmdas} that will be used throughout this section. For this define the operators $S_p$ over $\Sigma_p$ by 
\begin{equation} \label{def_SLBIO}
	S_p \varphi(x)=\int_{\Sigma_p} \cE_0(x-y)\varphi(y)\,\dd s(y), \quad x \in \Sigma_p, \quad \varphi \in C^\infty(\Sigma_p),\\
\end{equation}
where $\cE_0$ is given by~\eqref{fundamental_solution_Laplace}. It is known that $S_p$ gives rise to a compact self-adjoint operator in $L^2(\Sigma_p)$.
We denote the lower bound of the spectrum by $m_p:=\min \spec S_p\in (-\infty,0]$, the spectral projection corresponding to the spectrum of 
$S_p$ in $[m_p,0]$ by $\pi_p$,
and the corresponding spectral subspace by $\mathcal{H}_p^-=\ran\pi_p$.

\begin{lemma} \label{lemma_Lambda_p}
  For $S_p$ in \eqref{def_SLBIO} and $m_p,\pi_p,\mathcal{H}_p^-$ the following is true.
  \begin{itemize}
  \item[(a)] $S_p\in\Op\cS^{-1}_1(\Sigma_p)$
  \item[(b)] $\dim \mathcal{H}_p^- < \infty$ and $\pi_p \in \Op \cS^{-\infty}_1(\Sigma_p)$
  \item[(c)] The operator \begin{equation}
	\label{lambdap}
	\Lambda_p:=(S_p-(m_p-1) \pi_p)^{-\half}\in \Op\cS^{\half}_1(\Sigma_p)
  \end{equation}
  is elliptic, and, as an operator in $L^2(\Sigma_p)$, it is positive and self-adjoint with $\dom \Lambda_p=H^{\frac{1}{2}}(\Sigma_p)$; in particular,  $\Lambda_p$ is injective and symmetric in $L^2(\Sigma_p)$.  Its principal symbol is
\begin{equation}
	\label{sll}
	\sigma_{\Lambda_p}(x,\xi)=\sqrt{2|\xi|}.
\end{equation}
  \item[(d)] $\Lambda_p S_p \Lambda_p - I=(m_p-1)\Lambda_p^2\pi_p\in \Op \cS^{-\infty}_1(\Sigma_p)$
  \end{itemize}
\end{lemma}
\begin{proof}
First, recall that $S_p$
is the classical single layer boundary integral operator on $\Sigma_p$. In particular, one has (a) and
$S_p$ is elliptic, self-adjoint, and its principal symbol is $\sigma_{S_p}(x,\xi)=\frac{1}{2|\xi|}$, see, e.g., \cite[Section~10.2.2 and Section~11.4.1]{HW} and \cite[Theorem~6.22]{stein}. Moreover, again by \cite[Theorem~6.22]{stein} the space $\mathcal{H}_p^-=\ran\pi_p$ has finite dimension, which is the first claim in (b).
Let $(\varphi_j^p)_{j=1}^L \subset L^2(\Sigma_p)$ be an orthonormal basis of $\mathcal{H}_p^-$.
Then, by elliptic regularity $\varphi_j^p \in C^\infty(\Sigma_p)$, $j \in \{ 1, \dots, L \}$. Thus, with 
$$\Pi_p(x,y) = \sum_{j=1}^L \varphi_j^p(x) \overline{\varphi_j^p(y)} \in C^\infty(\Sigma_p \times \Sigma_p)$$ 
one can write
\begin{equation*}
  \pi_p \psi = \sum_{j=1}^L (\psi, \varphi_j^p)_{L^2(\Sigma_p)} \varphi_j^p = \int_{\Sigma_p} \Pi_p(\cdot, y) \psi(y) \dd s(y),
\end{equation*}
which implies the second claim in (b), that is, $\pi_p\in \Op\cS^{-\infty}_1(\Sigma_p)$. Next we verify (c).  As $m_p$ is the lower bound of the spectrum of $S_p$,  the operator $B_p:=S_p-(m_p-1) \pi_p \in \Op \cS^{-1}_1(\Sigma_p)$ is self-adjoint and positive: on $\mathcal{H}_p^-$ its spectrum is contained in $[1,1-m_p]$, and on the orthogonal complement of $\mathcal{H}_p^-$ it coincides with $S_p$, which is positive there. Moreover, $B_p$ is elliptic by (a) and (b), and hence its inverse $B_p^{-1}$ is a positive and self-adjoint elliptic operator in $\Op \cS^{1}_1(\Sigma_p)$ with $\dom B_p= H^1(\Sigma_p)$. Therefore, by Seeley's theorem on complex powers of elliptic operators \cite{Seeley},
\begin{equation*}
	\Lambda_p=(S_p-(m_p-1) \pi_p)^{-\half}=(B_p^{-1})^{\frac{1}{2}}\in \Op\cS^{\half}_1(\Sigma_p)
\end{equation*}
is elliptic, and it is positive and self-adjoint in $L^2(\Sigma_p)$ with $\dom \Lambda_p=H^{\frac{1}{2}}(\Sigma_p)$.
Since $\sigma_{S_p}(x,\xi)=\frac{1}{2|\xi|}$, we conclude that the principal symbol of $\Lambda_p$ is
\begin{equation*}
	\sigma_{\Lambda_p}(x,\xi)=\sqrt{2|\xi|}.
\end{equation*}
Finally, to see (d) observe that, since $S_p$ and $\pi_p$ commute,
\begin{equation}
\begin{split}
 \Lambda_p S_p \Lambda_p&=(S_p-(m_p-1) \pi_p)^{-\half}S_p(S_p-(m_p-1) \pi_p)^{-\half}\\
 &=(S_p-(m_p-1) \pi_p)^{-1}S_p\\
 &=I+ (S_p-(m_p-1) \pi_p)^{-1}(m_p-1) \pi_p\\
&=I+(m_p-1)\Lambda_p^2 \pi_p, 
\end{split}
 \end{equation}
 and since $\pi_p \in \Op \cS^{-\infty}_1(\Sigma_p)$ by (b) it follows that (d) holds.
\end{proof}

In our further analysis also the singular integral operators $W_p^\pm$ defined by
\begin{equation} \label{def_W_pm_op}
	W^\pm_p \varphi(x)=\pv \int_{\Sigma_p} w^\pm(x-y)\varphi(y)\,\dd s(y), \quad x \in \Sigma_p, \quad \varphi \in C^\infty(\Sigma_p),
	\end{equation}
where $\pv$ stands for the Cauchy principal value and
\begin{equation*} 
	w^\pm(x)=\begin{cases}
	\dfrac{\rmi}{2\pi}\dfrac{1}{x_1\pm \rmi x_2}, & n=2,\\[\bigskipamount]
	\dfrac{\rmi (\sigma\cdot x)}{4\pi |x|^3}, & n=3,
\end{cases}
\end{equation*}
will be important; cf. Lemma~\ref{lemma_expansion1}. In the next lemma, we study their pseudodifferential properties.	

For the three-dimensional case some further geometrical quantities are necessary. Consider for  $x\in\Sigma_p$ the shape operator
\begin{equation}
	\label{eq-vx}
	V_x:=\dd\nu^{\Sigma}|_x:T_x \Sigma_p\to T_x \Sigma_p.
\end{equation}	 Its eigenvalues $\kappa_1(x)$ and $\kappa_2(x)$ are the principal curvatures of $\Sigma_p$ in $x$. As $V_x$ is self-adjoint with respect to the scalar product inherited from $\RR^3$, one may find an orthonormal basis $(e_1,e_2)$ of $T_x \Sigma_p$ (we omit its dependence on $x$) such that $V_x e_j=\kappa_j(x)e_j$ for $j\in\{1,2\}$ and $\big(e_1,e_2,\nu^{\Sigma}(x)\big)$ is a positively oriented basis in $\RR^3$. With these notations we introduce the following special parametrizations of $\Sigma_p$.

\begin{definition}\label{def_special_parametrization}
Let $\Sigma_p$, $p \in \{ 1, \dots, P\}$, be one of the connected components of $\Sigma$.
\begin{itemize}
  \item[(a)] Let $n=2$ and $\ell_p$ be the length of $\Sigma_p$. We call an arc-length parametrization $\gamma: (-\frac{\ell_p}{2},\frac{\ell_p}{2})\to\Sigma_p$ (i.e. $|\gamma'(\cdot)|=1$) of $\Sigma_p$ \emph{properly oriented}, if $\gamma'(t)$ is obtained from the normal $\nu^{\Sigma}(\gamma(t))$ by the anti-clockwise rotation by $\frac{\pi}{2}$ for each $t$, i.e. $\nu^{\Sigma}(\gamma(t)) = (\gamma_2'(t), - \gamma_1'(t))$ for all $t \in (-\frac{\ell_p}{2},\frac{\ell_p}{2})$.
  \item[(b)] Let $n=3$, $x \in \Sigma_p$, and $e_1, e_2$ be as above. We say that a local chart $\varphi:\RR^2\supset U\to V\subset\Sigma_p$ corresponds to {\it normal coordinates} near $x$, if $0\in U$, $\varphi(0)=x$, and $\partial_j\varphi(0)=e_j$ for $j\in\{1,2\}$. 
\end{itemize}
\end{definition} 

In the next lemma the principal symbols of $W^\pm_p$ are specified.

\begin{lemma}\label{lem19}
One has $W_p^\pm \in \Op \cS_{\nhalf}^0 (\Sigma_p)$ and the following is true.
\begin{itemize}
  \item[(a)] For $n=2$ the principal symbols of $W^\pm_p$ at $x\in\Sigma_p$ with respect to any properly oriented arc-length
parametrization (see Definition~\ref{def_special_parametrization}~\textup{(a)}) are
\begin{equation}
	\label{ww-2d}
\sigma_{W^+_p}(x,\xi)=-\dfrac{\rmi}{2} \overline{\mu}(x) \sgn\xi,
\quad
\sigma_{W^-_p}(x,\xi)=\dfrac{\rmi}{2} \mu(x) \sgn\xi,
\end{equation}
where $\mu:\Sigma_p\to\CC$ is the standard complexification of $\nu^\Sigma$  given by
\begin{equation} \label{complexification_normal_2d}
\mu(x):=\nu^\Sigma_1(x)+\rmi \nu^\Sigma_2(x),\quad x\in\Sigma_p.
\end{equation}

  \item[(b)] For $n=3$ the principal symbols of $W^\pm_p$ at $x\in\Sigma_p$ with respect to normal coordinates near $x$ (see Definition~\ref{def_special_parametrization}~\textup{(b)}) are 
\begin{equation}
	\label{ww-3d}	
		\sigma_{W^\pm_p}(x,\xi)=\dfrac{1}{2}\sigma\cdot \left(\dfrac{\xi_1 e_1+\xi_2 e_2}{|\xi|} \right).
\end{equation}
\end{itemize}
\end{lemma}
\begin{proof}
We will compute the principal symbols of $W_p^\pm$; the fact that $W_p^\pm \in \Op \cS_{\nhalf}^0 (\Sigma_p)$ will follow as a byproduct of these considerations.

(a) In the following it is convenient to complexify $\RR^2$ by $(x_1, x_2)\simeq x_1+\rmi x_2$. Let	$x\in\Sigma_p$ and $\gamma:I\to \Sigma_p$
with $I=(-\frac{\ell_p}{2},\frac{\ell_p}{2})$
be a properly oriented complexified arc-length parametrization of $\Sigma_p$ with $\gamma(0)=x$. Due to the arc-length and orientation conditions one has for all $t\in I$ that
\[
\dfrac{1}{\gamma'(t)}=\overline{\gamma'(t)} \text{ and }
\gamma'(t)=\rmi \mu\big(\gamma(t)\big).
\]
 Let $J:= (-\frac{\ell_p}{4},\frac{\ell_p}{4})$ and $V:=\gamma(J)$,  which is an open neighborhood of $x$ in $\Sigma_p$, pick $\psi\in C^\infty_c(V)$, and let $f\in C^\infty(\Sigma_p)$. Then for any $s\in I$
 one has
\begin{align*}
	\big[W^+_p (\psi f)\big]\big(\gamma(s)\big)&=\dfrac{\rmi}{2\pi}\pv \int_I \dfrac{1}{\gamma(s)-\gamma(t)}(\psi f)\big(\gamma(t)\big)\,\dd t\\
	&=\dfrac{\rmi}{2\pi}\pv \int_I a(s,t)\dfrac{(\psi f)\big(\gamma(t)\big)}{s-t}\,\dd t
\end{align*}
with the function  $a: I\times I\to \CC $  defined by
\begin{equation*}
\begin{split}
a(s,t):&= \begin{cases} \frac{s-t}{\gamma(s)-\gamma(t)}, & s \neq t, \\ \frac{1}{\gamma'(s)},& s=t,\end{cases} \\
&= \begin{cases} \frac{s-t}{\gamma(s)-\gamma(t)}, & s \neq t, \\  -\rmi\overline{\mu}\big(\gamma (s)\big),& s=t.\end{cases}
\end{split}
\end{equation*}
Hence, we can further rewrite
\begin{equation} \label{decomposition_W_p}
\big[W^+_p (\psi f)\big]\big(\gamma(s)\big)=\dfrac{1}{2\pi}\overline{\mu}\big(\gamma (s)\big)\pv\int_I \dfrac{(\psi f)\big(\gamma(t)\big)}{s-t}\,\dd t
+ \int_I b(s,t)(\psi f)\big(\gamma(t)\big)\,\dd t
\end{equation}
with 
\begin{equation*}
  b(s,t):=\frac{\rmi}{2\pi} \cdot \begin{cases} \frac{a(s,t)-a(s,s)}{s-t}, &s \neq t, \\ -\partial_2 a(s,s), & s=t. \end{cases}
\end{equation*}
 As $\gamma$ is injective on $I$, $\gamma'\neq0$, and $\gamma(s)$ stays away from $\gamma(\overline{J})$ as $s\to \pm\frac{\ell_p}{2}$, one has $a\in C^\infty(I\times J)$. Hence,  $b\in C^\infty(I\times J)$ and thus, as $\psi\circ \gamma$ has compact support in $J$,  the second summand on the right hand side in~\eqref{decomposition_W_p} defines an operator in $\Op \cS^{-\infty}_1(\Sigma_p)$. Furthermore, the first summand in~\eqref{decomposition_W_p} can be written as the convolution
\[
\dfrac{1}{2\pi}\overline{\mu}\big(\gamma (s)\big)\pv\int_I \dfrac{(\psi f)\big(\gamma(t)\big)}{s-t}\,\dd t = \dfrac{1}{2\pi} \overline{\mu}(\gamma(s))\left(\pv\frac{1}{s}\right) \star ((\psi f)\circ\gamma).
\]
Recall that the Fourier transform of the distribution $\pv \frac{1}{s}$ on $\RR$ is given by
\[
\widehat{\pv \tfrac{1}{s}} =-\rmi \pi\sgn\xi.
\]
Therefore, we conclude that the principal symbol of $W^+_p$ at $x$ (corresponding to $\gamma(0)$) is
\[
\sigma_{W^+_p}(x,\xi)=\dfrac{1}{2\pi} \overline{\mu}(x) \widehat{\pv \dfrac{1}{s}}
\equiv -\dfrac{\rmi}{2} \overline{\mu}(x) \sgn\xi.
\]
From this and $W^-_p f=-\overline{W^+_p \Bar f}$ we obtain that the principal symbol of $W^-_p$ at $x$ is $\sigma_{W^-_p}(x,\xi)=-\overline{\sigma_{W^+_p}(x,-\xi)}$. Hence, all claims in (a) are shown.

(b) is shown in \cite[Theorem~8]{BP24}. It should be noted that in \cite[Remark~6]{BP24} there was a sign mistake: the curvatures were defined with the opposite sign, while all the subsequent computations actually used the definition corresponding to \eqref{eq-vx}.
\end{proof}

As a direct corollary of the expressions \eqref{ww-2d} and \eqref{ww-3d} we obtain the following facts.
\begin{corollary}\label{cor-ww14}
	$W^+_p W^-_p-\frac{1}{4}I, W^-_p W^+_p-\frac{1}{4}I \in \Op\cS^{-1}_\nhalf(\Sigma_p)$.
\end{corollary}

In order to compute the essential spectrum of $A^{\eta, \tau}$ for $\eta_p^2-\tau_p^2=4$ we need the following finer microlocal properties of its blocks. Recall that for $n=3$ the operators $W_p^\pm$ are equal, which allows us to set
\begin{equation*}
  W_p := W^+_p=W_p^-,
\end{equation*}
and that normal coordinates near $x$ were introduced in Definition~\ref{def_special_parametrization}~(b).

\begin{lemma}\label{lem26}
For the singular integral operators $W_p^\pm$ the following is true: 
\begin{itemize}
\item[(a)] Let $n=2$ and $\mu$ be the complexification of $\nu^\Sigma$ as in~\eqref{complexification_normal_2d}. Then, it holds that
\begin{gather}
	\label{www}
(W_p^+\mu)^2=-\dfrac{1}{4} I=(\Bar\mu W_p^-)^2,\\
	\label{www2}
W^+_p\mu+\Bar\mu W^-_p,\ W^+_pW^-_p-\dfrac{1}{4}I,\ W^-_pW^+_p-\dfrac{1}{4}I\in\Op \cS^{-\infty}_1(\Sigma_p).
\end{gather}

\item[(b)]  Let $n=3$. Then, it holds that
		\begin{equation}
			\label{w3d}
			\big((\sigma\cdot\nu^\Sigma)W_p\big)^2=\big(W_p (\sigma\cdot\nu^\Sigma)\big)^2=-\dfrac{1}{4}I_2.
		\end{equation}
		Furthermore, the operator
		\begin{equation}
			\label{eq-zp}
		Z_p:= W_p(\sigma\cdot\nu^\Sigma)+(\sigma\cdot\nu^\Sigma)W_p
		\end{equation}
		belongs to $\Op \cS^{-1}_{2}(\Sigma_p)$, and for any $x\in\Sigma_p$ its principal symbol
		with respect to normal coordinates near $x$ is given by
		\begin{equation}
			\label{eq-zsymb}
			\sigma_{Z_p}(x,\xi)= \frac{1}{2}\big(\kappa_1(x)-\kappa_2(x)\big)(\sigma\cdot \nu^\Sigma(x))\,	\frac{\xi_1 \xi_2}{|\xi|^3},
		\end{equation}
		where $\kappa_1(x)$ and $\kappa_2(x)$ are 
the principal curvatures of $\Sigma_p$ in $x$.
\end{itemize}
\end{lemma}

\begin{proof}
	(a) Let $\gamma:I\to \Sigma_p$ with $I=(-\frac{\ell_p}{2},\frac{\ell_p}{2})$
be a properly oriented complexified arc-length parametrization of $\Sigma_p$. Note that the operator
	$T:=\rmi W^+_p\mu$ acts on $\varphi \in C^\infty(\Sigma_p)$ as 
	\begin{equation*}
	  \begin{split}
  	  T \varphi(\gamma(s)) &= \frac{\rmi}{2\pi} \pv \int_{I} \frac{{\gamma}'_1(t) + \rmi \gamma'_2(t)}{\gamma_1(s) + \rmi \gamma_2(s) - (\gamma_1(t) + \rmi \gamma_2(t))} \varphi(\gamma(t)) \dd t = \frac{\rmi}{2\pi} \pv \int_{\widetilde{\Sigma}_p} \frac{\widetilde{\varphi}(\zeta)}{z - \zeta}  d \zeta,
  	  \end{split}
	\end{equation*}
	where $\widetilde{\Sigma}_p$, $\widetilde{\varphi}$, and $z$ are the complexified versions of $\Sigma_p$, $\varphi$, and $\gamma(s)$, respectively, i.e.
	\begin{equation*}
	  \widetilde{\Sigma}_p = \{ x_1 + \rmi x_2: (x_1, x_2) \in \Sigma_p \}, \quad \widetilde{\varphi}(x_1 + i x_2) = \varphi(x_1,x_2), \quad z = \gamma_1(s) + \rmi \gamma_2(s).
	\end{equation*}
	Hence, $T$ is the well-known Hilbert transform on $\Sigma_p$ (i.e. the direct value of the respective Cauchy integral), see, e.g., \cite[\S12]{musk}, which is known to satisfy $4T^2=I$, see e.g. \cite[\S27]{musk}. This gives the first identity in \eqref{www}. We further note that $\Bar\mu W^-_p$ is the adjoint of $W^+_p\mu$ in $L^2(\Sigma_p)$, which gives the second identity in \eqref{www}.
	The operator $K:=W^+_p\mu+\Bar\mu W^-_p$ is known  (modulo constant factors) as the Kerzman-Stein operator: it has a smooth integral kernel, see, e.g., \cite[Section~5]{bell}, which gives the inclusion $K\in \Op \cS^{-\infty}_1(\Sigma_p)$.
	Finally we have
	\begin{align*}
		 W^+_pW^-_p&=(W^+_p\mu) (\Bar \mu W^-_p)=W^+_p\mu(-W^+_p\mu +K)\stackrel{\eqref{www}}{=}\dfrac{1}{4}I +W^+_p \mu K,\\
		 W^-_pW^+_p&=\mu  (\Bar \mu W^-_p) (W^+_p\mu)\Bar\mu=\mu (-W^+_p\mu +K) (W^+_p\mu)\Bar\mu \stackrel{\eqref{www}}{=}
		 \dfrac{1}{4}I+\mu K W_p^+,
	\end{align*}
	while $W^+_p \mu K,\mu K W^+_p \in \Op \cS^{-\infty}_1(\Sigma_p)$ due to $W^+_p, \mu \in \Op \cS^{0}_1(\Sigma_p)$ and $K\in \Op \cS^{-\infty}_1(\Sigma_p)$. This shows~\eqref{www2}.
	
(b) For the relation \eqref{w3d} see \cite[Section~2, equation~(11)]{BP24}, and the expression~\eqref{eq-zsymb} is derived in \cite[Theorem~8]{BP24}; our situation corresponds to $\lambda=m$ in \cite{BP24}.
\end{proof}

In the following corollary we state the pseudodifferential properties of the operator $C_z$ defined by~\eqref{def_C_z}. 
Let $S_\Sigma$ and $W_\Sigma^\pm$ be the integral operators that are defined as in~\eqref{def_SLBIO} and~\eqref{def_W_pm_op} with the integral 
over $\Sigma$ instead of $\Sigma_p$, that is, 
\begin{equation} \label{def_SLBIO2}
	S_\Sigma\varphi(x)=\int_{\Sigma} \cE_0(x-y)\varphi(y)\,\dd s(y), \quad x \in \Sigma, \quad \varphi \in C^\infty(\Sigma),
\end{equation}
and 
\begin{equation} \label{def_W_pm_op2}
	W^\pm_\Sigma \varphi(x)=\pv \int_{\Sigma} w^\pm(x-y)\varphi(y)\,\dd s(y), \quad x \in \Sigma, \quad \varphi \in C^\infty(\Sigma).
	\end{equation}
It is important to note that one has 
	\begin{equation}
	\label{oplus}
	\begin{aligned}
		S_\Sigma&=\bigoplus_{p=1}^PS_p\mod \Op \cS^{-\infty}_{1}(\Sigma)\quad\text{and}\quad
		W^\pm_\Sigma&=\bigoplus_{p=1}^P W^\pm_p\mod \Op \cS^{-\infty}_{\frac{N}{2}}(\Sigma),
	\end{aligned}
\end{equation}	
since for $p\not=q$ operators of the form (with $\cE_0, w^\pm \in C^\infty(\RR^n \setminus \{ 0 \})$)
	\begin{align*}
	S_{p,q} \psi(x)&=\int_{\Sigma_p}\cE_0(x-y)\psi(y)\,\dd s(y),\quad x \in \Sigma_q, \quad \psi \in C^\infty(\Sigma_p),\\
	W^\pm_{p,q} \psi(x)&=\pv \int_{\Sigma_p} w^\pm(x-y)\psi(y)\,\dd s(y), \quad x \in \Sigma_q, \quad \psi \in C^\infty(\Sigma_p),
	\end{align*}
are smoothing. Hence, it follows from Lemma~\ref{lemma_Lambda_p} and Lemma~\ref{lem19} that $S_\Sigma\in\Op\cS^{-1}_1(\Sigma)$ and $W_\Sigma^\pm \in \Op \cS_{\nhalf}^0 (\Sigma)$. With these notations we have the following result:

\begin{corollary}\label{corollary_C_z_pseudo} Let $z\in\res A$. Then $C_{z} \in \Op \cS_{N}^0 (\Sigma)$ and
	\begin{equation}
		\label{cblock}
	C_{z}=\begin{pmatrix}
		(z+m)S_\Sigma\otimes I_\nhalf & W^+_\Sigma\\
		W^-_\Sigma & (z-m)S_\Sigma\otimes I_{\nhalf}
	\end{pmatrix}
	\mod \Op \cS^{-2}_N(\Sigma).
	\end{equation}
\end{corollary}

 \begin{proof}
 	Using the notation in Theorem~\ref{theo1} and Lemma~\ref{lemma_expansion1} one has for $x,y\in\Sigma$, $x\not= y$, 
 	\begin{equation*}
 	\begin{split}
 	 G_z(x,y)&=E_z(x-y)+T_z(x,y)\\
 	 &=\mathcal H_z(x-y)+\chi(x-y)\bigl(E_z(x-y)-\mathcal H_z(x-y)\bigr) +T_z(x,y)\\
 	 &=\mathcal H_z(x-y)+\mathcal K_z(x-y) +T_z(x,y),
 	\end{split}
 	\end{equation*}
 	where $\chi\in C_c^\infty(\RR^n)$ is chosen such that $\chi(x-y)=1$ for $x,y\in\Sigma$.
 	Let $K_{\Sigma,z}$ and $\widetilde{T}_{\Sigma,z}$ be defined as in~\eqref{def_int_ops}. Then it follows that
 	\begin{equation*}
 	  C_{z} =\begin{pmatrix}
		(z+m)S_\Sigma\otimes I_\nhalf & W^+_\Sigma\\
		W^-_\Sigma & (z-m)S_\Sigma\otimes I_{\nhalf}
	\end{pmatrix}  + K_{\Sigma, z} + \widetilde{T}_{\Sigma,z}.
 	\end{equation*}
 	As it was shown in the proof of Lemma~\ref{lem14} one has $K_{\Sigma,z} \in \Op \cS^{-2}_N(\Sigma)$ and 
 	$\widetilde{T}_{\Sigma,z} \in \Op \cS^{-\infty}_N(\Sigma)$. Since $S_\Sigma\in\Op\cS^{-1}_1(\Sigma)$ and $W_\Sigma^\pm \in \Op \cS_{\nhalf}^0 (\Sigma)$, this shows the claim.
 \end{proof}

Recall that $\Theta^{\eta, \tau}$ is given by~\eqref{def_Theta} and that $M_{z}$ denotes the Weyl function associated with the boundary triple $( L^2(\Sigma,\CC^N), \Pi_{0}, \Pi_{1})$, see~\eqref{Weyl_function}.
Next, we provide a decomposition of 
\begin{equation}\label{besserhier}
\Theta^{\eta,\tau} - M_{z}=-\Lambda\left[\frac{1}{\eta^2-\tau^2}(\eta I_N -\tau\beta)+C_{z}\right]\Lambda
\end{equation}
into parts, which will be easier to analyze in the following. Recall that $\Lambda_p\in \Op\cS^{\half}_1(\Sigma_p)$ and $ W^\pm_p \in\Op\cS^{0}_{\nhalf}(\Sigma_p)$ are defined by~\eqref{lambdap} and~\eqref{def_W_pm_op}, respectively. For $p\in\{1,\dots,P\}$ such that $\vert\eta_p\vert\not=\vert\tau_p\vert$ we consider the operators 
\begin{equation}
	\label{lql2}
\begin{aligned}
	L_{p}&:=-\begin{pmatrix}
		\left(\dfrac{1}{\eta_p+\tau_p}\Lambda_p^2+m\right)\otimes I_{\nhalf} & \Lambda_pW^+_p\Lambda_p\\
		\Lambda_pW^-_p\Lambda_p & \left(\dfrac{1}{\eta_p-\tau_p}\Lambda_p^2-m\right)\otimes I_{\nhalf}
	\end{pmatrix} \in \Op\cS^{1}_{N}(\Sigma_p)
\end{aligned}
\end{equation}
and have the following result.

\begin{proposition} \label{proposition_decomposition_Theta}
  Let $\Lambda = \bigoplus_{p=1}^P \Lambda_p$ with the special choice of $\Lambda_p$ from~\eqref{lambdap}, assume that 
  $\vert\eta_p\vert\not=\vert\tau_p\vert$ for all $p\in\{1,\dots,P\}$, and
  let $L_{p}$ be as above. Then, for any $z \in \res A$ one has
  \begin{equation} \label{decomposition_Theta}
    \Theta^{\eta,\tau} - M_{z} = \bigoplus_{p=1}^P (L_{p} - z) \mod \Op \cS^{-1}_N(\Sigma).
  \end{equation}
  In particular, each $L_p$ viewed as an operator in $L^2(\Sigma_p,\CC^N)$ on its maximal domain is self-adjoint, $H^1(\Sigma_p,\CC^N)$ is a core for $L_p$, and $\bigoplus_{p=1}^P L_{p}$ defined on $\dom\Theta^{\eta,\tau}$ is self-adjoint in $L^2(\Sigma,\CC^N)$.
\end{proposition}
\begin{proof}
Let $z \in \res A$ be fixed. 
Using the block representation from \eqref{cblock} and the fact that $\Lambda \in \Op\cS^{\half}_N(\Sigma)$ one gets that
\begin{equation} \label{decomposition1}
  \begin{split}
	\Theta^{\eta,\tau} - M_{z} &= -\Lambda\bigg[\dfrac{1}{\eta^2-\tau^2}(\eta I_N -\tau\beta)+C_{z}\bigg]\Lambda \\
	&=-\Lambda
	\begin{pmatrix}
		\left(\dfrac{1}{\eta+\tau}+(z+m)S_\Sigma\right)\otimes I_{\frac{N}{2}} & W^+_\Sigma\\
		W^-_\Sigma & \left(\dfrac{1}{\eta-\tau}+(z-m)S_\Sigma\right)\otimes I_{\frac{N}{2}}
	\end{pmatrix}
	\Lambda\\
	&\qquad \qquad \qquad \qquad \qquad \qquad \qquad \qquad \qquad \qquad \qquad \mod \Op\cS^{-1}_N(\Sigma).
  \end{split}
\end{equation}
According to \eqref{oplus} we have
\begin{equation}
	\label{oplus2}
	\begin{aligned}
		S_\Sigma&=\bigoplus_{p=1}^PS_p\mod \Op \cS^{-\infty}_{1}(\Sigma)\quad\text{and}\quad
		W^\pm_\Sigma&=\bigoplus_{p=1}^P W^\pm_p\mod \Op \cS^{-\infty}_{\frac{N}{2}}(\Sigma),
	\end{aligned}
\end{equation}	
and in view of \eqref{decomposition1} this implies the decomposition
\begin{equation}
	\label{oplus3}
	\begin{aligned}
		\Theta^{\eta,\tau} - M_{z} &=\bigoplus_{p=1}^P Q_{p,z} \mod \Op\cS^{-1}_N(\Sigma)
	\end{aligned}
\end{equation}
with
\begin{align*}
	Q_{p,z}
	&:=-\Lambda_p\begin{pmatrix}
		\left(\dfrac{1}{\eta_p+\tau_p}+(z+m)S_p\right)\otimes I_{\frac{N}{2}} & W^+_p\\
		W^-_p & \left(\dfrac{1}{\eta_p-\tau_p} +(z-m)S_p\right)\otimes I_{\frac{N}{2}}
	\end{pmatrix}\Lambda_p\\
	&=-\begin{pmatrix}
		\left(\frac{1}{\eta_p+\tau_p}\Lambda_p^2+(z+m)\Lambda_pS_p\Lambda_p\right)\otimes I_{\frac{N}{2}} & \Lambda_pW^+_p\Lambda_p\\
		\Lambda_pW^-_p\Lambda_p & \left(\frac{1}{\eta_p-\tau_p}\Lambda_p^2 +(z-m)\Lambda_pS_p\Lambda_p\right)\otimes I_{\frac{N}{2}}
	\end{pmatrix}.
\end{align*}
Since $\Lambda_p S_p \Lambda_p =I+ (m_p-1)\Lambda_p^2\pi_p$ and $(m_p-1)\Lambda_p^2\pi_p \in \Op \cS^{-\infty}_1(\Sigma_p)$
by Lemma~\ref{lemma_Lambda_p}~(d), one finds that
\begin{equation*}
\begin{aligned}
	Q_{p,z}&=L_{p} - z\mod \Op \cS^{-\infty}_N(\Sigma_p).
\end{aligned}
\end{equation*}
By combining this with \eqref{oplus3} one arrives at~\eqref{decomposition_Theta}. 
Finally,  $L_p$ is a symmetric first order pseudodifferential operator on $\Sigma_p$, and hence it is self-adjoint on the maximal domain due to Friedrichs lemma, see the dicussion in Appendix~\ref{appendix_pseudos}. Since, by the same reason, the adjoint of $L_p$ defined on $C^\infty(\Sigma_p,\CC^N)$ is the maximal realization of $L_p$, the self-adjointness of the latter implies that  $C^\infty(\Sigma_p,\CC^N)$, and hence also $H^1(\Sigma_p,\CC^N)$, is a core for $L_p$. 
Therefore, $\bigoplus_{p=1}^P L_p$ is self-adjoint in $L^2(\Sigma, \CC^N)$ on its maximal domain, which coincides, as $M_z$ and all operators in $\Op \cS^{-1}_N(\Sigma)$ are bounded in $L^2(\Sigma, \CC^N)$, with $\dom \Theta^{\eta, \tau}$ by~\eqref{decomposition_Theta}. This finishes the proof.
\end{proof}

\subsection{Self-adjointness and essential spectrum}\label{Section 4.2}

First, making use of Lemma~\ref{lem16}, the abstract theory of boundary triples from Appendix~\ref{sec-bt}, and the properties of $\Theta^{\eta,\tau}$ shown in Section~\ref{section_separation}, we prove the self-adjointness of $A^{\eta, \tau}$ defined in~Definition~\ref{adefjussi} and relate its essential spectrum to the essential spectra of the operators $L_p$ defined in~\eqref{lql2}.
Throughout this subsection we assume that we always have the special choice of $\Lambda_p$ from Lemma~\ref{lemma_Lambda_p}.

\begin{theorem}\label{thm17}
Assume that $|\eta_p|\ne |\tau_p|$ for all $p\in\{1,\dots,P\}$,
let $\Gamma_z$ and $M_z$, $z \in \res A$, be given by~\eqref{eq-gbz} and~\eqref{Weyl_function}, respectively, and let $\Theta^{\eta,\tau}$ be as in Lemma~\ref{lem16}. 
Then the operator $A^{\eta,\tau}$ in Definition~\ref{adefjussi} is 
self-adjoint in $L^2(\Omega,\CC^N)$ and the following holds:
\begin{itemize}
	\item[(a)] For any $z\in \res A\cap\res A^{\eta,\tau}$ the operator
	\[	
\Theta^{\eta,\tau} - M_{z}=-\Lambda\bigg[\dfrac{1}{\eta^2-\tau^2}(\eta I_N -\tau\beta)+C_z\bigg]\Lambda
	\]
	defined on the maximal domain $\dom\Theta^{\eta,\tau}$ is bijective in $L^2(\Sigma,\CC^N)$, and	the resolvent formula
	\[
	(A^{\eta,\tau}-z)^{-1}=(A-z)^{-1}+\Gamma_z(\Theta^{\eta,\tau}-M_z)^{-1}\Gamma_{\Bar z}^*
	\]
	holds.
	\item[(b)] For any $z\in\res A$ one has 
	\begin{equation*}
	 z\in \spec A^{\eta,\tau}  \quad\text{if and only if}\quad 0\in\spec  (\Theta^{\eta,\tau} - M_{z})
	\end{equation*}
	and 
	\begin{equation*}
	 z\in \spece A^{\eta,\tau}  \quad\text{if and only if}\quad 0\in\spece  (\Theta^{\eta,\tau} - M_{z}),
	\end{equation*}
	and the equality $\dim\ker (A^{\eta,\tau}-z)=\dim\ker (\Theta^{\eta,\tau} - M_{z})$ holds. In particular, 
    \begin{equation}
	\label{lol5}
	\spece A^{\eta,\tau} \cap \res A = \bigg[\bigcup_{p=1}^P\spece L_{p}\bigg] \cap \res A,
\end{equation}
 where the operator $L_p$ is defined as in \eqref{lql2}.
\end{itemize}
\end{theorem}

\begin{proof}
Recall that $A^{\eta,\tau}$ is the restriction of $S^*$ given by~\eqref{star1} on the functions $f$ satisfying 
\[
\Pi_0 f\in \dom \Theta^{\eta,\tau}\quad\text{and}\quad
\Pi_1 f= \Theta^{\eta,\tau}\Pi_0 f,
\]
where  $\big( L^2(\Sigma,\CC^N), \Pi_0, \Pi_1\big)$ is
the boundary triple
in Lemma~\ref{lem15}; cf. Lemma~\ref{lem16}.
Next, note that $\Lambda  \in \Op \cS^{\frac{1}{2}}_N(\Sigma)$ by construction, and $C_\zeta  \in \Op \cS^{0}_N(\Sigma)$ by Corollary \ref{corollary_C_z_pseudo}. Moreover, since $\zeta\in\RR$, the operator $C_\zeta$ is symmetric in $L^2(\Sigma,\CC^N)$ by \eqref{mk2}. Therefore, $\Theta^{\eta,\tau}$ is a symmetric first order pseudodifferential operator on $\Sigma$, and hence it is self-adjoint on the maximal domain due to Friedrichs' lemma, see the dicussion in Appendix~\ref{appendix_pseudos}. Therefore, $A^{\eta,\tau}$ is self-adjoint by Theorem~\ref{thm12}. The claims in (a) and the relations of the spectra of $A^{\eta, \tau}$  and $\Theta^{\eta,\tau} - M_{z}$ follow by using Theorem~\ref{thm12} with the specific expressions obtained in Lemma~\ref{lem15} and Lemma~\ref{lem16}.
Finally, 
as any operator in $\Op \cS^{-1}_N(\Sigma)$ gives rise to a bounded operator from $L^2(\Sigma, \CC^N)$ to $H^1(\Sigma, \CC^N)$ and hence, taking Rellich's embedding theorem into account, to a compact operator in $L^2(\Sigma, \CC^N)$, one concludes from \eqref{decomposition_Theta} that
\begin{equation*}
\spece  (\Theta^{\eta, \tau} - M_{z}) =\bigcup_{p=1}^P \spece (L_{p} - z).
\end{equation*}
This yields~\eqref{lol5} and finishes the proof.
\end{proof}

Next, we consider the \textit{non-critical case} 
$\eta_p^2-\tau_p^2\not=4$ and $\vert\eta_p\vert\not=\vert\tau_p\vert$. Here, it turns out that there is no contribution to the essential spectrum.

\begin{lemma}\label{lem21}
Let $p\in\{1,\dots,P\}$ and assume that $\eta_p^2-\tau_p^2\not=4$ and $\vert\eta_p\vert\not=\vert\tau_p\vert$.
Then, $L_{p}$ is a first order elliptic operator. In particular, when $L_p$ is viewed as a self-adjoint operator in $L^2(\Sigma_p,\CC^N)$, then
$\dom L_p=H^1(\Sigma_p,\CC^N)$ and 
$\spece L_{p}=\emptyset$.
\end{lemma}	

\begin{proof}
	By definition $L_{p}$ is of first order, so we simply need to show that its principal symbol is invertible for $|\xi|$ sufficiently large.
	For any choice of local coordinates
	near $x\in\Sigma_p$ the principal symbol $\sigma_{L_{p}}(x,\xi)$ of $L_{p}$
	is 
	\[
	\sigma_{L_{p}}(x,\xi)=
	-\begin{pmatrix}
		\dfrac{1}{\eta_p+\tau_p}I_{\frac{N}{2}} & \sigma_{W^+_p}(x,\xi)\\
		\sigma_{W^-_p}(x,\xi) & \dfrac{1}{\eta_p-\tau_p} I_{\frac{N}{2}}
	\end{pmatrix}\sigma_{\Lambda_p}(x,\xi)^2,
	\]
	and we recall that $\sigma_{\Lambda_p}(x,\xi)^2=2|\xi|$ by construction, see Lemma~\ref{lemma_Lambda_p}~(c). 
	
	Consider the case $n=2$, then $N=2$. Assume that the local coordinates are given by a properly oriented arc-length parametrization (see Definition~\ref{def_special_parametrization}~(a)). Then by Lemma~\ref{lem19} we obtain
	\begin{equation*}
	\sigma_{L_{p}}(x,\xi)=
	-2|\xi|\begin{pmatrix}
		\dfrac{1}{\eta_p+\tau_p} & -\dfrac{\rmi}{2} \overline{\mu}(x) \sgn\xi\\
		\dfrac{\rmi}{2} \mu(x) \sgn\xi & \dfrac{1}{\eta_p-\tau_p}
	\end{pmatrix}
	\end{equation*}
	and thus,
	\begin{equation*}
	\det \sigma_{L_{p}}(x,\xi)=4\xi^2\left(\dfrac{1}{\eta_p^2-\tau_p^2}-\dfrac{1}{4}\right).
	\end{equation*}
	Since $\eta_p^2-\tau_p^2\ne 4$ by assumption, the determinant is non-zero for all $\xi\ne 0$, implying that $\sigma_{L_{p}}(x,\xi)$ is invertible for $\xi \neq 0$ and therefore, $L_{p}$ is elliptic in the case $n=2$.
	
	Now let $n=3$, then $N=4$. Using normal coordinates near $x$ (see Definition~\ref{def_special_parametrization}~(b))
	we arrive at
	\[
	\sigma_{L_{p}}(x,\xi)=
	-2|\xi|\begin{pmatrix}
		\dfrac{1}{\eta_p+\tau_p}I_2 & \dfrac{1}{2}\sigma\cdot \left(\dfrac{\xi_1 e_1+\xi_2 e_2}{|\xi|} \right)\\[\bigskipamount]
		\dfrac{1}{2}\sigma\cdot \left(\dfrac{\xi_1 e_1+\xi_2 e_2}{|\xi|} \right) & \dfrac{1}{\eta_p-\tau_p} I_2
	\end{pmatrix}.
	\]
	Therefore, for $\xi \neq 0$ one has
	\[
	\sigma_{L_{p}}(x,\xi)\cdot \dfrac{1}{2|\xi|}\begin{pmatrix}
		\dfrac{1}{\eta_p-\tau_p}I_2 & -\dfrac{1}{2}\sigma\cdot \left(\dfrac{\xi_1 e_1+\xi_2 e_2}{|\xi|} \right)\\[\bigskipamount]
		-\dfrac{1}{2}\sigma\cdot \left(\dfrac{\xi_1 e_1+\xi_2 e_2}{|\xi|} \right) & \dfrac{1}{\eta_p+\tau_p} I_2
	\end{pmatrix}
	=\left(\dfrac{1}{4}-\dfrac{1}{\eta_p^2-\tau_p^2}\right)I_4,
	\]
	which shows that $\sigma_{L_{p}}(x,\xi)$ is invertible for all $\xi\ne 0$, if $\eta_p^2-\tau_p^2\ne 4$. Hence, $L_{p}$ is also elliptic in the case $n=3$.

	Finally, as $L_{p}$ is a first-order elliptic operator we have $\dom L_p=H^1(\Sigma_p,\CC^N)$, when $L_p$ is viewed as a self-adjoint operator in $L^2(\Sigma_p,\CC^N)$. Furthermore, it follows that the resolvent of $L_p$ is compact
	and thus $\spece L_{p} = \emptyset$.
\end{proof}
	
As an immediate consequence of the ellipticity of $L_{p}$ in the non-critical case one gets an improved Sobolev regularity near $\Sigma_p$ for functions belonging to $\dom A^{\eta, \tau}$. Recall that the space $H^{1,\Sigma_p}_\alpha(\Omega\setminus\Sigma,\CC^N)$ is defined by~\eqref{def_H_s_Sigma_p}.

\begin{corollary}\label{cor21}
Let $p\in\{1,\dots,P\}$ and assume that $\eta_p^2-\tau_p^2\not=4$ and $\vert\eta_p\vert\not=\vert\tau_p\vert$.
Then 
\[
\dom A^{\eta,\tau} \subset H^{1,\Sigma_p}_\alpha(\Omega\setminus\Sigma,\CC^N).
\]
\end{corollary}	

\begin{proof}
Let  $f\in\dom A^{\eta,\tau}$. Then $\Pi_{0}f\in \dom \Theta^{\eta,\tau}$.
According to Proposition~\ref{proposition_decomposition_Theta} one has
\begin{equation*}
  \dom \Theta^{\eta,\tau} = \bigoplus_{p=1}^P \dom L_{p},
\end{equation*}
which implies
\[
(\Pi_{0} f)|_{\Sigma_p} = \rmi \Lambda_p^{-1}(\alpha\cdot\nu^\Sigma)(\gamma^\Sigma_+ f-\gamma^\Sigma_- f)|_{\Sigma_p}\in \dom L_{p}= H^1(\Sigma_p,\CC^N),
\]
where the last equality holds by Lemma~\ref{lem21}.
This yields $(\gamma^\Sigma_+ f-\gamma^\Sigma_- f)|_{\Sigma_p}\in H^\half(\Sigma_p,\CC^N)$.
Next, note that the assumption $\vert\eta_p\vert\not=\vert\tau_p\vert$ implies that the constant matrix $\eta_p I_N + \tau_p \beta$ is invertible. Hence, due to the transmission condition for $f \in \dom A^{\eta, \tau}$ one gets  
\begin{equation*}
  (\gamma^\Sigma_+ f +\gamma^\Sigma_- f)|_{\Sigma_p} = -2 (\eta_p I_N + \tau_p \beta)^{-1} \rmi (\alpha \cdot \nu^\Sigma) (\gamma^\Sigma_+ f-\gamma^\Sigma_- f)|_{\Sigma_p} \in H^\half(\Sigma_p,\CC^N).
\end{equation*}
Overall, we obtain $(\gamma^\Sigma_\pm f)|_{\Sigma_p}\in H^\half(\Sigma_p,\CC^N)$, and then $f\in  H^{1,\Sigma_p}_\alpha(\Omega\setminus\Sigma,\CC^N)$ by Lemma~\ref{lemreg}.
\end{proof}

As a direct consequence we conclude the following statement.

\begin{corollary}\label{cor22}
	If $\eta_p^2-\tau_p^2\not=4$ and $\vert\eta_p\vert\not=\vert\tau_p\vert$ for all $p\in\{1,\dots,P\}$, then 
	\begin{multline}
		\dom A^{\eta,\tau}=\Big\{
		f\in H^{1,\Sigma}_\alpha(\Omega\setminus\Sigma,\CC^N)\cap H^{\half,\partial\Omega}_\alpha(\Omega\setminus\Sigma,\CC^N): \cB\gamma^{\partial\Omega}f=0\text{ on }\partial\Omega,\\
		\rmi (\alpha\cdot\nu^\Sigma)\big(\gamma^\Sigma_+f-\gamma^\Sigma_-f\big)+\dfrac{1}{2}(\eta I_N+\tau\beta)(\gamma^\Sigma_+f + \gamma^\Sigma_-f)=0 \text{ on }\Sigma
		\Big\}. \label{dom-pert2}
	\end{multline}	
\end{corollary}

Now, we turn our focus to the \textit{critical case}, i.e. when $\eta_p^2 - \tau_p^2 = 4$.
The following result on block operator matrices will be used below. It is an adapted formulation
of \cite[Theorem~2.2.18 and Theorem~2.4.7]{tretter}. Here and below, we understand the essential spectrum of a closed, but not necessarily self-adjoint operator $T$ to be the set of all $\lambda \in \CC$ such that $T - \lambda$ is not a Fredholm operator.

\begin{proposition}[Schur complement]\label{prop-schur} 
	Let $\cH$ be a Hilbert space and $A,B,C,D$ be densely defined and closable operators in $\cH$ such that $D$ has a bounded inverse,
	$\dom D\subset \dom B$, $\dom A\subset\dom C$, and $D^{-1}C$ is bounded on $\dom C$. Consider the linear operator $Q$ in $\cH\times\cH$,
	\[
	Q=\begin{pmatrix}
		A & B\\
		C & D
	\end{pmatrix},\quad
	\dom Q=\dom A\times \dom D,
	\]
	and its Schur complement $X:=A-BD^{-1}C$ defined on $\dom A$, which is viewed as a linear operator in $\cH$. Then, the following is true:
	\begin{enumerate}
		\item[(a)] $Q$ is closed (closable) if and only if $X$ is closed (closable, respectively), and in this case one has 
		\begin{equation*}
		 \overline Q=\begin{pmatrix}I & BD^{-1} \\ 0 & I\end{pmatrix}\begin{pmatrix}\overline X & 0 \\ 0 & D\end{pmatrix}\begin{pmatrix}I & 0 \\
		                                                                                                                  \overline{D^{-1}C} & I
		                                                                                                                 \end{pmatrix}.
		\end{equation*}
		\item[(b)] If $Q$ is closable, then $0\in\spece \overline Q$ if and only if $0\in\spece \overline {X}$.
	\end{enumerate}		
\end{proposition}

\begin{lemma}\label{lem28}
Let $p\in\{1,\dots,P\}$, assume that
$\eta_p^2-\tau_p^2=4$, and let $L_{p}$ be the self-adjoint operator in $L^2(\Sigma_p,\CC^N)$ in~\eqref{lql2}. Then
\begin{equation}
	\label{jpm}
J_{p}=\spece L_{p}=\begin{cases}
	\Big\{-\dfrac{\tau_p}{\eta_p}m\Big\}, & n=2,\\[\bigskipamount]
	-\dfrac{\tau_p}{\eta_p}m + \dfrac{1}{\eta_p}\Big[ -\dfrac{a_p}{2},\dfrac{a_p}{2}\Big], & n=3,
\end{cases}
\end{equation}
where, for $n=3$,  $a_p:=\max_{x\in\Sigma_p}\big|\kappa_1(x)-\kappa_2(x)\big|$ and $\kappa_1(x)$,  $\kappa_2(x)$ are the principal curvatures of $\Sigma_p$ in $x$.
\end{lemma}

\begin{proof}
Throughout the proof, if not declared differently, all operators will be considered to be defined on their maximal domains in the respective $L^2$-spaces.
Let $z\in\CC$ and note that $z\in \spece L_{p}$ if and only if $0\in\spece(z-L_{p})$. As the essential spectrum is invariant under compact perturbations,
for any $\Pi\in \Op\cS^{-\infty}_1(\Sigma_p)$, the condition $0\in\spece(z-L_{p})$
is equivalent to
\begin{equation}
	\label{ltilde}
0\in\spece\begin{pmatrix}
		\left(\dfrac{1}{\eta_p+\tau_p}\Lambda_p^2+(z+m)\right)\otimes I_{\frac{N}{2}} & \Lambda_pW^+_p\Lambda_p\\
		\Lambda_pW^-_p\Lambda_p & \left(\dfrac{1}{\eta_p-\tau_p}\Lambda_p^2+(z-m)+\Pi\right)\otimes I_{\frac{N}{2}}
	\end{pmatrix}.
\end{equation}
In particular, as $\Lambda_p$ has compact resolvent its spectrum is discrete and the same is true for $\Lambda_p^2$. Thus, one can choose
\[
\Pi:=\text{the orthogonal projector on }\ker \left( \frac{1}{\eta_p-\tau_p}\Lambda_p^2+(z-m)\right) \text{ in } L^2(\Sigma_p),
\]
which is a finite-dimensional projector with smooth integral kernel; cf. the proof of Lemma~\ref{lemma_Lambda_p}~(b) for a similar argument.  Since $\Pi$ is a spectral projection of the positive operator $\Lambda_p^2$, it commutes with $\Lambda_p$ and $\Lambda_p^{-1}$.  This choice of $\Pi$  guarantees that the bottom right block in \eqref{ltilde} has a bounded inverse and, as $\Lambda_p$ is elliptic,
\begin{equation} \label{mapping_properties_inverse}
  \left(\frac{1}{\eta_p-\tau_p}\Lambda_p^2+(z-m)+\Pi\right)^{-1} \in \Op \cS^{-1}_1(\Sigma_p).
\end{equation}
Here and below we write $a$ instead of $a\otimes I_{\frac{N}{2}}$ for better readability.  
Consider the block operator matrix
\begin{equation*}
\begin{split}
  \mathfrak{L}_p(z) &:= \begin{pmatrix}
		\dfrac{1}{\eta_p+\tau_p}\Lambda_p^2+(z+m) & \Lambda_pW^+_p\Lambda_p\\
		\Lambda_pW^-_p\Lambda_p & \dfrac{1}{\eta_p-\tau_p}\Lambda_p^2+(z-m)+\Pi
	\end{pmatrix}, \\
	\dom \mathfrak{L}_p(z) &:= H^1(\Sigma_p, \CC^N).
\end{split}
	\end{equation*}
 The operator $\mathfrak{L}_p(z)$ is a restriction of the closed operator 
 $$z-L_p+\begin{pmatrix} 0 & 0 \\ 0 & \Pi\end{pmatrix},$$ 
 where $L_p$ is defined on its maximal domain; in particular, $\mathfrak{L}_p(z)$  is closable. Since $H^1(\Sigma_p,\CC^N)$ is a core for $L_p$ by Proposition \ref{proposition_decomposition_Theta} and $\Pi$ is bounded, we conclude that 
$$\overline{\mathfrak{L}_p(z)} = z-L_p+\begin{pmatrix} 0 & 0 \\ 0 & \Pi\end{pmatrix}.$$
As $\Lambda_pW^-_p\Lambda_p \in \Op \cS^1_\nhalf(\Sigma_p)$, one obtains with~\eqref{mapping_properties_inverse} that
\begin{equation*}
  \left[ \frac{1}{\eta_p-\tau_p}\Lambda_p^2+(z-m)+\Pi\right]^{-1}\Lambda_p W^-_p\Lambda_p \in \Op \cS^0_\nhalf(\Sigma_p)
\end{equation*}
and hence, this operator is bounded in  $L^2(\Sigma_p,\CC^{\frac{N}{2}})$. In particular, for the block operator matrix  $\mathfrak{L}_p(z)$ all assumptions in Proposition~\ref{prop-schur} are fulfilled and because of  $\overline{\mathfrak{L}_p(z)} =z- L_p \mod \Op \cS^{-\infty}_N(\Sigma_p)$  one has that
\begin{equation}
	\label{schur001}
0\in\spece (z-L_{p})\iff 0\in \spece \overline{X_{p}(z)},
\end{equation}
where $X_{p}(z)$ is the  Schur complement associated with $\mathfrak{L}_p(z)$,
\begin{equation*}
  \begin{split}
    X_p(z) &= \frac{1}{\eta_p+\tau_p}\Lambda_p^2 +(z+m) -\Lambda_pW^+_p\Lambda_p\left[ \frac{1}{\eta_p-\tau_p}\Lambda_p^2+(z-m)+\Pi\right]^{-1}\Lambda_p W^-_p\Lambda_p, \\
    \dom X_p(z) &= H^1(\Sigma_p, \CC^\nhalf).
  \end{split}
\end{equation*}

In the following, we will consider $X_p(z)$ also as an element in $\Op \cS^1_\nhalf(\Sigma_p)$.  Since $\Pi$ commutes with $\Lambda_p$ we have 
\[
	\Lambda_p\left[ \frac{1}{\eta_p-\tau_p}\Lambda_p^2+(z-m)+\Pi\right]^{-1}\Lambda_p=(\eta_p-\tau_p)[I+T]^{-1}
\]	
with
\begin{equation*}
 T=(\eta_p-\tau_p)(z-m)\Lambda_p^{-2} +(\eta_p-\tau_p)\Pi\Lambda_p^{-2}\in \Op\cS^{-1}_{\frac{N}{2}}(\Sigma_p)
\end{equation*}
since $\Lambda_p^{-2}\in \Op\cS^{-1}_{\frac{N}{2}}(\Sigma_p)$ and $\Pi\Lambda_p^{-2}\in\Op\cS^{-\infty}_{\frac{N}{2}}(\Sigma_p)$.
Using $(I+T)^{-1}=I-T +T^2(1+T)^{-1}$, $T^2\in \Op\cS^{-2}_{\frac{N}{2}}(\Sigma_p)$, and $(1+T)^{-1}\in \Op\cS^{0}_{\frac{N}{2}}(\Sigma_p)$ we 
arrive at
\begin{align*}
	\Lambda_p\!\left[ \frac{1}{\eta_p-\tau_p}\Lambda_p^2+(z-m) + \Pi\right]^{-1}\!\!\Lambda_p&=(\eta_p-\tau_p)-(\eta_p-\tau_p)^2(z-m)\Lambda_p^{-2}	\mod \Op\cS^{-2}_{\frac{N}{2}}(\Sigma_p).
\end{align*}
Plugging this into the expression for the Schur complement, one gets
\begin{align*}
X_{p}(z)=\dfrac{1}{\eta_p+\tau_p}\Lambda_p^2 +(z+m)-(\eta_p-\tau_p)\Lambda_p W^+_pW^-_p\Lambda_p&\\
 +(\eta_p-\tau_p)^2(z-m)\Lambda_pW^+_p\Lambda_p^{-2}W^-_p\Lambda_p	&\mod \Op \cS^{-1}_{\frac{N}{2}}(\Sigma_p).
\end{align*}
As $\Lambda_p$ is a scalar operator, we have
$W^\pm_p \Lambda_p=\Lambda_p W^\pm_p \mod \Op \cS^{-\half}_{\frac{N}{2}}(\Sigma_p)$ (see Appendix~\ref{appendix_pseudos}).
Hence,
\begin{align*}
\Lambda_pW^+_p\Lambda_p^{-2}W^-_p\Lambda_p&=W^+_p\Lambda_p^{-1}W^-_p\Lambda_p \mod \Op\cS^{-1}_{\frac{N}{2}}(\Sigma_p)\\
&=W^+_pW^-_p \mod \Op\cS^{-1}_{\frac{N}{2}}(\Sigma_p)\\
&=\frac{1}{4}\mod \Op\cS^{-1}_{\frac{N}{2}}(\Sigma_p),
\end{align*}
where in the last step Corollary~\ref{cor-ww14} was employed.
Therefore, we conclude
\begin{equation*}
	X_{p}(z)=\Tilde X_{p}(z) \mod \Op \cS^{-1}_{\frac{N}{2}}(\Sigma_p),
\end{equation*}
where
\begin{equation*}
\begin{split}
	\Tilde X_{p}(z):&=\dfrac{1}{\eta_p+\tau_p}\Lambda_p^2 +(z+m)-(\eta_p-\tau_p)\Lambda_p W^+_pW^-_p\Lambda_p+\dfrac{(\eta_p-\tau_p)^2(z-m)}{4}\\
	&= (z+m)+\dfrac{(\eta_p-\tau_p)^2(z-m)}{4} +(\eta_p-\tau_p)\Lambda_p\left(\frac{1}{4}-W^+_pW^-_p\right)\Lambda_p;
\end{split}
\end{equation*}
in the last simplification the condition $\eta_p^2 - \tau_p^2 = 4$ was used and $\Tilde X_{p}(z)$ is viewed as the maximal realization of the associated pseudodifferential operator in $L^2(\Sigma_p, \CC^\nhalf)$ defined on $\dom \Tilde X_{p}(z) = \dom \overline{X_p(z)}$. 
In particular, as operators in $\Op \cS^{-1}_{\frac{N}{2}}(\Sigma_p)$ are bounded from $L^2(\Sigma_p, \CC^\nhalf)$ to $H^1(\Sigma_p, \CC^\nhalf)$ and hence, by Rellich's embedding theorem, compact in $L^2(\Sigma_p, \CC^\nhalf)$, one has
\begin{equation}
	\label{spectilde0}
0\in \spece \overline{X_{p}(z)}\iff 0\in\spece \Tilde X_{p}(z).
\end{equation}
A simple algebra using again $\eta_p^2 - \tau_p^2 = 4$ gives
\[
	(z+m)+\dfrac{(\eta_p-\tau_p)^2(z-m)}{4}
	=\dfrac{(\eta_p z+\tau_p m)}{2}(\eta_p-\tau_p),
\]
and thus,
\[
0\in \spece \Tilde X_{p}(z) \iff \eta_p z +\tau_p m\in \spece \left(2\Lambda_p\left(W^+_pW^-_p-\frac{1}{4}\right)\Lambda_p\right).
\]
By comparing with \eqref{spectilde0} and then with \eqref{schur001} we arrive at
\begin{equation}
	\label{schur111}
	\spece L_{p} = -\dfrac{\tau_p}{\eta_p}m + \dfrac{2}{\eta_p}\spece R
\quad\text{with} \quad
  R:= \Lambda_p\left(W^+_pW^-_p-\frac{1}{4}\right)\Lambda_p;
\end{equation}
observe that $\eta_p \not=0$ as we have $\eta_p^2 - \tau_p^2 = 4$.
Now we consider the cases $n=2$ and $n=3$ separately.

(a) For $n=2$ one has by Lemma~\ref{lem26}~(a) that $W^+_pW^-_p=\frac{1}{4}\mod \Op \cS^{-\infty}_1(\Sigma_p)$.
It follows that $R\in \Op \cS^{-\infty}_1(\Sigma_p)$. Since all operators in $\Op \cS^{-\infty}_1(\Sigma_p)$ are compact in $L^2(\Sigma_p)$, we conclude $\spece R =\{0\}$, and \eqref{schur111} yields
\[
\spece L_{p}=\left\{ -\frac{\tau_p}{\eta_p}m\right\}.
\]

(b) Consider now $n=3$. Recall that in this case $W^+_p=W^-_p=:W_p$. By Lemma~\ref{lem26}~(b) we have
\begin{align*}
R&= \Lambda_p\left(W_p^2-\frac{1}{4}\right)\Lambda_p=\Lambda_p\Big((\sigma\cdot\nu^\Sigma)^2W_p^2 +(\sigma\cdot\nu^\Sigma)W_p(\sigma\cdot\nu^\Sigma)W_p\Big)\Lambda_p\\
&=\Lambda_p(\sigma\cdot\nu^\Sigma)\Big[(\sigma\cdot\nu^\Sigma)W_p+W_p(\sigma\cdot\nu^\Sigma)\Big]W_p \Lambda_p\\
&=\Lambda_p (\sigma\cdot\nu^\Sigma)Z_pW_p\Lambda_p,
\end{align*}
where in the last step the definition of $Z_p$ from \eqref{eq-zp} was used.

Let $x\in \Sigma_p$ be fixed. Using Lemma~\ref{lem19}~(b), Lemma~\ref{lem26}~(b), and Lemma~\ref{lemma_Lambda_p}~(c) one gets for the principal symbol of $R$ with respect
to normal coordinates near $x$ (see Definition~\ref{def_special_parametrization}~(b))
\begin{align*}
	\sigma_R(x,\xi)&=\sigma_{\Lambda_p}(x,\xi)\sigma_{\sigma\cdot\nu^\Sigma}(x,\xi)\sigma_{Z_p}(x,\xi)\sigma_{W_p}(x,\xi)\sigma_{\Lambda_p}(x,\xi)\\
	&=\sqrt{2|\xi|}\big(\sigma\cdot \nu^\Sigma(x)\big)\frac{1}{2}\big(\kappa_1(x)-\kappa_2(x)\big)(\sigma\cdot \nu^\Sigma(x))\,	\frac{\xi_1 \xi_2}{|\xi|^3}\frac{1}{2}\sigma\cdot \left(\frac{\xi_1 e_1+\xi_2 e_2}{|\xi|}\right)\sqrt{2|\xi|}\\
	&=\dfrac{1}{2}\big(\kappa_1(x)-\kappa_2(x)\big)\dfrac{\xi_1\xi_2}{|\xi|^3}\sigma\cdot (\xi_1 e_1+\xi_2 e_2).
\end{align*}
Recall that $\spec (\sigma\cdot a)=\{-|a|,|a|\}$ for any $a\in\RR^3$. Thus, using $|\xi_1 e_1+\xi_2 e_2|=|\xi|$ we arrive at 
\[
\spec \sigma_R(x,\xi)=\left\{ -\dfrac{\big|\kappa_1(x)-\kappa_2(x)\big|}{2}\dfrac{|\xi_1\xi_2|}{|\xi|^2},\dfrac{\big|\kappa_1(x)-\kappa_2(x)\big|}{2}\dfrac{|\xi_1\xi_2|}{|\xi|^2}\right\}.
\]
Since the range of the function $(\xi_1,\xi_2)\mapsto  \tfrac{|\xi_1\xi_2|}{|\xi|^2}$, $\xi\not=0$, is $[0,\tfrac{1}{2}]$
the usual rule for zeroth order pseudodifferential operators (see, e.g., \cite[Lemma~11]{BP24}) gives
\[
\spece R=\bigcup_{x\in\Sigma_p}\bigcup_{\xi\ne 0}\spec \sigma_{R}(x,\xi)=\left[ -\dfrac{a_p}{4},\dfrac{a_p}{4}\right]
\]
with $a_p = \max_{x\in\Sigma_p}\big|\kappa_1(x)-\kappa_2(x)\big|$.
The substitution into \eqref{schur111} completes the proof.
\end{proof}

By combining Proposition~\ref{proposition_decomposition_Theta} with Lemmas \ref{lem21} and~\ref{lem28} we arrive at the following intermediate result.
\begin{corollary}\label{cor29}
If $\vert\eta_p\vert\not=\vert\tau_p\vert$ for all $p\in\{1,\dots,P\}$, then
	\begin{equation}
		\label{jpm2}
	\spece A^{\eta,\tau} \cap \res A=\bigg(\bigcup_{p:\ \eta_p^2-\tau_p^2=4}J_{p}\bigg) \cap \res A
	\end{equation}
	 with the sets $J_{p}$ given by \eqref{jpm}.
\end{corollary}

In order to complete the analysis of the essential spectrum, we have to show that \eqref{jpm2} remains valid if the set $\res A$ is removed on both sides.
More precisely, we have to
take care of the following remaining aspects:
\begin{itemize}
	\item[(I)] it is still unclear if an eigenvalue in $\spec A$ lying outside of $J_{p}$ can belong  to $\spece A^{\eta,\tau}$;
	\item[(II)] even if $\eta_p^2-\tau^2_p=4$ for some $p$, the right-hand side of \eqref{jpm2} can be empty, if $J_{p}$ degenerates to a point which coincides with an eigenvalue of $A$.
\end{itemize}
The open point (I) will be answered in Lemma~\ref{lem31} and (II) will be subject of Lemma~\ref{lem33}.

Let us first recall the following abstract result \cite[Theorem XIII.13]{RS4}.

\begin{proposition}[Meromorphic Fredholm theorem]\label{lem-fred}
Let $\cH$ be a Hilbert space and let $U\subset\CC$ be a connected open set.
Let $Z\subset U$ be a subset having no accumulation points in $U$
and $K:U\setminus Z\to\cL(\cH)$ be an analytic operator function
which is meromorphic in $U$ and such that
\begin{itemize}
	\item the operator $K(z)$ is compact for any $z\in U\setminus Z$,
	\item at each $z_0\in Z$ all coefficients of the principal part of the Laurent series of $K(z)$ at  $z_0$ are finite-dimensional operators.
\end{itemize}	
Then either
\begin{itemize}
	\item[(a)] $0\in\spec \big(I+K(z)\big)$ for all $z\in U\setminus Z$, 
\end{itemize}		
or
\begin{itemize}
	\item[(b)] there is a subset $Z'\subset U$ without accumulation points in $U$
	such that the bounded inverse $\big(I+K(z)\big)^{-1}$ exists for all $z\in U\setminus(Z\cup Z')$
	and extends to an operator function analytic in $U\setminus Z'$ and meromorphic in $U$,
	and for each 
	$z_0\in Z'$ all coefficients of the principal part of the Laurent series of $\big(I+K(z)\big)^{-1}$ at  $z_0$ are finite-dimensional operators.
\end{itemize}		
\end{proposition}

\begin{lemma}\label{lem31}
If $\vert\eta_p\vert\not=\vert\tau_p\vert$ for all $p\in\{1,\dots,P\}$, then
\begin{equation}
	\label{jpm3}
	\spece A^{\eta,\tau} \subset \bigcup_{p:\ \eta_p^2-\tau_p^2=4}J_{p}
\end{equation}
 with the sets $J_{p}$ given by \eqref{jpm}.
\end{lemma}

\begin{proof}
Denote the set on the right-hand side of \eqref{jpm3} by $J$ and note that the sought assertion
is a direct consequence of the following claim:
\begin{equation}
	\label{claim555}
\parbox{100mm}{The resolvent of $A^{\eta,\tau}$ is meromorphic in $\CC\setminus J$ and the 
coefficients of the principal part of the Laurent series of $(A^{\eta,\tau}-z)^{-1}$ at each pole are finite-dimensional operators.}
\end{equation}
In fact, \eqref{claim555} means that the spectrum of $A^{\eta,\tau}$ in $\CC\setminus J$ is purely discrete. Note also that due to the self-adjointness 
of $A^{\eta,\tau}$ each pole of the resolvent is necessarily of first order.

	Recall that we have the resolvent formula 
	\begin{equation}
		\label{res-krein2}
		(A^{\eta,\tau}-z)^{-1}=(A-z)^{-1}+\Gamma_{z}(\Theta^{\eta,\tau}-M_z)^{-1}\Gamma^*_{\Bar z},\quad z\in\res A^{\eta,\tau}\cap\res A,
	\end{equation}
	with $\Gamma_{z}$, $M_{z}$, and $\Theta^{\eta,\tau}$ given by \eqref{eq-gbz}, \eqref{Weyl_function}, and~\eqref{def_Theta}, respectively. 
	By abstract considerations the maps
	\[
	z\mapsto (A-z)^{-1},\quad
	z\mapsto \Gamma_{z},
	\quad
	z\mapsto M_{z}
	\]
	are meromorphic in $\CC$ having only first-order poles with finite-dimensional residues, see Remark~\ref{rem-qq}.
	Furthermore, in view of the representation \eqref{decomposition_Theta} we have
	\[
	\Theta^{\eta,\tau}-M_z=L-z +T(z),
	\]
	where $L:=\bigoplus_{p=1}^P L_{p}$ is self-adjoint in $L^2(\Sigma, \CC^N)$ and
	\[
	T:\ \res A\ni z\mapsto T(z)\in \cL\big(L^2(\Sigma,\CC^N)\big)
	\]
	is a meromorphic operator-valued function having only first-order poles with finite-dimensional residues and such that  for any $z$ the value $T(z) \in \Op \cS^{-1}_N(\Sigma)$ gives rise to a bounded operator from $L^2(\Sigma, \CC^N)$ to $H^1(\Sigma, \CC^N)$ and thus, it is a compact operator in $L^2(\Sigma, \CC^N)$. Due to Lemma~\ref{lem21} and Lemma~\ref{lem28} one has $\spece L=J$. Therefore, the resolvent $z\mapsto (L-z)^{-1}$ is meromorphic in $\CC\setminus J$ with only first-order poles and finite-dimensional residues. We represent
	\[
	\Theta^{\eta,\tau}-M_z=\big(I+K(z)\big) (L-z),\quad K(z):=T(z)(L-z)^{-1},
	\]
	and note that all $K(z)$ are compact operators and satisfy the assumptions of the meromorphic Fredholm theorem (Proposition~\ref{lem-fred}) on the domain $U:=\CC\setminus J$. The case  (a) in Proposition~\ref{lem-fred} can be excluded as the operators $\Theta^{\eta,\tau}-M_{z}$ and $L-z$ are boundedly invertible for $z\in\CC\setminus\RR$; for $\Theta^{\eta,\tau}-M_z$ this follows from Theorem~\ref{thm17}~(a), while for $L-z$ this is clear due to the self-adjointness of $L$.	 It follows that 	$z\mapsto \big(I+K(z)\big)^{-1}$  is meromorphic
	in $\CC\setminus J$ with finite-dimensional coefficients in the principal parts of the Laurent series at the poles. Then
	\[
	(\Theta^{\eta,\tau}-M_z)^{-1}=(L-z)^{-1}\big(I+K(z)\big)^{-1}
	\]
	has the same property, and the substitution into the resolvent formula \eqref{res-krein2}
	proves the required claim \eqref{claim555}.	
\end{proof}

Recall the following abstract convergence result \cite[Satz 9.24]{weid}.

\begin{proposition}\label{lem32}
	Let $T_j$ and $T$ be self-adjoint operators in a Hilbert space such that $T_j\rightarrow T$, $j\to\infty$, in the norm resolvent sense. 
	If $\lambda_j\in\spece T_j$ with $\lambda_j\rightarrow\lambda$, $j\to\infty$, then $\lambda\in\spece T$.
\end{proposition}

\begin{lemma}\label{lem33}
	Let $p\in\{1,\dots,P\}$ and assume that $\eta^2_p-\tau^2_p=4$. Then $-\tfrac{\tau_p}{\eta_p}m\in\spece A^{\eta,\tau}$.
\end{lemma}	

\begin{proof}
In the case that $-\tfrac{\tau_p}{\eta_p}m\in\res A$ the claim follows from  Corollary~\ref{cor29}, since $-\frac{\tau_p}{\eta_p}\in J_p$ by Lemma~\ref{lem28}. 
The objective of the next steps (a)-(c) is to cover the case $-\tfrac{\tau_p}{\eta_p}m\in\spec A$.

(a) In the following, we denote by $L_p$ the maximal realization as an operator in $L^2(\Sigma_p, \CC^N)$ associated with the pseudodifferential expression defined in~\eqref{lql2}. We start by proving the following auxiliary statement:
\begin{equation}
	\label{claim00}
	\parbox{100mm}{The operators $L_{p}$ depend continuously on $(\eta_p,\tau_p)$ on the set
	$\Xi:=\{(\eta_p,\tau_p):\eta_p^2-\tau_p^2=4\}$ with respect to the norm resolvent topology.}
\end{equation}
It is sufficient to consider the case $m=0$ (as adding a fixed bounded operator does not influence the norm resolvent convergence).
Consider the block operator matrix (with $a\otimes I$ shortened as $a$)
\begin{equation*}
\begin{split}
  \mathfrak{L}_p &:= \begin{pmatrix}
		\dfrac{1}{\eta_p+\tau_p}\Lambda_p^2+\rmi & \Lambda_pW^+_p\Lambda_p\\
		\Lambda_pW^-_p\Lambda_p & \dfrac{1}{\eta_p-\tau_p}\Lambda_p^2+\rmi
	\end{pmatrix}, \quad
	\dom \mathfrak{L}_p = H^1(\Sigma_p, \CC^N).
\end{split}
	\end{equation*}
As in the proof of Lemma~\ref{lem28} (see the considerations on the operator $\mathfrak{L}_p(z)$ there) one sees that $\mathfrak{L}_p$ is closable, $\overline{\mathfrak{L}_p} = \rmi-L_p$ (since $H^1(\Sigma_p,\CC^N)$ is a core for $L_p$ by Proposition~\ref{proposition_decomposition_Theta}), and $\mathfrak{L}_p$ satisfies all assumptions in Proposition~\ref{prop-schur}. For better readability we denote the right lower entry in $\mathfrak L_p$ by $D_p$. Then we conclude
\begin{equation}\label{auchdasbitte}
\overline{\mathfrak{L}_p} =
(\rmi-L_{p})=\begin{pmatrix}
	I & \Lambda_p W^+_p \Lambda_p D_p^{-1}\\
	0 & I
\end{pmatrix}
\begin{pmatrix}
	\overline{X(\eta_p,\tau_p)} & 0\\
	0 & D_p
\end{pmatrix} \begin{pmatrix}
	I & 0\\ 
	D_p^{-1}\Lambda_pW^-_p\Lambda_p & I
\end{pmatrix}
,
\end{equation}
where $X(\eta_p,\tau_p)$ is the Schur complement associated with $\mathfrak{L}_p$ defined as in Proposition~\ref{prop-schur},
\begin{align*}
	X(\eta_p,\tau_p)&=
	\dfrac{1}{\eta_p+\tau_p}\Lambda_p^2+\rmi-\Lambda_pW^+_p\Lambda_p D_p^{-1}\Lambda_pW^-_p\Lambda_p,\\
	\dom X(\eta_p,\tau_p) &= H^1(\Sigma_p, \CC^\nhalf).
\end{align*}
Note that $X(\eta_p,\tau_p)\in \Op \cS^0_\nhalf(\Sigma_p)$, since its principal symbol of order one is 
\begin{equation*}
  2|\xi|\left(\frac{1}{\eta_p+\tau_p}-\frac{\eta_p-\tau_p}{4}\right)=0
\end{equation*}
by Corollary~\ref{cor-ww14} and $\eta_p^2-\tau_p^2=4$; cf. the considerations in the proof of Lemma~\ref{lem28} leading to~\eqref{schur111}. Hence,  $\overline{X(\eta_p,\tau_p)}$ is the bounded extension of $X(\eta_p,\tau_p)$ to $L^2(\Sigma_p,\CC^\nhalf)$. 
With a computation one sees that 
\begin{align*}
	\overline{X(\eta_p,\tau_p)}&=
	\dfrac{1}{\eta_p+\tau_p}\Lambda_p^2+\rmi-\Lambda_pW^+_p\Lambda_p \left(\dfrac{1}{\eta_p-\tau_p}\Lambda_p^2+\rmi\right)^{-1}\Lambda_pW^-_p\Lambda_p\\
	&=\rmi+\dfrac{\eta_p-\tau_p}{4}\Lambda_p^2-(\eta_p-\tau_p)\Lambda_pW^+_p\Lambda_p \big(\Lambda_p^2+\rmi(\eta_p-\tau_p)\big)^{-1}\Lambda_pW^-_p\Lambda_p\\
	&=\rmi+\dfrac{\eta_p-\tau_p}{4}\Lambda_p^2-(\eta_p-\tau_p)\Lambda_pW^+_p\Big[I- \rmi(\eta_p-\tau_p) \big(\Lambda_p^2+\rmi(\eta_p-\tau_p)\big)^{-1}\Big]W^-_p\Lambda_p\\
	&=\rmi +(\eta_p-\tau_p)\Lambda_p\left( \dfrac{1}{4}-W^+_p W^-_p\right)\Lambda_p \\
	&\qquad +\rmi (\eta_p-\tau_p)^2 \Lambda_pW^+_p \big(\Lambda_p^2+\rmi(\eta_p-\tau_p)\big)^{-1}W^-_p\Lambda_p \\
	&=\rmi +(\eta_p-\tau_p) K_1 +\rmi (\eta_p-\tau_p)^2 K_2(\eta_p,\tau_p),
\end{align*}
where we have set
\begin{equation*}
  K_1 := \Lambda_p\left( \dfrac{1}{4}-W^+_p W^-_p\right)\Lambda_p
\end{equation*}
and
\begin{equation*}
  K_2(\eta_p,\tau_p) := \Lambda_pW^+_p \big(\Lambda_p^2+\rmi(\eta_p-\tau_p)\big)^{-1}W^-_p\Lambda_p.
\end{equation*}
We first remark that $\overline{X(\eta_p,\tau_p)}$ always has a bounded inverse, as with $W_p^+ = (W_p^-)^*$ and 
$$\Im D_p^{-1}=\Im \left(\dfrac{1}{\eta_p-\tau_p}\Lambda_p^2+\rmi\right)^{-1} \leq 0$$ one finds that
\begin{equation*}
  \begin{split}
    \Im \overline{X(\eta_p,\tau_p)} &= \Im \rmi- \Im \left[\Lambda_pW^+_p\Lambda_p D_p^{-1}\Lambda_pW^-_p\Lambda_p\right]\\
    &=1- \Lambda_pW^+_p\Lambda_p \left[\Im D_p^{-1}\right]\Lambda_pW^-_p\Lambda_p
    \ge 1
  \end{split}
\end{equation*}
in the sense of quadratic forms. In particular, zero does not belong to the closure of the numerical range of $\overline{X(\eta_p,\tau_p)}$ and thus, as $\overline{X(\eta_p,\tau_p)}$ is bounded, it is also boundedly invertible. The fixed operator $K_1$ belongs to $\Op \cS^0_{\frac{N}{2}}(\Sigma_p)$ by Corollary~\ref{cor-ww14}, hence, it is bounded in $L^2(\Sigma_p,\CC^\frac{N}{2})$. On the other hand, for any $s\in \RR$ one has fixed bounded maps
\[
\Lambda_p: H^s(\Sigma_p)\to H^{s-\half}(\Sigma_p),
\quad
W^\pm_p: H^s(\Sigma_p,\CC^\frac{N}{2})\to H^{s}(\Sigma_p,\CC^\frac{N}{2}),
\]
and the continuous map
\[
\RR\ni\lambda\mapsto \big(\Lambda_p^2+\rmi \lambda\big)^{-1}\in \cL\big(H^{s}(\Sigma_p),H^{s+1}(\Sigma_p)\big),
\]
which shows that the operators $D_p^{-1}$, $\Lambda_pW^+_p\Lambda_p D_p^{-1}$, $D_p^{-1}\Lambda_pW^-_p\Lambda_p$, and $K_2(\eta_p,\tau_p)$ 
depend continuously on $(\eta_p,\tau_p)\in\Xi$
in $\cL(L^2(\Sigma_p,\CC^\frac{N}{2}))$. Then, the same applies
to $\overline{X(\eta_p,\tau_p)}$ and hence also 
to $\overline{X(\eta_p,\tau_p)}^{-1}$. 
To complete the argument note that \eqref{auchdasbitte} implies
\begin{equation}\label{lokal1}
(\rmi-L_{p})^{-1}
=\begin{pmatrix}
	I & 0\\ 
	-D_p^{-1}\Lambda_pW^-_p\Lambda_p & I
\end{pmatrix}\\
\begin{pmatrix}
	\overline{X(\eta_p,\tau_p)}^{-1} & 0\\
	0 & D_p^{-1}
\end{pmatrix}
\begin{pmatrix}
	I & -\Lambda_p W^+_p \Lambda_p D_p^{-1}\\
	0 & I
\end{pmatrix},
\end{equation}
and the above considerations 
show that 
all blocks in \eqref{lokal1} depend  continuously on $(\eta_p,\tau_p)\in\Xi$
in $\cL\big(L^2(\Sigma_p,\CC^\frac{N}{2})\big)$. This yields 
the claim \eqref{claim00}.

(b) Assume that $m\ne 0$ and $-\tfrac{\tau_p}{\eta_p}m\in\spec A$.
For $(\eta'_p,\tau'_p)\in\RR^2$ with $(\eta'_p)^2-(\tau'_p)^2=4$ and $(\eta'_p,\tau'_p)\ne(\eta_p,\tau_p)$
define locally constant functions $\Tilde\eta,\Tilde\tau:\Sigma\to\RR$ by
\begin{align*}
\Tilde \eta&=\eta \text{ on }\Sigma\setminus\Sigma_p, & \Tilde \tau&=\tau \text{ on }\Sigma\setminus\Sigma_p,\\
\Tilde \eta&=\eta'_p \text{ on }\Sigma_p, & \Tilde \tau&=\tau'_p \text{ on }\Sigma_p,
\end{align*}
and consider the associated operator $A^{\Tilde\eta,\Tilde\tau}$. 
It follows from the definitions \eqref{def_Theta} and \eqref{lql2} and from Proposition~\ref{proposition_decomposition_Theta} (with $z=\zeta$, so that $M_\zeta=0$) that 
\[
\Theta^{\eta,\tau}=\bigoplus_{p=1}^P L_{p}+ B
\]
holds
with some bounded operator $B$ independent of $(\eta,\tau)$. Hence, part (a) implies that $\Theta^{\Tilde\eta,\Tilde\tau}$ converges to $\Theta^{\eta,\tau}$ in the norm resolvent sense
for $(\eta'_p,\tau'_p)\to(\eta_p,\tau_p)$.
Lemma~\ref{lemma12} implies the norm resolvent convergence of $A^{\Tilde\eta,\Tilde\tau}$
to $A^{\eta,\tau}$ in the same regime. Note that the map $(\eta'_p,\tau'_p)\mapsto \frac{\tau'_p}{\eta'_p}$ is injective on each of the two connected component of $\Xi$ defined in~\eqref{claim00}. As $m\neq0$, it follows that $\lambda' := -m\frac{\tau'_p}{\eta'_p}\neq -m\frac{\tau_p}{\eta_p}$ for all $(\eta'_p,\tau'_p)\neq(\eta_p,\tau_p)$ sufficiently close to $(\eta_p,\tau_p)$, and since $\spec A$ consists of isolated points, also $\lambda'\notin\spec A$ for all such 
 $(\eta'_p,\tau'_p)$. As in the first step of the proof, Corollary~\ref{cor29} applied to $A^{\Tilde\eta,\Tilde\tau}$ shows that  $\lambda'\in \spece A^{\Tilde\eta,\Tilde\tau}$ for all $(\eta'_p,\tau'_p)$ sufficiently close to $(\eta_p,\tau_p)$. Then Proposition~\ref{lem32} yields
\[
-\dfrac{\tau_p}{\eta_p}m=\lim_{(\eta'_p,\tau'_p)\to (\eta_p,\tau_p)}\lambda'\in \spece A^{\eta,\tau}.
\]

(c) Consider the remaining case $m=0$ and $-m\frac{\tau_p}{\eta_p} = 0\in\spec A$. We have obviously
 $A^{\eta,\tau} + M \beta\rightarrow A^{\eta,\tau}$, as $M\to0$, in the norm resolvent sense.
For all $M\ne 0$ we have $-\tfrac{\tau_p}{\eta_p}M\in \spece (A^{\eta,\tau} + M \beta)$ by the first step of the proof if $-\tfrac{\tau_p}{\eta_p}M\in \res (A+ M \beta)$, and by (b) if $-\tfrac{\tau_p}{\eta_p}M\in \spec (A+ M \beta)$; with the help of Proposition~\ref{lem32} we obtain
\[
 0=\lim_{M\to0,\, M\ne 0}-\dfrac{\tau_p}{\eta_p}M\in \spece A^{\eta,\tau}.\qedhere
\]
\end{proof}

We summarize all the constructions to arrive at the main result of the paper:

\begin{theorem}[Essential spectrum of shell interactions]\label{thm37}
Assume that $|\eta_p|\ne|\tau_p|$ for all $p\in\{1,\dots,P\}$. Then
	\begin{equation}
		\label{eqfinal}
	\spece A^{\eta,\tau}=\bigcup_{p=1}^P J_{p},
	\end{equation}
	where
\begin{equation*}
	J_{p}:=\begin{cases}
		\emptyset, & \eta_p^2-\tau_p^2\ne 4,\\[\medskipamount]
		\Big\{-\dfrac{\tau_p}{\eta_p}m\Big\}, & n=2,\ \eta_p^2-\tau_p^2= 4,\\[\bigskipamount]
		-\dfrac{\tau_p}{\eta_p}m + \dfrac{1}{\eta_p}\Big[ -\dfrac{a_p}{2},\dfrac{a_p}{2}\Big],& n=3,\ \eta_p^2-\tau_p^2= 4,
	\end{cases}
\end{equation*}
where $a_p:=\max_{x\in\Sigma_p}\big|\kappa_1(x)-\kappa_2(x)\big|$ and $\kappa_1(x)$,  $\kappa_2(x)$ are 
the principal curvatures of $\Sigma_p$ in $x$.
\end{theorem}

\begin{proof}
The inclusion $\subset$ in \eqref{eqfinal} holds by Lemma~\ref{lem31}, and
it remains to show the reverse inclusion.

Let $p\in\{1,\dots,P\}$ with $\eta_p^2-\tau_p^2=4$. If $J_{p}=\{-\tfrac{\tau_p}{\eta_p}m\}$, then $J_{p}\subset \spece A^{\eta,\tau}$ by Lemma~\ref{lem33}.
Assume now that $J_{p}$ has a non-empty interior, then by Corollary~\ref{cor29} it holds that
$J'_{p}:=J_{p}\cap \res A\subset\spece A^{\eta,\tau}$.
The essential spectrum of $A^{\eta,\tau}$
is a closed subset of $\RR$ and the spectrum of $A$ is discrete, therefore $J'_{p}$ is dense in $J_p$ and we conclude $J_{p}= \overline{J'_{p}}\subset\spece A^{\eta,\tau}$. 
\end{proof}	

\begin{remark}\label{remlast}
The case $|\eta_p|=|\tau_p|$ was only excluded for technical reasons, namely in order to ensure the invertibility of $\eta I_N+\tau\beta$ in the proof of Lemma~\ref{lem16}. It can still be covered by introducing additional projectors similarly to \cite[Proposition~4.3]{BHOP20},
and one can show that the result of Theorem~\ref{thm37} still holds with $J_{p}=\emptyset$ for  $|\eta_p|=|\tau_p|$.
\end{remark}
\begin{remark}\label{remlast2}  Throughout this paper the interaction strengths are constant on each connected component $\Sigma_p$. If $\eta$ and $\tau$ are allowed to vary along $\Sigma$, then a non-degenerate interval of essential spectrum can appear already for $n=2$. This was shown in \cite{BBP} for the Landau-Dirac operator in $\RR^2$ and smooth coefficients $\eta$ and $\tau$ on a closed curve $\Sigma$ satisfying $\eta^2-\tau^2=4$ on $\Sigma$: in this case the essential spectrum consists of the Landau-Dirac levels and the range of the function $-m\tfrac{\tau}{\eta}$ on $\Sigma$. It is also pointed out in \cite{BBP}  that an adaptation of the present approach should give the analogous result on bounded domains. In fact, if $n=2$ and $\eta,\tau\in C^\infty(\Sigma_p)$ with $\eta^2_p-\tau^2_p=4$ on $\Sigma_p$, then in the proof of Lemma~\ref{lem28} the constants $\tfrac{1}{\eta_p\pm\tau_p}$ have to be replaced by the corresponding multiplication operators, which commute with $W_p^\pm$ up to operators with smooth integral kernels, and one find that $\overline{X_p(z)}$ coincides, up to a compact operator, with the operator of multiplication by the function $\tfrac{1}{2}(\eta-\tau)(\eta z+\tau m)$. This indicates that for $n=2$ the set $J_p$ in Theorem~\ref{thm37} has to be replaced by the range of $-m\tfrac{\tau}{\eta}$ on $\Sigma_p$, which is a non-degenerate interval if $m\neq 0$ and $\eta,\tau$ are not constant on $\Sigma_p$. For $n=3$ additional work is needed, as the commutators of $W_p$ with non-constant multiplication operators are not smoothing.
\end{remark}

\begin{appendix}

\section{Notations on pseudodifferential operators and proof of Lemma~\ref{lemma_expansion1}} \label{appendix_pseudos}

In the following, let $M \subset \mathbb{R}^n$ be a smooth, closed, and compact hypersurface in $\mathbb{R}^n$, i.e. a curve in $\mathbb{R}^2$ or a surface in $\mathbb{R}^3$.
We use the following convention for the Fourier transform in $\RR^r$. For $f$ in the Schwartz class, it is given by
\[
\widehat{f}(\xi)=\int_{\RR^r}e^{-\rmi\langle x,\xi\rangle} f(x)\,\dd x.
\]
We denote by $\Op\cS^k_q(M)$ the class of classical pseudodifferential operators of order $\le k$ acting on sections of  $M\times\CC^q$ defined as in \cite[Chapter~7]{T1}; we refer to the latter reference for basic definitions
and properties of pseudodifferential operators on manifolds. Moreover, we set $\Op \cS^{-\infty}_q(M):=\bigcap_{\ell\in\RR}\Op \cS^{\ell}_q(M)$. For two pseudodifferential operators $A, B \in \Op \cS^k_q(M)$ and $\ell < k$ we write $A = B \mod \Op \cS^\ell_q(M)$, if $A - B \in \Op \cS^\ell_q(M)$. Furthermore, we will write $\sigma_{A}(x,\xi)$ for the principal symbol of $A \in \Op \cS^k_q(M)$.

In the following we recall some properties of pseudodifferential operators on $M$ that will be used frequently in the paper. Let $A \in \Op \cS^k_q(M)$. Then, for any $s \in \mathbb{R}$ the map $A$ extends uniquely to a bounded operator $A: H^s(M, \CC^q) \rightarrow H^{s-k}(M, \CC^q)$, see \cite[Chapter~7, Proposition~5.5 and Section~10]{T1}. In particular, if $k < 0$, then $A$ gives rise to a compact operator in $L^2(M, \CC^q)$ by Rellich's embedding theorem. If $k>0$, then it is a consequence of Friedrichs' lemma, see e.g. \cite[Chapter~II, Proposition~7.4]{taylor1}, that  $A: H^k(M, \CC^q) \rightarrow L^2(M, \CC^q)$ viewed as an operator in $L^2(M, \CC^q)$ is closable and its closure is defined on the maximal domain in $L^2(M, \CC^q)$; in particular, if $A$ is symmetric with respect to the inner product in $L^2(M, \CC^q)$, then its maximal realization is self-adjoint in $L^2(M, \CC^q)$. 
If $B \in \Op \cS^\ell_q(M)$ is another pseudodifferential operator, then by \cite[Chapter~7, Proposition~3.3 and Section~10]{T1} one has $A B \in \Op \cS^{k+\ell}_q(M)$ and its principal symbol is $\sigma_{A B}(x,\xi) = \sigma_{A}(x,\xi) \sigma_{B}(x,\xi)$. If $C \in \Op \cS^k_1(M)$, i.e. if $C$ is a scalar operator, then $(C\otimes I_q)B - B(C\otimes I_q) \in \Op \cS^{k+\ell-1}_q(M)$, see \cite[Chapter~7, equation~(3.24) and Section~10]{T1}.

Important examples of classical pseudodifferential operators are multiplication operators associated with smooth functions, which belong to $\Op \cS^0_q(M)$, and integral operators with smooth integral kernels, which are elements in $\Op \cS^{-\infty}_q(M):=\bigcap_{\ell\in\RR}\Op \cS^{\ell}_q(M)$.

Finally, we present the proof of Lemma~\ref{lemma_expansion1}, which will lead to the fact that certain integral operators can be viewed as pseudodifferential operators.

\begin{proof}[Proof of Lemma~\ref{lemma_expansion1}]
	Using the expression for $E_{z}$ from \eqref{fundamentla_solution_Dirac_2D} and \eqref{fundamentla_solution_Dirac_3D}, one sees that
	\begin{align*}
		E_{z}(x)&=\begin{pmatrix}
			(z+m) \cE_{z^2-m^2}(x)\otimes I_\nhalf & w^+_{z}(x)\\
			w^-_{z}(x) & (z-m)\cE_{z^2-m^2}(x)\otimes I_\nhalf
		\end{pmatrix},\quad x \in \mathbb{R}^n \setminus \{ 0 \},
	\end{align*}
	where $\cE_{z^2-m^2}$ is given by~\eqref{fundamental_solution_Laplace} and in the case $n=2$ we have set
	\[
	w^\pm_{z}(x)=\begin{cases}
		\dfrac{1}{2 \pi} k(z)|x| K_1 \big( - \rmi k(z) | x |\big)  \dfrac{1}{ x_1\pm \rmi x_2}, &z^2\ne m^2, \\[\bigskipamount]
		\dfrac{\rmi}{2 \pi}  \dfrac{1}{x_1\pm\rmi x_2}, & z^2=m^2,	
	\end{cases}
	\]
	whereas in the case $n=3$ we have set
	\[
	w^\pm_{z}(x)= \big( 1 - \rmi k(z) |x|\big) \dfrac{\rmi  (\sigma \cdot x)}{  4\pi | x |^3}  e^{\rmi k(z) |x|}.
	\]
	
	In the remaining part of the proof we distinguish between the two cases $n=2$ and $n=3$. 
	Let us first consider $n=2$.  Using the representations
	\begin{align} \label{Bessel_functions}
		K_1(t)&=\dfrac{1}{t}+tg_1(t^2)\log t+tg_2(t^2),&
		K_0(t)&=-\log t+g_3(t^2)+t^2g_4(t^2)\log t
	\end{align}
	with some analytic functions $g_1, g_2, g_3, g_4$, see the proof of \cite[Lemma~3.3]{BHOP20}, one has
	\begin{equation}
		\label{gmzxy2D}
		E_{z}(x)=\begin{pmatrix} -\dfrac{z+m}{2\pi}\log|x| & w^+(x) \\ w^-(x) & -\dfrac{z-m}{2\pi}\log|x| \end{pmatrix}+F_{z}(x),\quad x \in \mathbb{R}^n \setminus \{ 0 \},
	\end{equation}
	where $F_{z}$ is a matrix-valued function with components 
	\begin{equation*}
		\begin{aligned}
			F_{z}[1,1](x) &=\dfrac{z+m}{2\pi}\Big[ -\log (-\rmi k(z))
			+g_3(-k(z)^2 | x |^2)\\
			&\qquad -k(z)^2 | x |^2 g_4(-k(z)^2 | x |^2)\log (-\rmi k(z) |x|)\Big], \\
			F_{z}[1,2](x) &= -\frac{\rmi k(z)^2}{2 \pi} (x_1 - \rmi x_2) \Bigl[ g_1(-k(z)^2 | x |^2)\log (-\rmi k(z) |x|) + g_2(-k(z)^2 | x |^2) \Bigr],\\
			F_{z}[2,1](x) &= -\frac{\rmi k(z)^2}{2 \pi} (x_1 + \rmi x_2) \Bigl[ g_1(-k(z)^2 | x |^2)\log (-\rmi k(z) |x|) + g_2(-k(z)^2 | x |^2) \Bigr], \\
			F_{z}[2,2](x) &= \dfrac{z-m}{2\pi}\Big[-\log (-\rmi k(z))+g_3(-k(z)^2 | x |^2)\\
			&\qquad-k(z)^2 | x |^2 g_4(-k(z)^2 | x |^2)\log (-\rmi k(z) |x|)\Big],
		\end{aligned}
	\end{equation*}
	if $z^2 \neq m^2$, and $F_{z} = 0$, if $z^2 = m^2$. 
	Since $g_1, g_2, g_3, g_4$ are analytic, $F_z$ is smooth outside the origin and it has no singularity at $0$ (note that $|x|^2\log|x| \to 0$ and $(x_1\pm \rmi x_2)\log|x|\to 0$ as $x\to 0$). In particular, $\mathcal{K}_{z} = \chi F_z \in L^1(\mathbb{R}^n) \cap C(\mathbb{R}^n)$, which shows (a) and the first claim in (d).  Moreover, as $\nabla\log|x|=x/|x|^2$, the explicit form of $F_z$ shows 
	\begin{equation*}
	  |\nabla F_z(x)|\leq c(1+|\log|x||)
	\end{equation*}
	for $0<|x|\leq R$, and the same bound holds for $\nabla \mathcal K_z=(\nabla\chi)F_z+\chi\nabla F_z$; this finishes the proof of (d).  Moreover, as $\chi$ is compactly supported, $\mathcal{K}_{z}$ is rapidly decreasing (with all of its derivatives) at infinity, i.e. (b) is true. 
	It remains to show that $\mathcal{K}_z$ satisfies~\eqref{expansion_K_z}; here and below, smooth contributions to $\mathcal K_z$ can always be absorbed into the $C^L$-remainder in the expansion. We show that the term $\chi F_z[1,2]$ has such an expansion, the other terms can be treated in a similar way.  Note first that 
	\begin{equation*}
		f(x) := -\frac{\rmi k(z)^2}{2 \pi} (x_1 - \rmi x_2) \Bigl[ g_1(-k(z)^2 | x |^2)\log (-\rmi k(z)) + g_2(-k(z)^2 | x |^2) \Bigr]
	\end{equation*}
	is smooth. Hence, using the power series representation $g_1(t) = \sum_{j \geq 0} a_j t^j$ of the analytic function $g_1$, one finds that
	\begin{equation} \label{expansion_F12}
		\begin{split}
			F_z[1,2](x) &= -\frac{\rmi k(z)^2}{2 \pi} \sum_{j \geq 0} a_j \bigl( -k(z)^2 \bigr)^{j} (x_1 - \rmi x_2) |x|^{2j} \log |x|  + f(x) \\
			&= \sum_{j \geq 0} \widetilde{p}_{2j+1}(x) \log |x|  + f(x),
		\end{split}
	\end{equation}
	where each
	\begin{equation*}
		\widetilde{p}_{2j+1}(x) := \frac{\rmi}{2 \pi} a_j \bigl( -k(z)^2 \bigr)^{j+1} (x_1 - \rmi x_2) |x|^{2j}, \qquad j \in \mathbb{N}_0,
	\end{equation*}
	is a homogeneous polynomial of order $2j+1$. To conclude the asymptotic expansion in~\eqref{expansion_K_z}, fix $L \in \mathbb{N}$ and note that one has the Taylor series expansion
	\begin{equation*}
		\chi(x) = \sum_{|\gamma| \leq L} c_\gamma x^\gamma + \sum_{|\gamma|=L+1} x^\gamma \widetilde{\chi}_\gamma(x),
	\end{equation*}
	where the sums are over multi-indices $\gamma \in \mathbb{N}_0^2$, $c_\gamma$ are suitable constants, and $\widetilde{\chi}_\gamma$ is smooth. Hence, we can write
	\begin{equation*}
		\begin{split}
			\chi(x) F_z[1,2](x) &= \left( \sum_{|\gamma| \leq L} c_\gamma x^\gamma + \sum_{|\gamma|=L+1} x^\gamma \widetilde{\chi}_\gamma(x) \right) \cdot \left(  \sum_{j \geq 0} \widetilde{p}_{2j+1}(x) \log |x| \right) + \chi(x) f(x)  \\
			&= \sum_{\ell=0}^L p_\ell(x) \log|x| +  \sum_{\substack{|\gamma| \leq L \\ |\gamma| + 2 j + 1 > L+1}} c_\gamma x^{\gamma} \widetilde{p}_{2j+1}(x)  \log |x|  \\
			&\qquad \qquad \qquad +  \sum_{|\gamma| = L + 1} x^\gamma \widetilde{\chi}_\gamma(x) \log |x| \sum_{j \geq 0} \widetilde{p}_{2j+1}(x) + \chi(x) f(x)
		\end{split}
	\end{equation*}
	with 
	\begin{equation*}
		p_\ell(x) = \sum_{|\gamma| + 2 j + 1 = \ell+1} c_\gamma x^\gamma \widetilde{p}_{2 j + 1}(x), \qquad \ell \in \{ 0, \dots, L \},
	\end{equation*}
	being a homogeneous polynomial of order $\ell+1$. Since $x^\gamma \log|x| \in C^L(\mathbb{R}^2)$ for $|\gamma| > L$, one concludes that
	\begin{equation*}
		\chi(x) F_z[1,2](x) - \sum_{\ell=0}^L p_\ell(x) \log|x| \in C^L(\mathbb{R}^2),
	\end{equation*}
	which yields~\eqref{expansion_K_z} with $q_\ell \equiv 0$ for all $\ell$;  note that the smooth function $\chi f$ contributes only to the $C^L$-remainders. Thus, also (c) is true for $n=2$.

	In a similar way, one has for $n=3$ 
	\begin{equation} 
		\label{gmzxy3D}
			E_{z}(x) =\begin{pmatrix}
			\dfrac{z+m}{4\pi|x|} \otimes I_2  & w^+(x)\\
			w^-(x) & \dfrac{z-m}{4\pi|x|} \otimes I_2
		\end{pmatrix}+F_{z}(x),\quad x \in \mathbb{R}^n \setminus \{ 0 \},
	\end{equation}
	where $F_{z}$ is a block matrix function
	\begin{equation*}
		F_{z}(x) = \begin{pmatrix} F_{z}[1,1](x) & F_{z}[1,2](x) \\ F_{z}[2,1](x) & F_{z}[2,2](x) \end{pmatrix}
	\end{equation*}
	with blocks $F_{z}[j,k] \in \mathbb{C}^{2 \times 2}$ given by
	\begin{equation*}
		\begin{split}
			F_{z}[1,1](x) &=
			\dfrac{z+m}{4\pi|x|} \bigl( e^{\rmi k(z) |x|} - 1 \bigr) \otimes I_2, \\
			F_{z}[1,2](x) &= F_{z}[2,1](x) = \dfrac{\rmi  (\sigma \cdot x)}{  4\pi | x |^3} \left( \big( 1 - \rmi k(z) |x|\big)  e^{\rmi k(z) |x|} - 1 \right),\\
			F_{z}[2,2](x) &= \dfrac{z-m} {4\pi|x|} \bigl( e^{\rmi k(z) |x|} - 1 \bigr) \otimes I_2.
		\end{split}
	\end{equation*}
	Taking the power series representation of the exponential function into account,
	one finds first that 
	\begin{equation*}
		F_{z}[1,1](x) =
		\dfrac{z+m}{4\pi|x|} \bigl( e^{\rmi k(z) |x|} - 1 \bigr) \otimes I_2 
		= \dfrac{z+m}{4\pi} \left( \sum_{\ell=0}^\infty \frac{(\rmi k(z))^{\ell+1}}{(\ell+1)!} |x|^\ell\right) \otimes I_2.
	\end{equation*}
	Note that each term in the above sum is a homogeneous function (and hence, a homogeneous distribution) of degree $\ell$ that is smooth outside the origin.
	A similar expansion holds for $F_{z}[2,2]$. Using again the power series representation of the exponential function, one obtains 
	\begin{equation*} 
		F_{z}[1,2](x)=\sum_{\ell=0}^\infty\rmi^{\ell+3}\left(\frac{1}{(\ell+2)!}-\frac{1}{(\ell+1)!}\right)k(z)^{\ell+2} |x|^{\ell-1} \frac{ (\sigma\cdot x) }{4\pi}.
	\end{equation*}
	Again, each term in the above sum is a homogeneous function (and hence, a homogeneous distribution) of degree $\ell$ that is smooth outside the origin.  In particular, the term with $\ell=0$ in the last sum equals $\frac{\rmi k(z)^2}{8\pi}\frac{\sigma\cdot x}{|x|}$, which is bounded, homogeneous of degree $0$, discontinuous at the origin if $k(z)\neq0$, and whose gradient is bounded by $c|x|^{-1}$. All other terms, i.e., $|x|^{\ell}$ for $\ell\geq 0$ in the diagonal blocks and $|x|^{\ell-1}(\sigma\cdot x)$ for $\ell\geq 1$ in the off-diagonal blocks, are continuous on $\RR^3$ with gradients bounded by $c\ell R^{\ell-1}$ on $B_R(0)\setminus\{0\}$, and their coefficients decay factorially; hence these series and the series of their gradients converge uniformly on $B_R(0)\setminus\{0\}$. Consequently, $F_z$ is bounded near $0$, $\mathcal K_z=\chi F_z\in L^1(\RR^n)\cap L^\infty(\RR^n)$, and $|\nabla \mathcal K_z(x)|\leq c|x|^{-1}$ for $0<|x|\leq R$.  This function is smooth outside the origin and rapidly decreasing (with all of its derivatives) at infinity. This shows (a), (b), and (e). Moreover, using the Taylor series expansion of $\chi$ and the last two displayed formulas, one finds via a Cauchy product that 
	$\mathcal{K}_z$ satisfies~\eqref{expansion_K_z} (with $p_\ell \equiv 0$ for all $\ell$), i.e. (c) is true.  
\end{proof}

\section{A brief introduction to boundary triples and Weyl functions}\label{sec-bt}

Let us recall some basic facts about boundary triples in a minimalistic form strictly adapted to the needs of the present paper. We follow \cite[Section 1]{BGP},  where the proofs of all statements can be found.
We also refer the interested reader to the original papers \cite{DM91,DM95} and the monographs \cite{BHS-book,DM-book} for a more detailed discussion.

\subsection{General overview}

Let $S$ be a closed densely defined symmetric operator in a Hilbert space $\cH$.

\begin{definition}
	A \emph{boundary triple} for $S^*$ is a triple $(\cG,\Pi_0, \Pi_1)$
	consisting of a Hilbert space $\cG$ and two linear maps $\Pi_0,\Pi_1: \dom S^* \rightarrow \cG$
	satisfying the following two conditions:
	\begin{itemize}
		\item[\textup{(i)}] For all $f, g \in \dom S^*$ it holds that
		\begin{equation*}
			\langle S^* f, g \rangle_\cH - \langle f, S^* g \rangle_\cH = \langle\Pi_1 f, \Pi_0 g\rangle_\cG - \langle\Pi_0 f, \Pi_1 g\rangle_\cG.
		\end{equation*}
		\item[\textup{(ii)}] The map $\dom S^*\ni f\mapsto (\Pi_0 f, \Pi_1 f )\in \cG\times\cG$ is surjective. 
	\end{itemize}
\end{definition}
A boundary triple for $S^*$ exists if and only if  $S$ admits self-adjoint extensions in $\cH$. From now on we assume that this is satisfied, and we pick a boundary triple $(\cG,\Pi_0,\Pi_1)$. Then a number of additional associated objects arise naturally. First, the operator
\[
Q := S^* |_{\ker \Pi_0}
\]
turns out to be self-adjoint, and for any $z\in \res Q$ one has the direct sum decomposition
\begin{equation}\label{deco33}
	\dom S^* = \dom Q \dsum \ker (S^* - z) = \ker \Pi_0 \dsum \ker (S^* - z),
\end{equation}
showing that the restriction
$\Pi_0: \ker (S^*-z)\to \cG$
is bijective. The $\gamma$-field $z\mapsto \Gamma_z$ and the Weyl function $z\mapsto M_z$
associated with the chosen boundary triple are the holomorphic maps defined by
\begin{align*}
	\res Q \ni z &\mapsto \Gamma_z := \big[\Pi_0|_{\ker (S^*-z)}\big]^{-1}\in\cL(\cG,\cH),\\
	\res Q \ni z &\mapsto M_z :=\Pi_1 \,\Gamma_z\in \cL(\cG).
\end{align*}
In addition, one may consider $\Gamma_z$ as a bounded linear operator from $\cG$ to $\dom S^*$ endowed with the graph norm, because $\ran \Gamma_z\subset \ker (S^*-z)$ and the graph norm on $\ker(S^*-z)$ is equivalent to the norm of $\cH$. The adjoints of $\Gamma_z$ and $M_z$ are given by
\begin{equation*} 
	\Gamma_z^* = \Pi_1 (Q - \overline z)^{-1}, \qquad \quad M_z^* = M_{\Bar{z}}.
\end{equation*}

The following theorem describes the main application of boundary triples.

\begin{theorem}\label{thm12}
For any linear operator $\Theta$ in $\cG$, denote by $Q_\Theta$ the restriction of $S^*$ on the domain
\begin{equation*} 
	\label{qtheta}
	\dom Q_\Theta = \big\{ f \in \dom S^*: \Pi_0f\in\dom\Theta,\ \Pi_1 f = \Theta \Pi_0 f\big\}.
\end{equation*}
Then, the following is true:
\begin{itemize}
	\item[\textup{(a)}] The operator $Q_\Theta$ is self-adjoint in $\cH$ if and only if the operator $\Theta$ is self-adjoint in $\cG$. 
\end{itemize}	
If this assumption is fulfilled, then for any  $z \in \res Q$ the following assertions are true:
	\begin{itemize}
		\item[\textup{(b)}] $z\in \spec Q_\Theta$ if and only if $0 \in \spec(\Theta-M_z)$.
		\item[\textup{(c)}] $\ker (Q_\Theta-z) = \Gamma_z  \ker  (\Theta-M_z)$. In particular,
		 $z\in \specp Q_\Theta$ if and only if $0 \in \specp(\Theta-M_z)$ with the same multiplicities.
		\item[\textup{(d)}] $z\in \spece Q_\Theta$ if and only if $0 \in \spece(\Theta-M_z)$.
	\end{itemize}
Moreover, for  all $z \in \res Q_\Theta \cap \res Q$ one has the resolvent formula
\[
	(Q_\Theta-z)^{-1}=(Q-z)^{-1}+\Gamma_z(\Theta-M_z)^{-1} \Gamma_{\Bar{z}}^*.
\]
\end{theorem}

It seems that the following simple result is not explicitly stated anywhere in the required form, so we provide a short complete proof:
\begin{lemma}\label{lemma12}
	If $\Theta_j,\Theta$ are self-adjoint operators in $\cG$ such that $\Theta_j\rightarrow\Theta$, $j\rightarrow\infty$, in the norm resolvent sense, then
	$Q_{\Theta_j}\rightarrow Q_\Theta$, $j\rightarrow\infty$, in the norm resolvent sense.
\end{lemma}

\begin{proof}
Let us fix some $z\in\CC\setminus\RR$. In view of the resolvent formula in Theorem~\ref{thm12} it is sufficient to show that $(\Theta_j-M_z)^{-1}$ converges in the norm to $(\Theta-M_z)^{-1}$ as $j\rightarrow\infty$.
Recall that $\Theta_j-M_z$ and $\Theta-M_z$ have bounded inverses by Theorem~\ref{thm12}.

Note first that $(\Theta-M_z)^{-1}-(\Theta+\rmi)^{-1}=(\Theta-M_z)^{-1}(M_z+\rmi)(\Theta+\rmi)^{-1}$, and
hence
\[
(\Theta-M_z)^{-1}\big( I-(M_z+\rmi)(\Theta+\rmi)^{-1}\big)=(\Theta+\rmi)^{-1}.
\]
It follows from this that the operator $T:=I-(M_z+\rmi)(\Theta+\rmi)^{-1}$ has an everywhere defined bounded inverse, and
\[
	(\Theta-M_z)^{-1}=(\Theta+\rmi)^{-1}T^{-1}.
\]
The same conclusion and identity hold if one replaces the operators $\Theta$ and $T$ by $\Theta_j$ and  $T_j:=I-(M_z+\rmi)(\Theta_j+\rmi)^{-1}$, respectively. Our assumption $\Theta_j\rightarrow\Theta$, $j\rightarrow\infty$, in the norm resolvent sense
implies $T_j\rightarrow T$, $j\rightarrow\infty$, in the operator norm. As both $T$ and $T_j$ have bounded inverses, one also
has
\[
T_j^{-1}=\big(T+(T_j-T)\big)^{-1}=\big(I+T^{-1}(T_j-T)\big)^{-1}T^{-1}\rightarrow T^{-1},\quad j\rightarrow\infty,
\]
in the operator norm. Now it follows that
\[
(\Theta_j-M_z)^{-1}=(\Theta_j+\rmi)^{-1}T_j^{-1}\rightarrow (\Theta+\rmi)^{-1}T^{-1}=(\Theta-M_z)^{-1},\quad j\rightarrow\infty,
\]
which yields the claim.
\end{proof}

For further aspects on the dependence of $Q_\Theta$ on $\Theta$ see, e.g. \cite[Section 2.8]{BHS-book}, \cite[Section~8.5.2]{DM-book}, and \cite{lat}.

\subsection{Construction of boundary triples using abstract trace maps} Finally, we recall an approach to the construction of boundary triples using abstract trace maps as initially proposed in~\cite{Pos04,Pos08}, see also \cite[Section~1.4.2]{BGP}.

Let $Q$ be a self-adjoint operator in the Hilbert space $\cH$. Let $\cG$ be another Hilbert space and
\begin{equation*}
	\gamma: \dom Q \rightarrow \cG
\end{equation*}
be a surjective linear operator which is bounded with respect to the graph norm of $Q$
and such that $\ker \gamma$ is a dense subspace of $\cH$ (such $\gamma$ are often referred to as abstract trace maps). Then
\begin{equation*}
	S := Q|_{\ker\gamma}
\end{equation*}
is a closed densely defined symmetric operator in $\cH$. For  $z\in\res Q$ define
\begin{equation}
	\label{eq-gz1a}
	\Gamma_z:= \big(\gamma(Q - \Bar{z})^{-1}\big)^*\in \cL(\cG,\cH),
\end{equation}
then
\[
\ran \Gamma_z=\ker(S^*-z), \qquad \dom S^* = \dom Q \dsum \ran \Gamma_z,
\]
i.e. each $f \in \dom S^*$ is uniquely represented
as
\begin{equation}\label{jabitte2}
f = f_z + \Gamma_z \xi,\quad f_z\in\dom Q,\quad \xi\in\cG,
\end{equation}
and one can show that $\xi$ in the above representation is independent of $z$.

\begin{theorem} \label{thm13}
	Pick an arbitrary $\zeta \in \res Q$ and define $\Pi_0, \Pi_1: \dom S^* \rightarrow \cG$
	by
	\begin{equation*}
		\Pi_0 f := \xi \quad\text{and}\quad \Pi_1 f := \frac{1}{2}\, \gamma (f_\zeta + f_{\Bar{\zeta}}).
	\end{equation*}
	Then $(\cG,\Pi_0,\Pi_1)$ is a boundary triple for $S^*$ such that $Q=S^*|_{\ker\Pi_0}$. Moreover, the associated $\gamma$-field $\Gamma_z$ is given by~\eqref{eq-gz1a} and the associated Weyl function $M_z$
	is 
	\begin{equation}
		\label{eqmz22}
	M_z = \gamma\Big(\Gamma_z - \frac{1}{2} (\Gamma_\zeta + \Gamma_{\Bar{\zeta}}) \Big).
	\end{equation}
\end{theorem}

\begin{remark}\label{rem-qq}
	If the operator $Q$ has purely discrete spectrum, then its resolvent $z\mapsto (Q-z)^{-1}$
	is meromorphic in $\CC$, and the only possible singularities are poles of first order at the eigenvalues of $Q$ with finite-dimensional residues equal  to the negatives of the orthogonal projectors on the respective
	eigensubspaces. From the construction of the maps $\Gamma_z$ and and $M_z$ in \eqref{eq-gz1a} and \eqref{eqmz22} it follows
	that the maps
	\[
	z\mapsto \Gamma_z \in \cL(\cG,\cH), \quad z\mapsto \Gamma_{\Bar z}^* \in \cL(\cH,\cG),
	\quad z\mapsto M_z\in\cL(\cG)
	\]
	are also meromorphic in $\CC$ and only have poles of first order with finite-dimensional residues.
\end{remark}

\end{appendix}


\begin{thebibliography}{999}
		
		
		\normalsize
		
		\itemsep 3pt
		
	\bibitem{ammann}
	 B. Ammann, E. Humbert, B. Morel, \emph{Mass endomorphism and spinorial Yamabe type problems on conformally flat manifolds.} Comm. Anal. Geom. {\bf 14} (2006) 163--182.
		
	\bibitem{AMV14} N. Arrizabalaga, A. Mas, L. Vega, {\em Shell interactions for Dirac operators}. J. Math. Pures Appl. {\bf 102} (2014) 617--639. 
     \bibitem{AMV15} N. Arrizabalaga, A. Mas, L. Vega, {\em Shell interactions for Dirac operators: on the point spectrum and the confinement.} SIAM J. Math. Anal. {\bf 47} (2015) 1044–-1069.

	\bibitem{baer}
C. B\"ar, W. Ballmann, \emph{Guide to boundary value problems for Dirac-type operators.}
W.~Ballmann et al. (Eds.): \emph{Arbeitstagung Bonn 2013}. Progress in Mathematics, Vol.~319, Birkh\"auser, 2016, pp.~~43--80.




\bibitem{baer2}
C. B\"ar, W. Ballmann, \emph{Boundary value problems for elliptic differential operators of first order.}
H.-D.~Cao, S.-T.~Yau (Eds.): \emph{Algebra and Geometry: In Memory of C.~C.~Hsiung.} 
Surveys in Differential Geometry, Vol.~17, International Press, 2012, pp.~1--78.
	


\bibitem{BEHL18}
J.~Behrndt, P.~Exner, M. Holzmann,  V.~Lotoreichik,
\emph{On the spectral properties of {D}irac operators with electrostatic
  $\delta$-shell interactions.}
 J. Math. Pures Appl. {\bf 111} (2018) 47--78.


\bibitem{BEHL19} J. Behrndt, P. Exner, M. Holzmann, V. Lotoreichik, {\em On Dirac operators in $\RR^3$ with electrostatic and Lorentz scalar $\delta$-shell interactions.} Quantum Stud. Math. Found. {\bf 6} (2019) 295--314.


\bibitem{BEHT25} J. Behrndt, P. Exner, M. Holzmann,  M. Tusek, 
{\em On two-dimensional Dirac operators with delta-shell interactions supported on unbounded curves with straight ends.}
B.~Cassano et al. (Eds), \emph{Singularities, asymptotics, and limiting models.} Springer INdAM Series, Vol.~64, Springer, 2025, pp.~35--77

	

	\bibitem{BHS-book}
	J.~Behrndt, S.~Hassi, H.~S.~V.~de Snoo, {\it Boundary value problems, Weyl functions, and differential operators}.
	Monographs in Mathematics, Vol.~108, Birkh\"auser, 2020. 
	
		\bibitem{BH20} J. Behrndt, M. Holzmann, {\em On Dirac operators with electrostatic delta-shell interactions of critical strength.} J. Spectral Theory {\bf 10} (2020) 147--184.
	
	\bibitem{BHM}
	J. Behrndt, M. Holzmann, A. Mas, \emph{Self-adjoint Dirac operators on domains in $\mathbb{R}^3$.}
	Ann. Henri Poincar\'e {\bf 21} (2020) 2681--2735.
	
	
	\bibitem{BHSLS24}
	J.~Behrndt, M.~Holzmann, C.~Stelzer-Landauer, G.~Stenzel,
	\emph{Boundary triples and Weyl functions for Dirac operators with singular interactions.}
	Rev. Math. Phys. {\bf 36} (2024) 2350036.
%
	\bibitem{BHOP20}
	J.~Behrndt, M.~Holzmann, T.~Ourmi\`eres-Bonafos, K.~Pankrashkin,
	\emph{Two-dimensional Dirac operators with singular interactions supported on closed curves}.
	J. Funct. Anal. {\bf 279} (2020) 108700.

\bibitem{bell} S. R. Bell, \emph{The Cauchy transform, potential theory and conformal mapping.}
CRC Press, 2016.

	\bibitem{BFSV}
	R.~D.~Benguria, S.~Fournais, E.~Stockmeyer, H.~Van~Den~Bosch,
	\emph{Self-adjointness of two-dimensional Dirac operators on domains.}
	Ann. Henri Poincar\'e {\bf 18} (2017) 1371--1383.




\bibitem{BB222} B. Benhellal, {\em Spectral properties of the Dirac operator coupled with $\delta$-shell interactions}.  Lett. Math. Phys. {\bf 112} (2022) 52.
%
\bibitem{BBP} B. Benhellal, V. Bruneau, P. Miranda, {\em The Landau-Dirac operator with shell interactions: self-adjointness and clustering}. Preprint arXiv:2608.20665 (2026).

\bibitem{BBZ} B. Benhellal, V. Bruneau, M. Zreik, {\em A Poincaré-Steklov map for the MIT bag model.} Anal. PDE
{\bf 17} (2024) 2923--2970.

\bibitem{BKMP26} B. Benhellal, N. K\"orner, D. Machado, Konstantin Pankrashkin, {\em MIT bag model and infinite mass limit in non-smooth domains}. Preprint arXiv:2609.08791 (2026).

\bibitem{BP24} B. Benhellal, K. Pankrashkin, {\em Curvature contribution to the essential spectrum of Dirac operators with critical shell interactions}. Pure Appl. Anal. {\bf 6} (2024) 237--252.

\bibitem{BPZ25} B. Benhellal, K. Pankrashkin, M. Zreik, {\em On the self-adjointness of two-dimensional relativistic shell interactions.} J. Operator Theory {\bf 93} (2025) 569--591.

\bibitem{BCF26} W. Borrelli, P. Carimati, D. Fermi,
{\em On two-dimensional Dirac operators with critical delta-shell interactions.}
J. Math. Anal. Appl. {\bf 564} (2026) 130880.


\bibitem{BGP} J.~Br\"uning, V.~Geyler, K.~Pankrashkin, {\it Spectra of self-adjoint extensions and applications to solvable Schr\"odinger operators}. 	Rev. Math. Phys. {\bf 20} (2008) 1--70.

\bibitem{CLMT21}
B.~Cassano, V.~Lotoreichik, A.~Mas, M.~Tu\v{s}ek,
{\it General $\delta$-shell interactions for the two-dimensional Dirac operator: self-adjointness and approximation.}
Rev. Mat. Iberoam. {\bf 39} (2023) 1443--1492.


\bibitem{DM91}
V.~A.~Derkach, M.~M.~Malamud,
\emph{Generalized resolvents and the boundary value problems for
	Hermitian operators with gaps}. J. Funct. Anal. {\bf 95} (1991) 1--95.

\bibitem{DM95}
V.~A.~Derkach, M.~M.~Malamud,
{\it The extension theory of Hermitian operators and the moment problem}.
J. Math. Sci. (N. Y.) {\bf 73} (1995) 141--242.


\bibitem{DM-book}
V.~Derkach, M.~Malamud,
\emph{Extension theory of symmetric operators and boundary value problems.}
Operator Theory: Advances and Applications, Vol.~309, Birkh\"auser, 2025.


\bibitem{DES89}
J. Dittrich, P. Exner, P. \v{S}eba, {\it Dirac operators with a spherically symmetric $\delta$-shell
interaction.} J. Math. Phys. {\bf  30} (1989) 2875--2882.



%
%
%
%


\bibitem{grosse} N.~Grosse, A.~Uribe, H.~Van Den Bosch, \emph{Local boundary conditions for Dirac-type operators.}  Ann. Henri Poincar\'e {\bf 27} (2026) 3635--3672.

\bibitem{HW} G. C. Hsiao, W. L. Wendland, {\em Boundary integral equations. Second edition.} 
Springer, 2021.

\bibitem{lat} Y. Latushkin, S. Sukhtaiev, \emph{First-order asymptotic perturbation theory for extensions of symmetric operators.} J. Lond. Math. Soc. {\bf 110} (2024) e13005.

\bibitem{M00}
W.~McLean,
\newblock {\em Strongly Elliptic Systems and Boundary Integral Equations}.
\newblock Cambridge University Press, Cambridge, 2000.

\bibitem{musk} N. I. Muskheshvili, \emph{Singular integral equations}. Wolter-Noordogg Publishing, 1958.


%
%
%
%
%
%

\bibitem{DLMF}
	F.\,W.\,J.~Olver, A.\,B.~{Olde Daalhuis}, D.\,W.~Lozier, B.\,I.~Schneider, R.\,F.~Boisvert, C.\,W.~Clark, B.\,R.~Miller, B.\,V.~Saunders, H.\,S.~Cohl, M.\,A.~McClain (Eds.),
	\textit{NIST Digital Library of Mathematical Functions.} Available online at 	\url{https://dlmf.nist.gov}.

	\bibitem{OP21}
	T.~Ourmi{\`e}res-Bonafos, F. Pizzichillo, {\it Dirac operators and shell interactions: a survey.} A. Michelangeli (Ed.), {\it Mathematical challenges of zero-range physics.} Springer INDAM Series, Vol.~42, Springer, 2021, pp.~105–-131.
	
	
	\bibitem{OV}
	T.~Ourmi{\`e}res-Bonafos, L.~Vega,
	{\em A strategy for self-adjointness of {D}irac operators: application to the {MIT} bag model and $\delta$-shell interactions}. 
	Publ. Mat. {\bf 62} (2018) 397--437.


\bibitem{KPconvex} K. Pankrashkin, \emph{MIT bag in non-smooth convex domains.} Rev. Mat. Iberoam.
{\bf 42} (2026) 2377--2402.

%
%
	\bibitem{Pos04} A. Posilicano, \emph{Boundary triples and Weyl functions for
		singular perturbations of self-adjoint operators.}
	Meth. Funct. Anal. Topol. {\bf 10} (2004)  57--63. 
	
	\bibitem{Pos08} A. Posilicano, \emph{Self-adjoint extensions of
		restrictions.} Oper. Matr. {\bf 2} (2008)  483--506.
	

	
	

	
	
\bibitem{R21-1}
V. Rabinovich, {\em Two-dimensional Dirac operators with interactions on unbounded smooth curves.} Russ. J. Math. Phys. {\bf 28} (2021) 524--542.


\bibitem{R21-2}
V. Rabinovich, {\em Dirac operators with singular potentials supported on unbounded surfaces in $\mathbb R^3$.}
Funct. Anal. Appl. {\bf 55} (2021) 245--249.
	

\bibitem{R23}
V. Rabinovich, 
{\em Method of potential operators for interaction problems on unbounded hypersurfaces in $\mathbb R^n$ for Dirac operators.}
Russ. J. Math. Phys. {\bf 30} (2023) 674--690.
	
		\bibitem{raulot} S. Raulot, \emph{Green functions for the Dirac operator under local boundary conditions and applications.}
	Ann. Glob. Anal. Geom. {\bf 39} (2011) 337--359.
	%
	\bibitem{RS4} M. Reed, B. Simon, {\em Methods of modern mathematical physics. IV: Analysis of operators.}
	Academic Press, New York, 1978. 
	
	\bibitem{Seeley} R. T. Seeley, {\em Complex powers of an elliptic eperator},   Singular integrals (Proc. Sympos. Pure Math., Chicago, IL, 1966), American Mathematical Society, Providence, R.I., 1967, pp.~288–-307.
	
	\bibitem{stein} O. Steinbach, {\em Numerical approximation methods for elliptic boundary value problems.} Springer, 2008.
%
%
%




	\bibitem{taylor1} M. Taylor, \emph{Pseudodifferential operators}. Princeton University Press, 1981.


	\bibitem{T1} M. Taylor, {\em Partial differential equations. II. Qualitative studies of linear equations}.  Volume 116 of Applied Mathematical Sciences. Springer-Verlag, New York, 1996.
%
%
%
	\bibitem{tretter} C. Tretter, {\it Spectral theory of block operator matrices and applications}. Imperial College Press, 2008.
	
		\bibitem{weid} J. Weidmann, \emph{Lineare Operatoren in Hilberträumen. Teil 1. Grundlagen.}
	Teubner, 2000.
	
	\bibitem{zworski} M. Zworski, {\it Semiclassical analysis}. Graduate Studies in Mathematics, Vol.~138, AMS, 2012.
\end{thebibliography}
	\end{document}